\documentclass[aip,reprint]{revtex4-2} 

\usepackage{amssymb,amsmath}
\usepackage{array,multirow,graphicx}
\usepackage{enumitem}
\usepackage{multirow}
\usepackage{adjustbox}
\usepackage{bbm}
\usepackage{makecell}
\usepackage{tikz}
\usepackage{graphicx}
\usepackage{tabularx}
\newcolumntype{Y}{>{\centering\arraybackslash}X}
\usepackage{forest}
\usepackage{tcolorbox}
\usepackage{makecell}
\newcolumntype{s}{>{\centering\arraybackslash}X} 
\usepackage{tikz}           
\usetikzlibrary{
    arrows.meta,            
    positioning,            
    calc,                   
    fit,                    
    backgrounds             
}
\usepackage[utf8]{inputenc}
\usepackage[T1]{fontenc}
\usepackage{hyperref}
\usepackage{dirtytalk}
\usepackage{mathtools}
\usepackage[english]{babel}
\usepackage{amsthm}
\usepackage{float}
\usepackage{amsfonts,mathrsfs}
\usepackage{bm}
\usepackage{mathrsfs}
\usepackage{makecell}
\usepackage{ragged2e}

\usepackage[table]{xcolor}
\usepackage{tabularx, booktabs}

\usepackage{array}
\usepackage{colortbl}
\tcbuselibrary{skins}
\usepackage{tikz}
\usetikzlibrary{matrix}

\usepackage{epstopdf}
\usepackage{subcaption}
\usepackage{booktabs}
\usepackage{tablefootnote}
\graphicspath{ {images/} }
\usepackage{comment}

\theoremstyle{definition}

\def\R{\mathbb{R}}

\def\L{\mathcal{L}}

\def\N{\mathcal{N}}

\def\P{\mathcal{P}}

\def\Pi{\mathbf{P}}

\tcbset{tab2/.style={enhanced,fonttitle=\bfseries,fontupper=\normalsize\sffamily,
colback=yellow!10!white,colframe=red!50!black,colbacktitle=Salmon!40!white,
coltitle=black,center title}}
\newcolumntype{M}{>{\centering\arraybackslash}m{1.5cm}}
\newcolumntype{P}{>{\centering\arraybackslash}m{5.5cm}}

\begin{document}

\title{RELift: Learned Coarse-to-Fine Propagators for Time-Dependent PDEs with Applications to Electron Dynamics}

\newcommand{\eqmark}{\textsuperscript{\dag}}

\author{Hardeep Bassi}
\email{hbassi2@ucmerced.edu, hbassi@lbl.gov}
\altaffiliation{Work done as a visiting student researcher at Lawrence Berkeley National Laboratory.}
\affiliation{University of California, Merced, Merced, CA, USA}

\author{Yuanran Zhu}
\author{Erika Ye}
\author{Pu Ren}
\author{Alec Dektor}

\affiliation{Lawrence Berkeley National Laboratory, Berkeley, CA, USA}
\author{Michael W. Mahoney} 
\affiliation{Lawrence Berkeley National Laboratory, Berkeley, CA, USA}
\affiliation{International Computer Science Institute, Berkeley, CA, USA}
\affiliation{University of California, Berkeley, Berkeley, CA, USA}

\author{Harish S. Bhat}
\affiliation{University of California, Merced, Merced, CA, USA}

\author{Chao Yang}
\email{cyang@lbl.gov}
\affiliation{Lawrence Berkeley National Laboratory, Berkeley, CA, USA}

\date{\today}
\begin{abstract}
We present \textbf{RELift} (\textbf{R}estrict, \textbf{E}volve, \textbf{Lift}), a two-phase learning framework that couples coarse-grid numerical solvers with neural operators to super-resolve and forecast fine-grid dynamics for time-dependent partial differential equations (PDEs).  
In Phase~1, RELift learns a super-resolution operator that maps the solution on a coarse grid to a fine grid.  
In Phase~2, this learned operator is composed with a coarse-grid numerical integrator to construct an effective fine-grid propagator for the governing equation.  
We benchmark RELift on three canonical two-dimensional PDEs of increasing dynamical complexity---the heat equation, the wave equation, and the incompressible Navier--Stokes equations---and we further demonstrate its performance on a kinetic electron dynamics case study via the 1D1V Vlasov--Poisson system.  
Across all examples, RELift delivers high-fidelity super-resolution (Phase~1) and accurate long-horizon rollouts (Phase~2), outperforming standard super-resolution and neural operator baselines in both field-level error metrics and physics-relevant diagnostics.  
Finally, we provide error analysis of the effective fine-grid propagator, characterizing how approximation errors accumulate over time and explaining the observed numerical stability of the RELift framework.

\end{abstract}


\maketitle

\section{Introduction} 
Partial differential equations (PDEs) form the mathematical backbone for a variety of problems in the physical sciences and engineering. Reliable numerical solutions, therefore, underpin progress in areas as diverse as climate forecasting~\cite{Warner2011NWP}, photonics~\cite{Oskooi2010MEEP}, and materials design~\cite{Steinbach2013PhaseField}. A persistent obstacle in being able to generate solutions to these PDEs, however, is the scale at which their dynamics evolve \cite{lei2023machine, xie2024ab}. This is closely tied to the notion of \textit{resolution}; important physical effects for a system of interest often emerge only after one refines the resolution of the computational grid to be fine enough. Performing such a refinement, however, proves to be costly even with modern day supercomputing power.

A natural response to this computational bottleneck has been the emergence of machine learning (ML) based accelerations for PDE solvers \cite{brunton2023machine, brunton2024promising}.  These can be broadly grouped into two complementary strands.  The first strand projects the full high-dimensional dynamics onto a reduced set of variables, yielding a lower-order evolution equation that often includes an unresolved closure term representing the missing subgrid-scale physics. A neural network is then trained to approximate this term, thereby retaining the structure of the governing equations while supplying the missing physics that simplified models lack \cite{girimaji2024turbulence, duraisamy2021perspectives, brunton2022data,bhat2022dynamic, bassi2024learning, zhu2025predicting}. The second strand adopts a fully equation-free, data-driven approach, opting to learn the complete propagator of the PDE and any refinements directly from training snapshots~\cite{long2018pde,cao2023machine, lu2019deeponet, li2021fourier, mi2025conservation}.  While this approach allows equation-agnostic training and inference, the lack of explicit physical constraints or governing equations can limit the fidelity and robustness of the learned propagators, especially when extrapolating the dynamics beyond the training regime \cite{neurde_TR}.

In this paper, we introduce a two-phase operator-learning framework---called \textit{Restrict, Evolve, Lift} (RELift)---for learning and predicting PDE dynamics at fine scales from coarse scale data.  
In Phase~1, we train a neural operator that maps coarse-grid PDE solutions to their fine-grid counterparts over short time windows. In Phase~2, the learned super-resolution (SR) operator is composed with a coarse-grid numerical integrator to construct an effective fine-grid propagator for the governing PDE. As such, the dynamics prediction procedure within RELift can be understood as a predictor–corrector loop: at each step, the coarse-grid numerical solver advances the dynamics, and this coarse-scale prediction is then corrected by applying the super-resolution operator to lift it back to the fine grid. Because the coarse-grid integrator already respects the governing PDE and the neural operator serves only as a corrective spatial lift, we expect the global approximation error to grow slower than in fully data-driven propagators or approaches that merely upsample predictions from the coarse-grid numerical solvers~\cite{ing2007accumulated,pasini2024continuous,chattopadhyay2023long}. In this light, RELift has features of both strands of ML-based PDE acceleration, and it occupies a middle ground between fully equation-free and strictly PDE-constrained approaches (technically different than but somewhat akin to the Neural Discrete Equilibrium (NeurDE) framework \cite{neurde_TR}).
In particular, because Phase~2 embeds the physical laws through the coarse-grid time integrator, the learned super-resolution operator remains flexible and purely data-driven, with no explicit restriction or enforcement of physical constraints at either resolution.


Furthermore, RELift is not tied to a single data generation paradigm, and it can accommodate both (i) the commonly used restriction/synthetic-downsampling (sensor-style) setting~\cite{brunton2024promising, kutz2024shallow, ren2023physr, ren2025superbench} and (ii) the simulation--to--simulation setting. In the former, coarse inputs are obtained by applying a fixed restriction (e.g., decimation or averaging) to fine solutions at each time; and the network learns to invert this artificial map, representative of imperfect measurements yielded by real-world senors. In the latter, we generate coarse trajectories using an actual coarse-grid time integrator and fine trajectories using a fine-grid time integrator; this more closely matches deployment when only a coarse solver is available at inference. Because the coarse solver introduces numerical effects (e.g., truncation, dispersion, and aliasing), the coarse trajectory generally does not coincide with a pointwise downsampling of the fine solution, making simulation--to--simulation super-resolution a stringent regime for both reconstruction and time forecasting. As a result, RELift is especially valuable when fine solutions are unavailable or too expensive to compute or store beyond those generated for training. Conceptually, RELift also shares the same restrict–evolve–lift paradigm as the equation-free multiscale computation (EFMC)~\cite{kevrekidis2009equation,samaey2006coarse}, and the SR portion can be understood as a data-driven coarse-grid correction technique within the celebrated multi-grid method~\cite{wesseling1995introduction}.

\begin{figure*}[htbp!]
    \centering
    \includegraphics[width=1.0\linewidth]{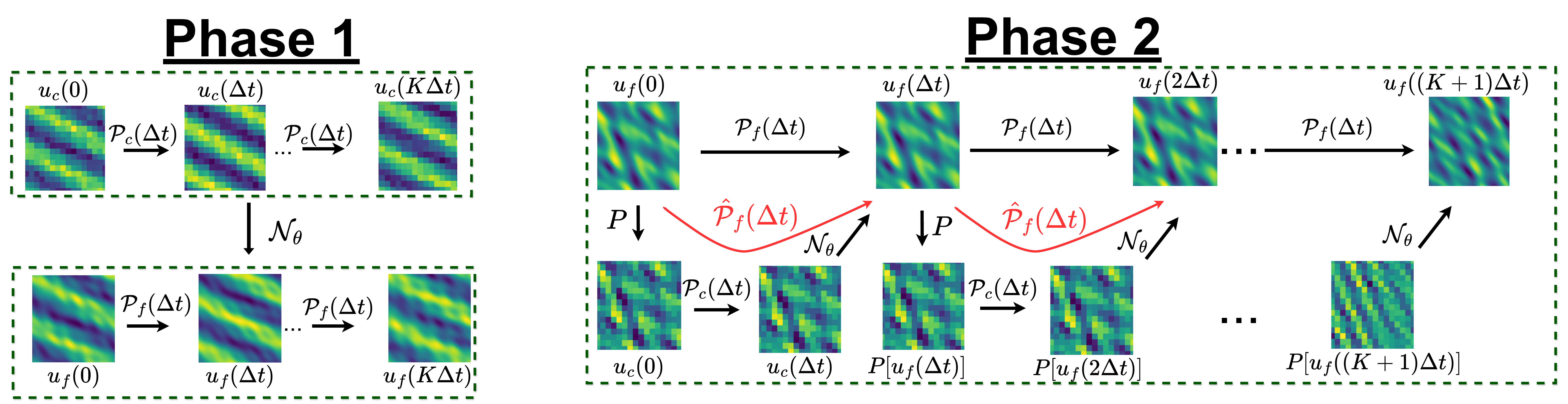}
    \caption{\textbf{Phase 1 and Phase 2 training in the RELift framework.} In Phase 1, we learn a super-resolution operator that maps the PDE solution at the coarse grid to the fine grid. In Phase 2, we use this learned operator to compose an effective propagator $\hat{\mathcal P}_{f}(\Delta t) :=\mathcal{N}_\theta\circ\mathcal{P}_c(\Delta t)\circ P$ that approximates the true fine-grid propagator $\mathcal P_f(\Delta t)$. }
    \label{fig:srafte-training}
\end{figure*}


\subsection{Main contributions}

The main contributions of this work are twofold.  
First, through systematic numerical experiments on both linear and nonlinear PDEs, we show that the effective fine-grid dynamics propagator obtained by composing a learned super-resolution operator with a coarse-grid numerical integrator enables accurate long-time prediction of time-dependent PDE dynamics.  
The resulting trajectories exhibit substantially greater stability than those produced by purely data-driven neural operators trained to learn the full solution operator.  
Second, we provide numerically verifiable local and global approximation error bounds for the effective fine-grid propagator, and we further demonstrate that these bounds closely track the observed errors in simulation, thereby explaining the improved numerical stability of RELift in long-time dynamical predictions.


\subsection{Organization}
The remainder of the paper is organized as follows.  
Section~\ref{sec:srafte} introduces the two-phase RELift framework.  
Section~\ref{sec:models} outlines the super-resolution operators used within RELift and the canonical baselines for comparison.  
Section~\ref{sec:data} describes the three PDE examples studied in this work and the evaluation metrics used to assess model performance.  
Section~\ref{sec:results} presents Phase~1 reconstruction and Phase~2 long-time forecasting results and provides further analysis of the predicted trajectories.  
Section~\ref{sec:phase2-error} conducts a systematic error analysis of RELift. In particular, we derive theoretical local and global approximation error bounds, and we numerically verify their validity and tightness for a representative nonlinear PDE.  Section~\ref{sec:eldyapp} applies RELift to a kinetic electron-dynamics case study via the 1D1V Vlasov--Poisson system, reporting both reconstruction and extrapolation performance together with physics-based diagnostics.  
Finally, Section~\ref{sec:conclusion} concludes by summarizing key findings, limitations, and directions for future research.  
Additional material may be found in the appendices; in particular, for clarity, Table~\ref{tab:notation-core} summarizes all notation and the associated definitions used in this~work.


\section{RELift: Phase 1 and Phase 2}\label{sec:srafte}
Let $\Omega \subset \mathbb{R}^D$ be a bounded domain and $X(\Omega)$ be a Banach space of functions on $\Omega$. We consider two spatial discretizations of the bounded domain: $
\Omega_f = \{x_i^f\}_{i=1}^{N_f} \subset \Omega$ and $
\Omega_c = \{x_i^c\}_{i=1}^{N_c} \subset \Omega,
$ which correspond to the fine and coarse grids, respectively. For the time-dependent PDE
\begin{align}\label{eqn:PDE}
\partial_t u = \mathcal{L}u \quad \text{in } \Omega,
\end{align}
the operator $\mathcal{L}$ is the infinitesimal generator of a strongly continuous semigroup associated with the PDE. Let $u(\cdot,t)\in X(\Omega)$ denote the solution of the initial value problem (IVP) for \eqref{eqn:PDE}.  
Its discrete fine-grid and coarse-grid representations are defined as
\begin{align*}
u_f(t) &:= \big(u(x_1^f,t),\dots,u(x_{N_f}^f,t)\big) \in \mathbb{R}^{N_f},
\\
u_c(t) &:= \big(u(x_1^c,t),\dots,u(x_{N_c}^c,t)\big) \in \mathbb{R}^{N_c}.
\end{align*}
Spatial discretization of $\mathcal{L}$ on the two grids yields approximated semi-discrete evolution equations for $u_f(t)$ and $u_c(t)$:
\begin{align}
\partial_t u_f(t) &\approx \mathcal{L}_f\, u_f(t), 
& u_f(0) &= u_f^0, \label{eq:fine}\\[4pt]
\partial_t u_c(t) &\approx\mathcal{L}_c\, u_c(t), 
& u_c(0) &= u_c^0, \label{eq:coarse}
\end{align}
where $\L_f,\L_c$ denote the fine-grid and coarse-grid discretization of $\mathcal{L}$ under the chosen boundary conditions, and $u_f^0$ and $u_c^0$ are the discrete initial conditions. Consider an equispaced time grid $t_n = n\Delta t$, $0 \le n \le N$.  
Defining 
\(
u_f^{\,n} := u_f(t_n), 
u_c^{\,n} := u_c(t_n),
\)
the exact solution can therefore be approximated by the discrete one-step evolution equation:
\begin{align*}
u_f^{\,n+1} &\approx \mathcal{P}_f(\Delta t)\,u_f^{\,n}, 
\\
u_c^{\,n+1} &\approx \mathcal{P}_c(\Delta t)\,u_c^{\,n}.
\end{align*}
We write 
\(
\mathcal{P}_f(\Delta t) = e^{\,\Delta t\,\mathcal{L}_f},
\mathcal{P}_c(\Delta t) = e^{\,\Delta t\,\mathcal{L}_c},
\)
which will be called, respectively, the fine-grid and coarse-grid (one-time) propagators for the discretized dynamics. 

RELift can be broken down into two separate learning phases. As shown in Figure \ref{fig:srafte-training}, Phase 1 can be understood as the standard super-resolution task \cite{ren2025superbench, wang2020deep, gao2021super, fukami2024single, gao2024bayesian,song2024forecasting,sun2020surrogate}: given $K$ pairs of coarse and fine grid snapshots of the PDE solution at different times $\{(u_{c}^{(i)},u_{f}^{(i)})\}_{i=0}^{K}$, we learn the neural operator $\mathcal{N}_\theta: u_c \to u_f$ by minimizing the loss function
\begin{equation}\label{eqn:L_reconstruction}
\mathcal{J}_{\mathrm{P1}}(\theta)
 =\frac{1}{K}\sum_{i=0}^{K}
  \|
      u_f^{(i)}-\mathcal{N}_\theta\bigl[u_c^{(i)}\bigr]
  \|_{\ell_1(\R^{N_f})}.
\end{equation}
Here, $\theta$ is the collection of all trainable parameters of the neural operator $\mathcal{N}_\theta$; and $\|\cdot\|_{\ell_1(\R^{N_f})}$ denotes the averaged $\ell_1$-norm, which is defined by 
\[
\|v\|_{\ell_1(\mathbb{R}^{N_f})}
:= \frac{1}{H W} \sum_{i=1}^{H} \sum_{j=1}^{W} |v_{ij}|.
\] 
The lattice grid points are $\{v_{ij}\}_{i=1,j=1}^{H,W}$, with $N_f = H W$ and entries. 
For the entire Phase~1 and Phase~2 learning tasks, the choice of neural operator can be rather flexible, and several architectures that we tested will be described later.
For future time extrapolation, only knowing the super-resolution neural operator is insufficient. In particular, we note that 
\[
\bigl(\mathcal{N}_\theta\circ\mathcal{P}_c(\Delta t)\bigr)[u_c(t)]
 \neq
\bigl(\mathcal{P}_f (\Delta t)\circ\mathcal{N}_\theta \bigr)[u_c(t)]].
\]
Therefore, simply composing $\mathcal{N}_\theta$ with the coarse propagator $\mathcal{P}_c(\Delta t)$ does not yield a good approximation of the fine‐grid propagator.  
The approximation error incurred at each step accumulates over time, causing the predicted fine‐grid dynamics to drift rapidly away from the ground truth.  Controlling this compounding error is precisely the goal of Phase~2. To this end, as shown in Figure~\ref{fig:srafte-training}, we compose the learned super‐resolution operator $\mathcal{N}_\theta$ with the coarse one‐step propagator $\mathcal{P}_c(\Delta t)$ and a projection operator
\(
P:\mathbb{R}^{N_f}\rightarrow \mathbb{R}^{N_c},
\)
to form an \emph{effective fine‐grid propagator}:
\begin{equation}\label{eqn:opcomp}
    \hat{\mathcal{P}}_f(\Delta t;\theta)
    :=\mathcal{N}_\theta \circ \mathcal{P}_c(\Delta t)\circ P.
\end{equation}
Here $\mathcal{N}_\theta$ is the super‐resolution neural operator trained in Phase~1, while  
$\mathcal{P}_c(\Delta t)$ and $P$ are fixed operators. Thus in Phase~2, we use the Phase~1 training results as initialization for $\N_{\theta}$, and we then \emph{fine‐tune} the parameters $\theta$ via a new loss function.

Specifically, we use paired snapshots of fine‐grid dynamics  
\(
\{(u_f^{(i)},u_f^{(i+1)})\}_{i=0}^{K}
\)
as the training data. Note that 
this is \textit{not} a new dataset. By (\ref{eqn:opcomp}), each fine snapshot is projected to the coarse grid, advanced by one coarse‐grid step, and then lifted back to the fine grid using the neural operator. Accordingly, we minimize the discrepancy between the true fine‐grid evolution and the effective propagator:
\begin{equation}\label{eqn:L_propagator}
\mathcal{J}_{\mathrm{P2}}(\theta)
 =\frac{1}{K}\sum_{i=0}^{K}
 \|
      \mathcal{P}_f(\Delta t)\,u_f^{(i)}
      -\hat{\mathcal{P}}_f(\Delta t;\theta)\,[u_f^{(i)}]
  \|_{\ell_1(\R^{N_f})}  .
\end{equation}  
Once trained, the learned effective propagator is used for autoregressive forecasting:
\[
\hat{u}_f^{\,n+1}
    =\hat{\mathcal{P}}_f(\Delta t;\theta)\,\hat{u}_f^{\,n},
    \qquad n\ge K.
\]
After $m$ additional steps (i.e., at endtime $T=m\Delta t$), the fine‐grid forecast is
\[
\hat{u}_f(K\Delta t+T)
= \bigl[\hat{\mathcal{P}}_f(\Delta t;\theta)\bigr]^{m}\,u_f^{K}.
\]
In principle, the effective fine-grid propagator introduces only minimal computational overhead.  
The projection operator $P$ is a fixed matrix, and the evaluation of the trained neural operator $\mathcal{N}_{\theta}$ is inexpensive since as it can be seen as a sequence of function compositions. Thus, the main computational cost at each time step arises from the coarse-grid update $\mathcal{P}_c(\Delta t)$. In the whole process, numerical discretization errors and model errors still accumulate over time. 

Nevertheless, because the composed operator incorporates the exact coarse-grid propagator---which preserves the macro-scale features of the underlying PDE---we expect the resulting dynamics to be more stable in future dynamics prediction, when compared to pure data-driven autoregressive models commonly used in the literature \cite{li2021fourier, chattopadhyay2023long, darman2025fourier}, which lack such physically grounded structure and which therefore tend to suffer from rapid error growth during extrapolation.


From a training perspective, training via RELift reduces the computational burden of training the neural network, as it is solely responsible for SR. First, the same neural operator architecture is used in both training stages; only the loss function changes between Phase~1 and Phase~2. Second, no backpropagation through the coarse solver is ever required.  
Once the coarse solver is determined, it remains fixed, and gradients from backpropagation flow solely through the neural operator $\mathcal{N}_\theta$, whose role is to perform super-resolution.  
This significantly simplifies training and reduces computational overhead compared with end-to-end approaches that differentiate through multi-step PDE solvers.

\textbf{Remark.}
The effective fine-grid propagator $\mathcal{P}_f(\Delta t)$ in RELift can be interpreted as providing a Markovian closure for the fine-grid dynamics.  A related approach, SHRED~\cite{kutz2024shallow}, is similar in spirit but instead learns a finite-memory, history-dependent closure $\Phi_h(u_c(t;\tau))$ using a shallow recurrent encoder.  
In contrast, RELift delegates temporal evolution entirely to the coarse-grid time propagator and restricts learning to a \emph{time-local} super-resolution operator.  
This perspective suggests a separation of variables: $\mathcal{P}_c(\Delta t)$ governs the (approximately) low-dimensional temporal evolution, while $\mathcal{N}_\theta$ provides the high-dimensional spatial reconstruction. Additional intuition, including linear and nonlinear heuristics, a Mori--Zwanzig interpretation~\cite{zhu2018estimation, zhu2019mori, zhu2021effective}, and a discussion of when short-memory assumptions are justified or may fail, is provided in Appendix~\ref{sec:stmem-details}.


\section{Models and baselines}\label{sec:models}

RELift imposes no \textit{a priori} constraints on the architecture of the neural operator $\mathcal{N}_{\theta}$. In principle, any neural operator framework developed for super-resolution tasks can be adapted into the two-phase RELift pipeline, which is used first to learn and super-resolve the coarse-scale time dynamics, and then to form an effective one-time propagator for predicting future fine-scale dynamics.  In this section, we introduce four super-resolution operators that we use within RELift: U-Net\cite{ronneberger2015u}, FUnet (our architecture), EDSR\cite{lim2017enhanced}, and FNO-SR\cite{li2021fourier} and the baselines for benchmarking the numerical results. Complete model setups and hyperparameter definitions are deferred to Appendix~\ref{sec:model-summaries}.

\subsubsection{\textbf{U-Net}}
The U-Net \cite{ronneberger2015u} is a classical encoder-decoder architecture widely used in computer vision for image super-resolution. We use U-Net as a classical convolutional based networks within RELift and test its applicability in Phase 2 dynamics predictions.

\subsubsection{\textbf{Fourier U-Net (FUnet)}}
FUnet, a new architecture introduced in this work, augments the U-Net with a spectral convolutional layer. Specifically, between the encoder and decoder in the bottleneck, we add a spectral convolutional block. The FUnet lifts the features $z \in \mathbb{R}^{C_0}$ from the encoder and computes their discrete Fourier transform: $\widehat{z} = \mathcal{F}[z]$. Subsequently, the lowest $R$ modes are retained and transformed back into the time domain. The decoder then upsamples to the target resolution in the decoder using pixel-shuffle layers with skip connections from the encoder \cite{shi2016real}. The FUnet architecture fuses spectral features with convolutional layers, making it effective both for super-resolution and for extracting scientifically meaningful structure from high-resolution data.

\subsubsection{\textbf{Enhanced Deep Super-Resolution (EDSR)}}
The image super-resolution EDSR model \cite{lim2017enhanced} is adopted here with minimal modifications and placed within RELift. Like the U-Net, the RGB channels used in image processing are replaced by field channels representing the temporal trajectories of the PDE solution. The low-resolution input is first processed by 16 residual blocks, followed by a two-stage pixel-shuffle to achieve the desired resolution increase. A final convolution layer then produces the high-resolution output predictions. EDSR is another convolutional-based model tested within RELift.

\subsubsection{\textbf{Fourier neural operator (FNO)}}
The FNO was originally proposed in \citet{li2021fourier} to learn a one-step propagator for PDEs (and its zero-shot super-resolution behavior has recently been characterized~\cite{FalsePromizeZeroShot_TR}). However, as a general neural operator framework, it can also be adapted as a modeling ansatz for a super-resolution operator $\mathcal{N}_{\theta}$. Specifically, it can be trained to take the numerical solution of a PDE on a coarse grid as input and return the corresponding fine grid solution. This super-resolution variant of FNO (referred to as FNO-SR) can then be used in Phase~2 to form an effective propagator for PDE dynamics. Note that this usage of FNO differs from applying FNO to directly approximate the one-step propagator, which will be referred to as the autoregressive FNO (FNO-AR) in the following discussion.

\subsubsection{\textbf{Upsampling}}
Upsampling refers to interpolation methods that directly maps the corase-grid solution into fine-grid, typically via interpolation kernels such as nearest neighbor, bilinear, or bicubic. For this work, we adopt bicubic interpolation as a simple baseline for assessing the performance of RELift in fine-scale dynamics reconstruction and prediction. In Phase~1, it maps the coarse grid PDE solution at each timestep to a fine grid solution by evaluating the standard cubic convolution kernel in all spatial directions. In Phase~2, the PDE is first solved on the coarse grid, and the resulting coarse-grained solution at each timestep is then upsampled using the same interpolation. Although simplistic, this baseline is helpful for quantifying how much of RELift’s advantage stems from operator learning versus mere smooth spatial upsampling. 

\medskip
\textbf{Rationale. }The selected models and baselines include simple nonparametric interpolation, established super-resolution methods, and autoregressive models for dynamics forecasting. Together, these baselines allow us to compare RELift against the two main alternatives:  
(i) spatial-only refinement methods and  
(ii) purely data-driven autoregressive models for temporal dynamics prediction.  
Within this context, the chosen baselines are sufficient to demonstrate the advantage of RELift's operator-composition framework for dynamics predictions.



\section{Datasets and evaluation metrics}\label{sec:data}
We now introduce the PDEs we study, briefly describing the generation of the training datasets for each example and the evaluation metrics used to assess the predicted dynamics. More details are provided in Appendices~\ref{app:data_generation} and~\ref{sec:metrics}.

All systems use periodic boundary conditions, a uniform discretization of the square domain $[0,1] \times [0,1]$ for both the coarse and fine grids, and a timestep $\Delta t=0.01$. We take the projection operator $P$ to subsample every $k$-th point from the fine grid (i.e., traditional point sampling/spatial decimation).

\paragraph{\textbf{2D heat equation.}}
Consider the 2D heat equation $\partial_t u = \nu \Delta u + \alpha f$ with parameters $\nu = 10^{-3}$, $\alpha = 10^{-2}$, and forcing term $f(x) = \sin(2\pi x)\sin(2\pi y)$. The equation is solved using a Fourier spectral method.

\paragraph{\textbf{2D wave equation.}}
As an example of an energy-conserving PDE, we consider the linear wave equation $\partial_{tt} u = c^2 \Delta u$ with $c = 0.5$. The equation is solved using an explicit finite-difference leapfrog scheme on both coarse and fine spatial grids to generate the training~datasets.

\paragraph{\textbf{Navier–Stokes equations (NSE).}} 
As an example of a nonlinear PDE, we consider the 2D incompressible NSE in vorticity form $\partial_t\omega + u \!\cdot\! \nabla\omega = \nu \Delta\omega + f$, $\nabla \cdot u = 0$ with $\nu = 10^{-4}$ and forcing term $f(x) = 0.025\bigl[\sin 2\pi(x{+}y) + \cos 2\pi(x{+}y)\bigr]$. The NSE is solved using a Fourier spectral method.
To close the system, we introduce a streamfunction $\psi$ and set
$u=\nabla^\perp \psi := (\partial_{y}\psi,-\partial_{x}\psi)$ with $-\Delta \psi=\omega$.


\begin{figure*}[htbp!] \vspace{-0.5cm} 
  \centering
  
  \begin{subfigure}[b]{\linewidth}
    \begin{minipage}[c]{1.4cm}
      \centering
      \rotatebox{90}{\textbf{\textcolor{black}{Heat}}}
    \end{minipage}%
    \begin{minipage}[c]{\dimexpr\linewidth-1.4cm\relax}
      \centering
      \includegraphics[width=\linewidth]{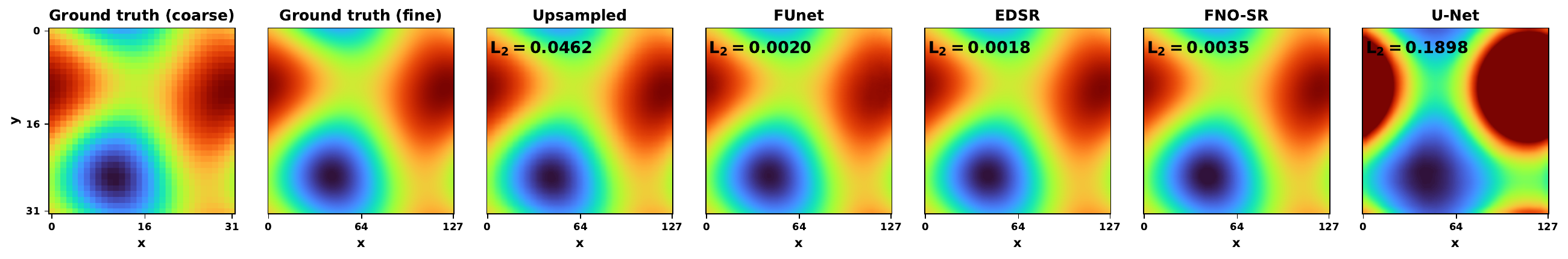}
    \end{minipage}
    \label{fig:heat-field-080}
  \end{subfigure}

  \vspace{0.6em}

  \begin{subfigure}[b]{\linewidth}
    \begin{minipage}[c]{1.4cm}
      \centering
      \rotatebox{90}{\textbf{\textcolor{black}{Wave}}}
    \end{minipage}%
    \begin{minipage}[c]{\dimexpr\linewidth-1.4cm\relax}
      \centering
      \includegraphics[width=\linewidth]{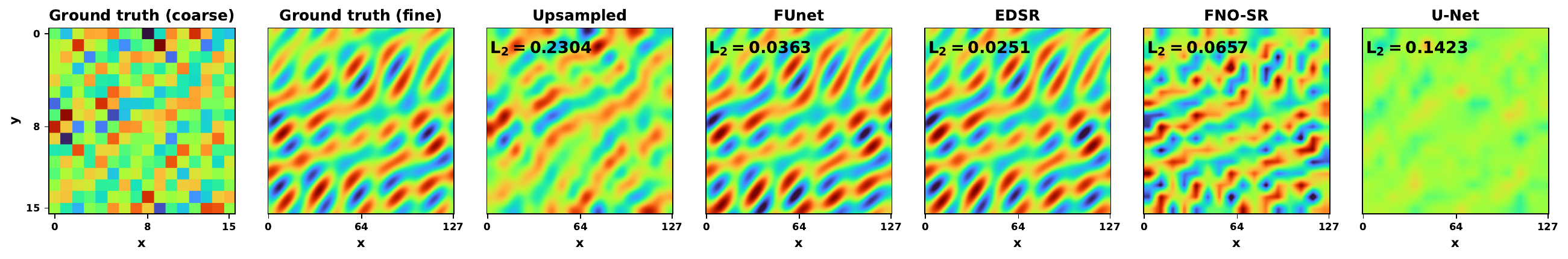}
    \end{minipage}
    \label{fig:wave-field-070}
  \end{subfigure}

  \begin{subfigure}[b]{\linewidth}
    \begin{minipage}[c]{1.4cm}  
      \centering
      \rotatebox{90}{\textcolor{black}{\textbf{Navier--Stokes}}}
    \end{minipage}%
    \begin{minipage}[c]{\dimexpr\linewidth-1.4cm\relax}
      \centering
      \includegraphics[width=\linewidth]{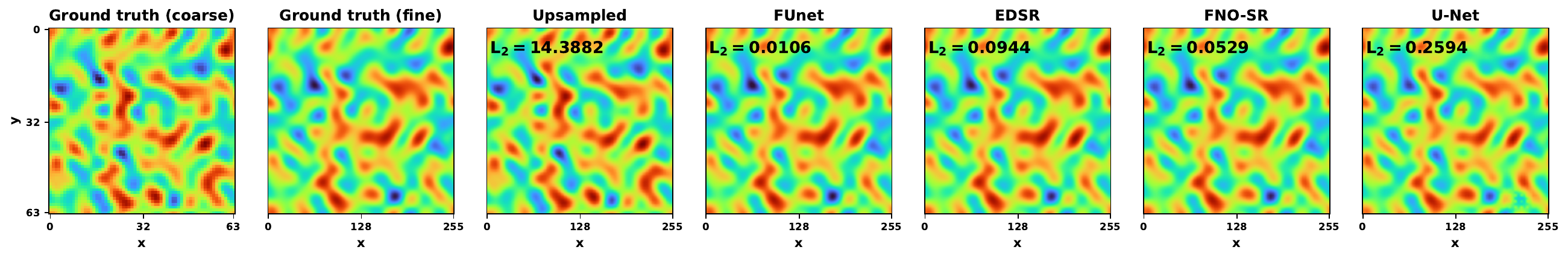}
    \end{minipage}
    \label{fig:ns-field-090}
  \end{subfigure}
\caption{
\justifying
\textbf{Phase~1 super-resolution results across different PDEs.} For each PDE, we show results for an unseen test trajectory, with representative snapshots from the \textbf{(Top)} 2D heat equation, \textbf{(Middle)} 2D wave equation, and \textbf{(Bottom)} 2D Navier--Stokes equations (NSE).  
The first column shows the true fine-grid dynamics.  
The next two columns—labeled ``Upsampled'' and ``FNO-AR''—correspond to the bicubic interpolation and autoregressive FNO baselines, respectively.  
All subsequent columns display predictions from different RELift models.  
The inset in each panel reports the relative $L_2$ reconstruction error.}
\label{fig:phase1-fields}
\end{figure*}

\squeezetable
\begin{table*}[ht!]
\setlength{\tabcolsep}{2pt}        
\renewcommand{\arraystretch}{0.95} 
\centering
\scriptsize
\newcolumntype{b}{>{\columncolor{gray!8}}c}
{\begin{tabularx}{\linewidth}{@{}l b c c c c@{}}
\toprule
& \multicolumn{1}{c}{\cellcolor{gray!15}\textbf{Baseline}} 
& \multicolumn{4}{c}{\cellcolor{gray!6}\textbf{RELift models}} \\
\cmidrule(lr){2-2}\cmidrule(lr){3-6}
\textbf{System} & \cellcolor{gray!15}\textbf{Upsampled} & \textbf{FUnet} & \textbf{U-Net} & \textbf{EDSR} & \textbf{FNO-SR} \\
\midrule

\textbf{Heat equation}
  & $5.37\times10^{-2}$
  & $2.08\times10^{-3}\pm 1.39\times10^{-4}$
  & $2.23\times10^{-1}\pm 1.39\times10^{-4}$
  & $\mathbf{6.02\times10^{-4}}\pm 8.49\times10^{-5}$
  & $3.85\times10^{-3}\pm 4.69\times10^{-4}$ \\[0.3em]

\textbf{Wave equation}
  & $2.23\times10^{-1}$
  & $3.88\times10^{-2}\pm 2.96\times10^{-3}$
  & $2.53\times10^{-1}\pm 1.08\times10^{-2}$
  & $\mathbf{3.43\times10^{-2}}\pm 2.61\times10^{-3}$
  & $7.02\times10^{-2}\pm 2.33\times10^{-3}$ \\[0.3em]
  
\textbf{Navier--Stokes}
  & $1.47\times10^{1}$
  & $\mathbf{1.25\times10^{-2}}\pm 4.82\times10^{-3}$
  & $1.55\times10^{-1}\pm 6.19\times10^{-2}$
  & $4.20\times10^{-2}\pm 2.89\times10^{-3}$
  & $3.07\times10^{-2}\pm 4.96\times10^{-3}$ \\[0.3em]
\bottomrule
\end{tabularx}}
\caption{
\justifying
\textbf{Relative $L_2$ errors of Phase~1 super-resolution results.} Each error is obtained by averaging over 20 new test trajectories. The shaded column corresponds to the bicubic upsampling baseline, while the remaining columns represent different RELift models. The lowest error in each row is highlighted in bold.}
\label{tab:phase1_toy}
\end{table*}

\subsection{{\textbf{Training data}}} For all three examples, we generate 1000 trajectories using random initial conditions detailed in Appendix \ref{app:data_generation} to form the training datasets. The neural operator $\mathcal{N}_{\theta}$ is trained on a paired dataset of coarse-fine dynamics with an 80/20 training/validation split. Specifically, unless otherwise stated, the following coarse-to-fine grid mappings are used: $32\times 32 \to 128\times 128$ (heat equation), $16\times 16 \to 128\times 128$ (wave equation) and $64\times 64 \to 256\times 256$ (NSE). These ratios are chosen after experimenting with a range of downsampling factors, which clearly demonstrate the benefits of super-resolution and future-time extrapolation, particularly in regimes with pronounced multiscale structure. 

In Phase~2, the same neural operator $\mathcal{N}_{\theta}$ is used to construct an effective propagator $\hat{\mathcal{P}}_f(\Delta t; \theta)$, which is retrained using the fine grid data from Phase 1 over a slightly shifted time window $t \!\in\! [0, 1+\Delta t]$ and subsequently used for dynamics extrapolation. For the heat equation, extrapolation is performed until physical time $t = 2$, when the system reaches a steady state. For the wave equation and NSE, the PDEs are evolved up to physical time $t = 10, 7$, respectively, to probe strongly out-of-distribution (OOD) dynamics.

\subsection{\textbf{Metrics}} 
To assess the generalizability of the learned super-resolution operator and effective propagator for PDE dynamics prediction, we evaluate the model-predicted dynamics in both phases using a new set of solutions generated from initial conditions randomly sampled from the same distribution. Specifically, we test 20 new trajectories, and we report the mean and standard deviation across three independent training runs. For Phase~1, we calculate the relative $L_2$ error across all predictions, while for Phase~2, we compute each of the metrics described below. Full details of each metric can be found in Appendix~\ref{sec:metrics}.

\paragraph{\textbf{Why multiple metrics?}}
As suggested in \citet{ren2025superbench}, simple pixel-wise losses alone can be misleading for scientific data, as lower pixel-wise error does not guarantee physically faithful solutions. Thus, it is necessary to use different metrics to evaluate the performance of our models, as no single scalar score can capture every facet of fidelity that matters for high-resolution scientific data. Below we characterize the key features of each metric. 
\begin{itemize}
    \item Relative $L_2$ error (\textbf{L}$_\mathbf{2}$) is \textit{scale-agnostic} but \textit{location-sensitive}: any pixel-wise mismatch is penalized,
    so large-amplitude errors anywhere in the field dominate. In ML contexts, this is one of the most standard error metrics.
    
    \item The Structural Similarity Index Measure (\textbf{SSIM}) is a perceptual metric designed to assess local, pixel-wise agreement in texture, contrast, and structural patterns.  Two fields may be close in the \(L_2\) norm yet differ significantly in visual or structural appearance; SSIM is specifically sensitive to such discrepancies.  By definition, $\text{SSIM} \in [-1,1]$, with a value of \(1\) indicating perfect structural fidelity between the reconstructed and reference fields.

    \item Spectral mean-squared error (\textbf {Spectral MSE}) only measures the \textit{distribution} of energy across Fourier modes. It is particularly sensitive to aliasing, mode transfer distortions, or over-smoothing effects that may not be strongly reflected in either the \(L_2\) error or SSIM.
    
    \item Pearson correlation (\textbf{Corr}) measures the \textit{linear correlations} between the predicted and reference fields.  Because it is computed after centering each field, \(r\) is \textit{scale-invariant}: two fields that differ only by a global multiplicative factor still achieve \(r=1\).  Thus, \(r\) complements the other metrics by highlighting how well the spatial patterns vary, irrespective of absolute amplitude.
\end{itemize}

\section{Numerical results}
\label{sec:results}

\subsection{Dynamics super-resolution (Phase 1 test results)}
For all the datasets considered, the Phase~1 dynamics reconstruction results and the corresponding evaluation metrics are summarized in Figure~\ref{fig:phase1-fields} and Table~\ref{tab:phase1_toy}. Specifically, Figure~\ref{fig:phase1-fields} shows a representative snapshot of the PDE solution at a chosen time ($t = 0.8$, $0.7$, and $0.9$ for the heat equation, wave equation, and NSE, respectively). Quantitatively, Table~\ref{tab:phase1_toy} reports the $L_2$ errors averaged over all 20 new test trajectories. These results clearly show that the super-resolution models consistently outperform the naive bicubic upsampling baseline by $1$–$3$ orders of magnitude, confirming the utility of the ML-based super-resolution approach for dynamics reconstruction.

For the three examples considered, the heat equation exhibits the simplest dynamics because it is a linear parabolic PDE with constant coefficients, corresponding to a gradient flow that smooths and decouples Fourier modes via pure exponential decay. As a result, all tested super-resolution models successfully capture the main features of the flow and produce quantitatively more accurate predictions compared to the other two PDEs. The dynamics become more complex for the wave equation, and they are further complicated in the Navier–Stokes equation due to nonlinearity. Consequently, as shown in Figure~\ref{fig:phase1-fields} and Table~\ref{tab:phase1_toy}, the discrepancy between super-resolution models and baseline bicubic interpolation increases from the heat equation to the wave equation and then to the NSE. Among all tested RELift models, no significant difference is observed in the Phase~1 results between FUnet, EDSR, and FNO-SR; these three models consistently yield more accurate predictions of the PDE solutions than the purely convolutional-based U-Net model.

\begin{figure*}[t]
    \centering
    \begin{minipage}[t]{0.49\textwidth}
        \centering
        \includegraphics[width=\linewidth]{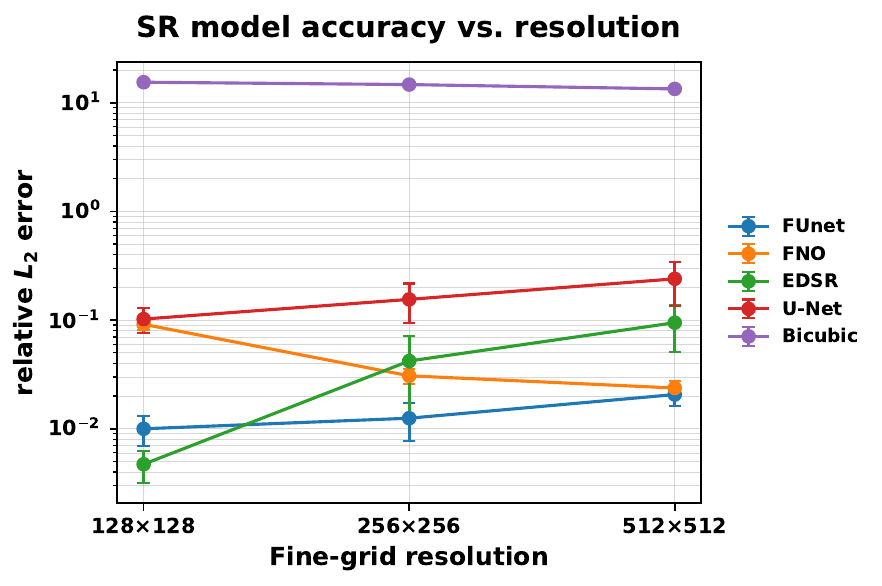}
        \vspace{0.25em}
        {\small (a) Accuracy of the super-resolution results for different models with constant downsampling factor of $\times 4$.}
    \end{minipage}\hfill
    \begin{minipage}[t]{0.49\textwidth}
        \centering
        \includegraphics[width=\linewidth]{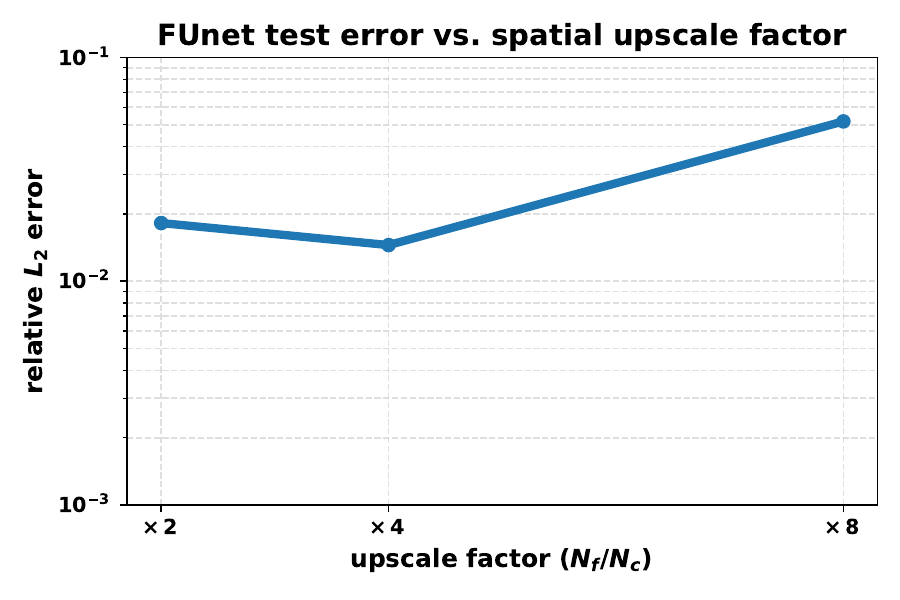}
        \vspace{0.25em}
        {\small (b) Accuracy of the FUnet super-resolution results with different downsamplings factors.}
    \end{minipage}
    \caption{\justifying 
    \textbf{Phase~1 results across different resolutions and scaling factors.}  
\textbf{(a)} Relative $L_2$ errors for all models and baselines evaluated at multiple fine-grid resolutions.  
Each curve is obtained by averaging over 20 unseen trajectories.  
Models that explicitly exploit spectral structure (FUnet, FNO) achieve the best accuracy, with FUnet exhibiting the slowest growth in error as the fine-grid resolution increases.  
\textbf{(b)} Relative $L_2$ error of FUnet at a fixed fine-grid resolution of $256\times256$ while varying the coarse-grid input size $N_c \in \{128,64,32\}$ (corresponding to upscale factors $\times 2, \times 4, \times 8$).
 }
    \label{fig:phase1-rescurve}
\end{figure*}

To test the scale-robustness of Phase 1 results, we repeat the NS experiment at two additional coarse-fine pair resolutions, $32\times32\to128\times128$ and $128\times128\to512\times512$ for all models and baselines. The resulting resolution sweep, presented in Figure~\ref{fig:phase1-rescurve}, tracks the mean and standard deviation of the relative $L_2$ error for every model as fine grid resolution increases. Across all three resolutions, we see that the FUnet remains the most accurate and shows the slowest error growth. Interestingly, the FNO-SR curve bends downward as we move from the lowest target resolution to the highest, indicating that the model becomes more accurate. 
This could be because the FNO is inherently biased toward low- and mid-frequency content because each spectral layer retains only a prescribed band of modes \cite{qin2024toward, khodakarami2025mitigating}. When the fine grid is enlarged while the down-sampling ratio remains fixed, most dynamically important structures now fall inside that retained band, so the proportion of energy that the network must extrapolate beyond its spectral cut-off is reduced.


To further assess the sensitivity of the Phase~1 super-resolution results to input resolution, we evaluate the best-performing model (FUnet) at a fixed fine-grid output resolution of $256\times256$ while varying the coarse-grid input over $N_c \in \{128,64,32\}$ (corresponding to upscale factors $2,4,8$).  
As shown in the right panel of Figure~\ref{fig:phase1-rescurve}, the relative $L_2$ error increases as the coarse resolution decreases, reflecting the expected difficulty of reconstructing fine-scale features from coarser inputs.  
However, the growth in error remains slow, indicating that the super-resolution performance remains accurate and stable across a broad range of spatial scales.  

As the final probe for generalization of Phase~1, we evaluate the Phase~1 FUnet trained at $\nu=10^{-4}$ on unseen viscosity datasets with  $\nu=10^{-5}$ and $\nu=10^{-3}$. As shown in Table~\ref{tab:phase1_ns_funet_multi_nu}, the learned lift remains substantially more accurate than the bicubic upsampling baseline across all three viscosities, with the strongest performance at the training viscosity and only moderate degradation under viscosity shift; full details are given in Section~\ref{sec:ablation}.

Taken together, these results suggest that the Phase~1 trained FUnet is capable of learning reliable fine-grid reconstructions from coarse-grid inputs across a variety of resolutions.

\newcommand{\nonnswidth}{1.0\linewidth}
\newcommand{\pdetime}[3]{
  \begin{minipage}[c]{1.6cm}           
    \centering
    \begin{minipage}[c]{0.8cm}         
      \centering
      \rotatebox{90}{\textbf{\textcolor{#1}{#2}}}
    \end{minipage}%
    \begin{minipage}[c]{0.8cm}         
      \centering
      \rotatebox{90}{\textcolor{black}{#3}}
    \end{minipage}
  \end{minipage}%
}

\begin{figure*}[t]
  \centering
  \begin{subfigure}[b]{\linewidth}
    \pdetime{black}{Heat}{$t=1.5$}%
    \begin{minipage}[c]{\dimexpr\linewidth-1.6cm\relax}
      \centering
      \includegraphics[width=\nonnswidth]{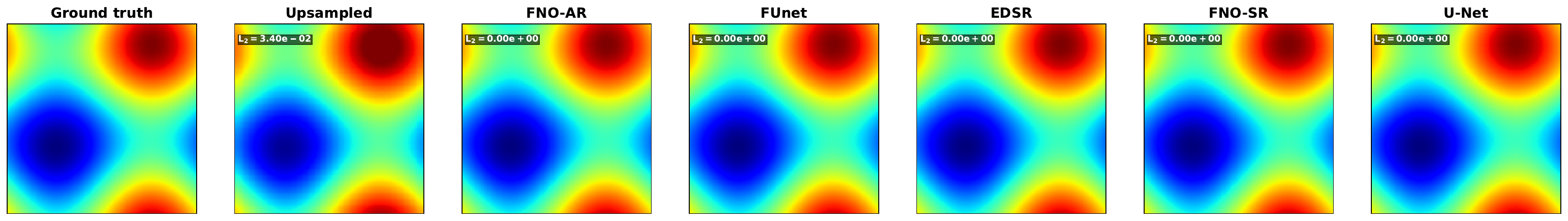}
    \end{minipage}
  \end{subfigure}

  \vspace{0.6em}

  \begin{subfigure}[b]{\linewidth}
    \pdetime{black}{Heat}{$t=2.0$}%
    \begin{minipage}[c]{\dimexpr\linewidth-1.6cm\relax}
      \centering
      \includegraphics[width=\nonnswidth]{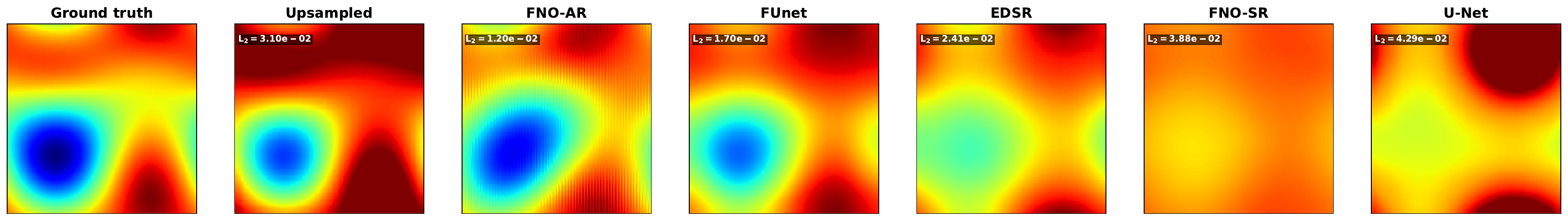}
    \end{minipage}
  \end{subfigure}

  \vspace{0.6em}

  \begin{subfigure}[b]{\linewidth}
    \pdetime{black}{Wave}{$t=2.0$}%
    \begin{minipage}[c]{\dimexpr\linewidth-1.6cm\relax}
      \centering
      \includegraphics[width=\nonnswidth]{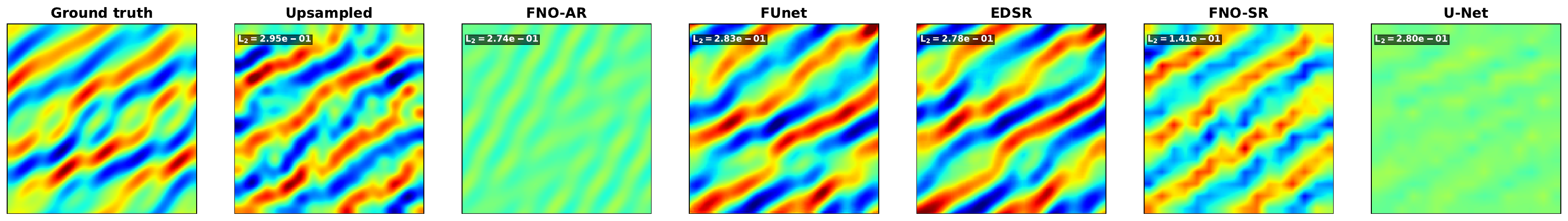}
    \end{minipage}
  \end{subfigure}

  \vspace{0.6em}

  \begin{subfigure}[b]{\linewidth}
    \pdetime{black}{Wave}{$t=5.0$}%
    \begin{minipage}[c]{\dimexpr\linewidth-1.6cm\relax}
      \centering
      \includegraphics[width=\nonnswidth]{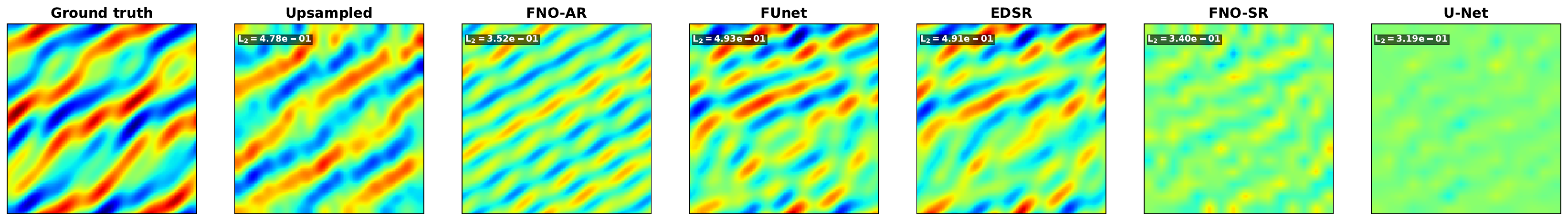}
    \end{minipage}
  \end{subfigure}

  \begin{subfigure}[b]{\linewidth}
    \pdetime{black}{Navier--Stokes}{$t=2.0$}%
    \begin{minipage}[c]{\dimexpr\linewidth-1.6cm\relax}
      \centering
      \includegraphics[width=\linewidth]{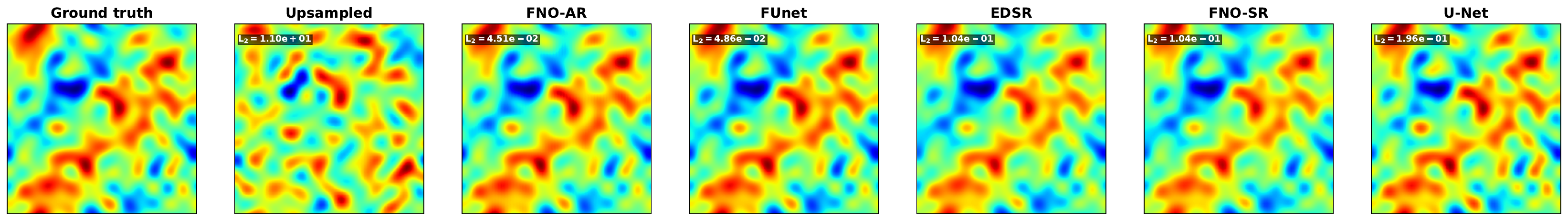}
    \end{minipage}
  \end{subfigure}

  \vspace{0.6em}

  \begin{subfigure}[b]{\linewidth}
    \pdetime{black}{Navier--Stokes}{$t=7.0$}%
    \begin{minipage}[c]{\dimexpr\linewidth-1.6cm\relax}
      \centering
      \includegraphics[width=\linewidth]{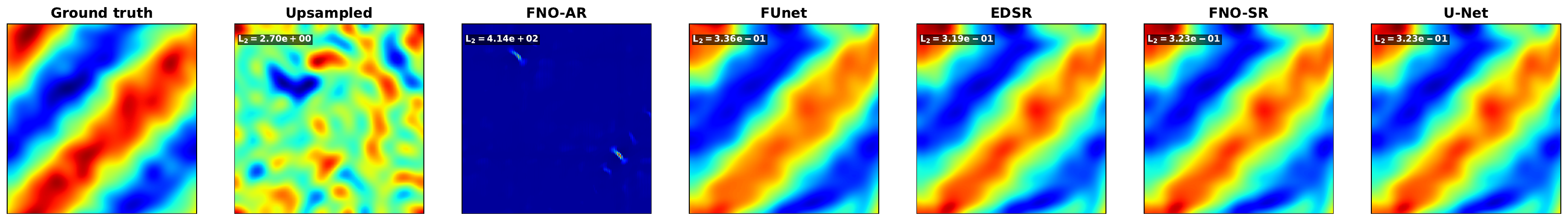}
    \end{minipage}
  \end{subfigure}

  \vspace{0.6em}

  \caption{\justifying 
  \textbf{Phase~2 future dynamics prediction results accross different PDEs}. For each PDE, we show results for an unseen test trajectory, with representative snapshots from the \textbf{(Top)} 2D heat equation, \textbf{(Middle)} 2D wave equation, and \textbf{(Bottom)} 2D Navier--Stokes equations (NSE).  The first column shows the true fine-grid dynamics.  The next two columns—labeled ``Upsampled'' and ``FNO-AR''—correspond to the bicubic interpolation and autoregressive FNO baselines, respectively.  All subsequent columns display predictions from different RELift models.  The inset in each panel reports the relative $L_2$ error.}
  \label{fig:phase2-fields}
\end{figure*}

\squeezetable 

\begin{table*}[ht!]
\setlength{\tabcolsep}{2pt}        
\renewcommand{\arraystretch}{0.95} 
\centering
\scriptsize
\begin{adjustbox}{max width=\textwidth}
\begin{tabular}{@{}l l >{\columncolor{gray!10}}c >{\columncolor{gray!10}}c c c c c@{}}
\toprule
& & \multicolumn{2}{c}{\cellcolor{gray!20}\textbf{Baselines}} & \multicolumn{4}{c}{\cellcolor{gray!8}\textbf{RELift models}} \\
\cmidrule(lr){3-4}\cmidrule(lr){5-8}
\textbf{System} & \textbf{Metric} &
\cellcolor{gray!20}\textbf{Upsampled} &
\cellcolor{gray!20}\textbf{FNO-AR} &
\textbf{FUnet} & \textbf{EDSR} & \textbf{FNO-SR} & \textbf{U-Net} \\
\midrule
\multirow[t]{4}{*}{\makecell[l]{\textbf{Heat}\\\textbf{equation}}}
 & rel.\,$L_2$ ($\downarrow$)
   & $3.66\!\times\!10^{-2}$
   & $\mathbf{7.07\!\times\!10^{-4}}\pm 2.9\!\times\!10^{-5}$
   & $1.11\!\times\!10^{-3}\pm 4.5\!\times\!10^{-4}$
   & $1.47\!\times\!10^{-3}\pm 6.0\!\times\!10^{-4}$
   & $2.10\!\times\!10^{-3}\pm 8.6\!\times\!10^{-4}$
   & $2.23\!\times\!10^{-3}\pm 9.1\!\times\!10^{-4}$ \\[0.15em]
 & SSIM ($\uparrow$)
   & $\mathbf{0.9962}$
   & $0.9923\pm 0.0316$
   & $0.9950\pm 0.0215$
   & $0.9891\pm 0.0447$
   & $0.9844\pm 0.0639$
   & $0.9806\pm 0.0794$ \\[0.15em]
 & Spec.\ MSE ($\downarrow$)
   & $4.17\!\times\!10^{-2}$
   & $\mathbf{1.45\!\times\!10^{-4}}\pm 5.9\!\times\!10^{-5}$
   & $5.89\!\times\!10^{-4}\pm 2.4\!\times\!10^{-5}$
   & $8.80\!\times\!10^{-4}\pm 3.6\!\times\!10^{-5}$
   & $2.08\!\times\!10^{-3}\pm 8.5\!\times\!10^{-5}$
   & $1.55\!\times\!10^{-3}\pm 6.3\!\times\!10^{-4}$ \\[0.15em]
 & Corr ($\uparrow$)
   & $\mathbf{0.9987}$
   & $0.9964\pm 0.0148$
   & $0.9982\pm 0.0079$
   & $0.9935\pm 0.0266$
   & $0.9955\pm 0.0183$
   & $0.9773\pm 0.0922$ \\[0.3em]
\midrule
\multirow[t]{4}{*}{\textbf{Wave}}
 & rel.\,$L_2$ ($\downarrow$)
   & $3.31\!\times\!10^{-1}$
   & $2.00\!\times\!10^{-1}\pm 1.1\!\times\!10^{-2}$
   & $2.22\!\times\!10^{-1}\pm 1.4\!\times\!10^{-2}$
   & $2.52\!\times\!10^{-1}\pm 1.5\!\times\!10^{-3}$
   & $\mathbf{1.74\!\times\!10^{-1}}\pm 1.2\!\times\!10^{-3}$
   & $2.06\!\times\!10^{-1}\pm 1.2\!\times\!10^{-1}$ \\[0.15em]
 & SSIM ($\uparrow$)
   & $0.1697$
   & $0.2578\pm 0.0340$
   & $0.3836\pm 0.0418$
   & $0.3667\pm 0.0419$
   & $\mathbf{0.4820}\pm 0.0383$
   & $0.3285\pm 0.0371$ \\[0.15em]
 & Spec.\ MSE ($\downarrow$)
   & $1.13\!\times\!10^{2}$
   & $1.91\!\times\!10^{2}\pm 1.4\!\times\!10^{1}$
   & $8.75\!\times\!10^{1}\pm 8.2\!\times\!10^{0}$
   & $8.68\!\times\!10^{1}\pm 1.1\!\times\!10^{1}$
   & $\mathbf{7.52\!\times\!10^{1}}\pm 7.1\!\times\!10^{0}$
   & $1.04\!\times\!10^{2}\pm 2.1\!\times\!10^{1}$ \\[0.15em]
 & Corr ($\uparrow$)
   & $0.0479$
   & $0.2628\pm 0.0461$
   & $0.4072\pm 0.0536$
   & $0.3830\pm 0.0553$
   & $\mathbf{0.5393}\pm 0.0512$
   & $-0.0764\pm 0.0606$ \\
\midrule
\multirow[t]{4}{*}{\makecell[l]{\textbf{Navier--}\\\textbf{Stokes}}}
 & rel.\,$L_2$ ($\downarrow$)
   & $7.92\!\times\!10^{0}$
   & $1.58\!\times\!10^{1}\pm 5.4\!\times\!10^{-1}$
   & $\mathbf{9.36\!\times\!10^{-2}}\pm 7.0\!\times\!10^{-3}$
   & $1.13\!\times\!10^{-1}\pm 6.5\!\times\!10^{-3}$
   & $1.13\!\times\!10^{-1}\pm 6.5\!\times\!10^{-3}$
   & $1.68\!\times\!10^{-1}\pm 7.0\!\times\!10^{-2}$ \\[0.15em]
 & SSIM ($\uparrow$)
   & $0.0053$
   & $0.453\pm 0.409$
   & $\mathbf{0.9410}\pm 0.0518$
   & $0.9210\pm 0.0565$
   & $0.9209\pm 0.0573$
   & $0.8772\pm 0.0582$ \\[0.15em]
 & Spec.\ MSE ($\downarrow$)
   & $2.76\!\times\!10^{4}$
   & $5.87\!\times\!10^{6}\pm 3.0\!\times\!10^{5}$
   & $1.05\!\times\!10^{1}\pm 1.4\!\times\!10^{0}$
   & $\mathbf{9.11\!\times\!10^{0}}\pm 1.0\!\times\!10^{0}$
   & $9.19\!\times\!10^{0}\pm 1.0\!\times\!10^{0}$
   & $1.94\!\times\!10^{1}\pm 1.4\!\times\!10^{0}$ \\[0.15em]
 & Corr ($\uparrow$)
   & $0.3317$
   & NaN
   & $\mathbf{0.9951}\pm 0.0046$
   & $0.9924\pm 0.0058$
   & $0.9924\pm 0.0059$
   & $0.9887\pm 0.0067$ \\[0.3em]
\bottomrule
\end{tabular}
\end{adjustbox}
\caption{\justifying 
\textbf{Phase~2 future time dynamics prediction across different PDEs.}  
Each error value is averaged over 20 unseen test trajectories. The shaded baseline columns correspond to bicubic Upsampling and FNO-AR, while the remaining columns show results from various RELift models.  
The best value in each row is highlighted in bold.}
\label{tab:phase2_toy}
\end{table*}

\begin{figure*}[t]
    \centering
    \setlength{\tabcolsep}{4pt}
    \begin{tabular}{cc}
        \begin{minipage}[t]{0.48\linewidth}
            \centering
            \includegraphics[width=\linewidth]{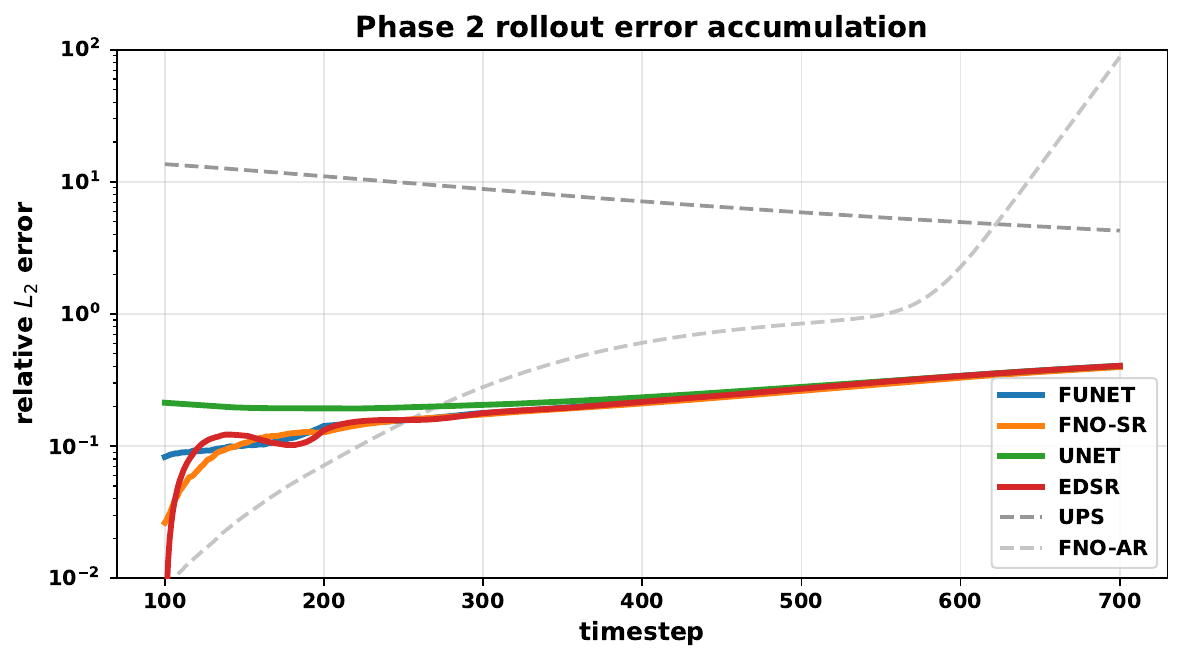}
            \vspace{3pt}
            \subcaption{Error accumulation of different models for an new test trajectory at a fixed downsampling scale factor of $4\times$.}
            \label{fig:err-acc-a}
        \end{minipage}
        &
        \begin{minipage}[t]{0.48\linewidth}
            \centering
            \includegraphics[width=\linewidth]{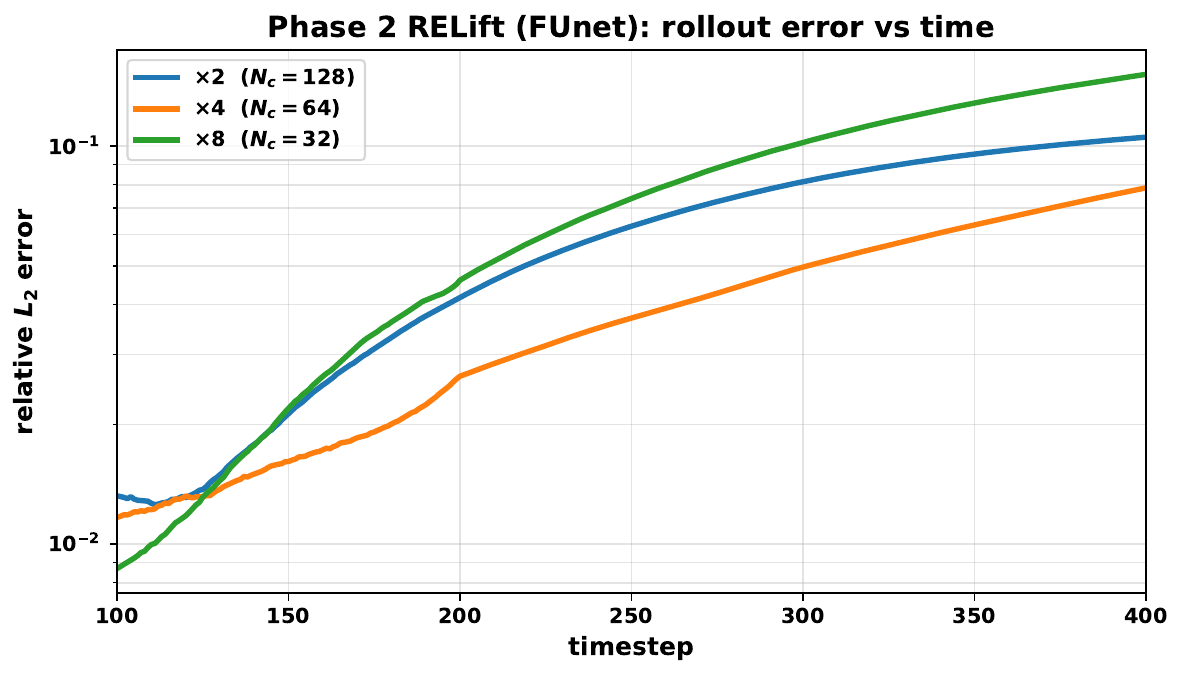}
            \vspace{3pt}
            \subcaption{FUnet error accumulation by different downsampling scale factor.}
            \label{fig:err-acc-b}
        \end{minipage}
    \end{tabular}
    \caption{\justifying
\textbf{Phase~2 error accumulation for NSE.}  
\textbf{(Left)} Relative \(L_2\) error for a representative unseen test trajectory.  
At a fixed \(4\times\) downsampling factor, the RELift Phase~2 predictions (solid, colored curves) are compared against the baselines (dashed, gray curves).  
Beyond timestep 400, the autoregressive FNO baseline rapidly diverges, whereas the RELift predictions exhibit substantially slower error growth.  
\textbf{(Right)} Phase~2 rollout error of RELift with the FUnet architecture, illustrating how prediction accuracy scales with the coarse input resolution \(N_c \in \{32,64,128\}\) while keeping the fine grid fixed at \(256\times256\).}
    \label{fig:err-acc}
\end{figure*}

\subsection{Future dynamics prediction (Phase 2 test results)}

For the prediction of future dynamics in each PDE example, we use the effective propagator (\ref{eqn:opcomp}) trained in a short time-window $[0,1]$ to approximate the fine-scale solution. Figure~\ref{fig:phase2-fields} presents the dynamics prediction results for a test run with a new random initial condition. Specifically, for the heat, wave, and NS equation, we extrapolate the dynamics from time $t=1$ to $t=2,10,7$, respectively.  Table~\ref{tab:phase2_toy} summarizes the averaged statistics over the 20 new test trajectories for each RELift model. Taken together, the results in Table~\ref{tab:phase2_toy} and Figure~\ref{fig:phase2-fields} lead to the following conclusions:

\paragraph{\textbf{ Numerical advantages for nonlinear PDE.}}
From Figure~\ref{fig:phase2-fields} and the metric values in Table~\ref{tab:phase2_toy}, RELift shows clear numerical advantages over baseline methods in predicting the dynamics of the NSE, a nonlinear PDE. This trend is consistent with  Phase~1 reconstruction results. The main reason is the rapid accumulation of prediction errors that is especially present in nonlinear systems\cite{chattopadhyay2023long}. As shown in the left panel of Figure~\ref{fig:err-acc}, for an unseen test trajectory, the FNO-AR method quickly accumulates approximation errors during the dynamics extrapolation phase, leading to large deviations from the true trajectory after only a short time. This effect is also evident in Figure~\ref{fig:phase2-fields}, where at the final time $t=7$, the predicted solution blows up. 


Across the three PDE examples, the Phase~2 results exhibit a clear and consistent trend: the numerical advantage of RELift increases with the spectral and dynamical complexity of the underlying system.  
For the heat equation, the fine-scale solution is dominated by low-frequency content, and the autoregressive FNO baseline already performs well; in this regime, RELift models achieve accuracy comparable to FNO-AR. For the linear wave equation, however, the presence of persistent medium-frequency structure makes the prediction task more challenging. Here, coupling the coarse propagator with the learned super-resolution operator in Phase~2 yields a more pronounced improvement over the autoregressive and upsampling baselines. 
Finally, for the incompressible NSE, whose nonlinear dynamics contain high-frequency vorticity features, RELift provides the largest gains.  
The learned super-resolution operator reconstructs fine-scale structures that are absent on the coarse grid, enabling the effective propagator to outperform autoregressive baselines substantially in long-time forecasting.

To test the \text{scale robustness} of the Phase~2 results, the right panel of Figure~\ref{fig:err-acc} repeats the rollout using RELift with the FUnet across three coarse-grid resolutions \(N_c\in\{32,64,128\}\) while fixing the fine grid at \(256\). Together with the left panel of Figure~\ref{fig:err-acc}, this indicates that RELift’s effective propagator is numerically more accurate and stable than autoregressive (AR) baselines for the long-time predictions of the dynamics of nonlinear NSE. Moreover, it is more resilient to moderate changes in spatial sampling of the input, which is crucial for practical usage of the developed approach.

\paragraph{\textbf{Composition propagator vs. equation-free propagator.}}  
Comparing the numerical results of FNO-SR and FNO-AR, we can assess the effectiveness of learning an effective propagator via the composition operator~\eqref{eqn:opcomp}. FNO-AR uses the FNO architecture from its original formulation~\cite{li2021fourier}  to learn an equation-free fine propagator; in contrast, FNO-SR learns an effective fine propagator that consists of the composition of a projection, a coarse propagator, and a learned super-resolution operator. 

From Figure~\ref{fig:phase2-fields} and the metrics in Table~\ref{tab:phase2_toy}, we observe that FNO-AR is more accurate for the heat equation, whereas FNO-SR performs better for more complex systems such as the non-dissipative wave equation and the nonlinear NSE, across all metrics.  It is important to note that our extrapolation time window (up to 700 timesteps for NS) are substantially longer than those typically reported in the original FNO literature\cite{li2021fourier}, where rollouts of 10 to 50 steps are  common. At these extended time windows, autoregressive error accumulation becomes significantly more pronounced, especially for the nonlinear case, as each step compounds the approximation error from all previous steps, despite adequate hyperparameter tuning (see Figure~\ref{fig:err-acc} and analysis of Section~\ref{sec:phase2-error}).  Furthermore, we use a nonrestrictive spectral cutoff that preserves high-frequency content (Appendix~\ref{app:data_generation}), meaning that the inherent low-frequency bias of FNO-AR becomes especially prevalent at long forecasting time windows~\cite{khodakarami2025mitigating, qin2024toward, kalimuthu2025loglo}. In contrast, RELift's composition structure, where temporal evolution is handled by the coarse solver and the neural operator provides only spatial correction, exhibits the slow-growth additive regime that remains stable over these long time windows. In such cases, the composition operator~\eqref{eqn:opcomp} serves as a better ansatz for the fine-scale dynamics solver. This trend is consistent with the numerical performance of the other RELift models discussed above. 

This behavior is further supported by the viscosity-shift NS rollout study in Figure~\ref{fig:ns-relL2-time-multi-nu}, where all models are trained only at $\nu=10^{-4}$ from the training window $[0,1]$ and extrapolated to physical time $T=51$ (5000 timesteps). We see that RELift predictions remain bounded and stable for $\nu=10^{-5}$ and $\nu=10^{-3}$, while the autoregressive FNO baseline becomes unstable much earlier under the same long-horizon setting. The small transient spikes in error observed at early times in the rollout are also consistent with the chosen forcing term. After an initial adjustment period, the forcing drives the dynamics toward a smoother long-time regime, thereby reducing sensitivity to accumulated local prediction errors, yielding stable predictions, as also evidenced by the decrease in error of upsampling in time. Ultimately, this provides additional evidence that the composition-based effective propagator is more robust than a fully equation-free autoregressive ansatz, even outside the training viscosity regime; full details are given in Section~\ref{sec:ablation}.

\paragraph{\textbf{ SSIM and spectral MSE for perceptual and spectral fidelity.}}
Relative \(L_2\) error tends to conflate large-scale shape alignment with fine-scale features, so model predictions that  \text{visually} differ can have similar \(L_2\) error. Thus, we employ two additional metrics, SpecMSE and SSIM.  SSIM compares the structural difference between the predicted flow and the ground truth; hence it is able to track the human-level perceived fidelity of the predictions much more closely. SpecMSE compares the error in the energy spectra of the predicted PDE solutions, a physically relevant metric. Thus, when relative $L_2$ errors are similar (e.g., for the wave equation), SSIM and SpecMSE offer more discriminative, task-relevant diagnostics to identify the best performing model. This is particularly pronounced in the Phase~2 wave equation results of Table \ref{tab:phase2_toy}.

\paragraph{\textbf{Runtime and storage benefits.}} 
We further test the generalizability of the learned super-resolution operator $\N_{\theta}$ for dynamics prediction. 
For the NSE example, we vary a newly introduced coarse timescale factor, denoted by $s\geq 1$, and use it in the Phase~2 dynamics prediction scheme as: $u^{n+s}\approx \hat{\mathcal P}_f(s\Delta t;\theta)\,u^n$.  The obtained numerical solution will be benchmarked with the fine-grid one on this coarse time grid. Note that changing $s$ only alters the coarse integrator’s step size and thus the temporal spacing from $\Delta t$ to $s\Delta t$ at test time, and the $\N_{
\theta},P$ operator will be fixed. Concretely, the fine baseline advances with $600$ steps of size $\Delta t$, whereas RELift advances with $600/s$ steps of size $s\Delta t$. For a fixed end time of the simulation, we compare the total wall-clock time and the numerical error.

\begin{figure*}[t]
  \centering
  \begin{minipage}[b]{0.32\textwidth}
    \includegraphics[width=\linewidth]{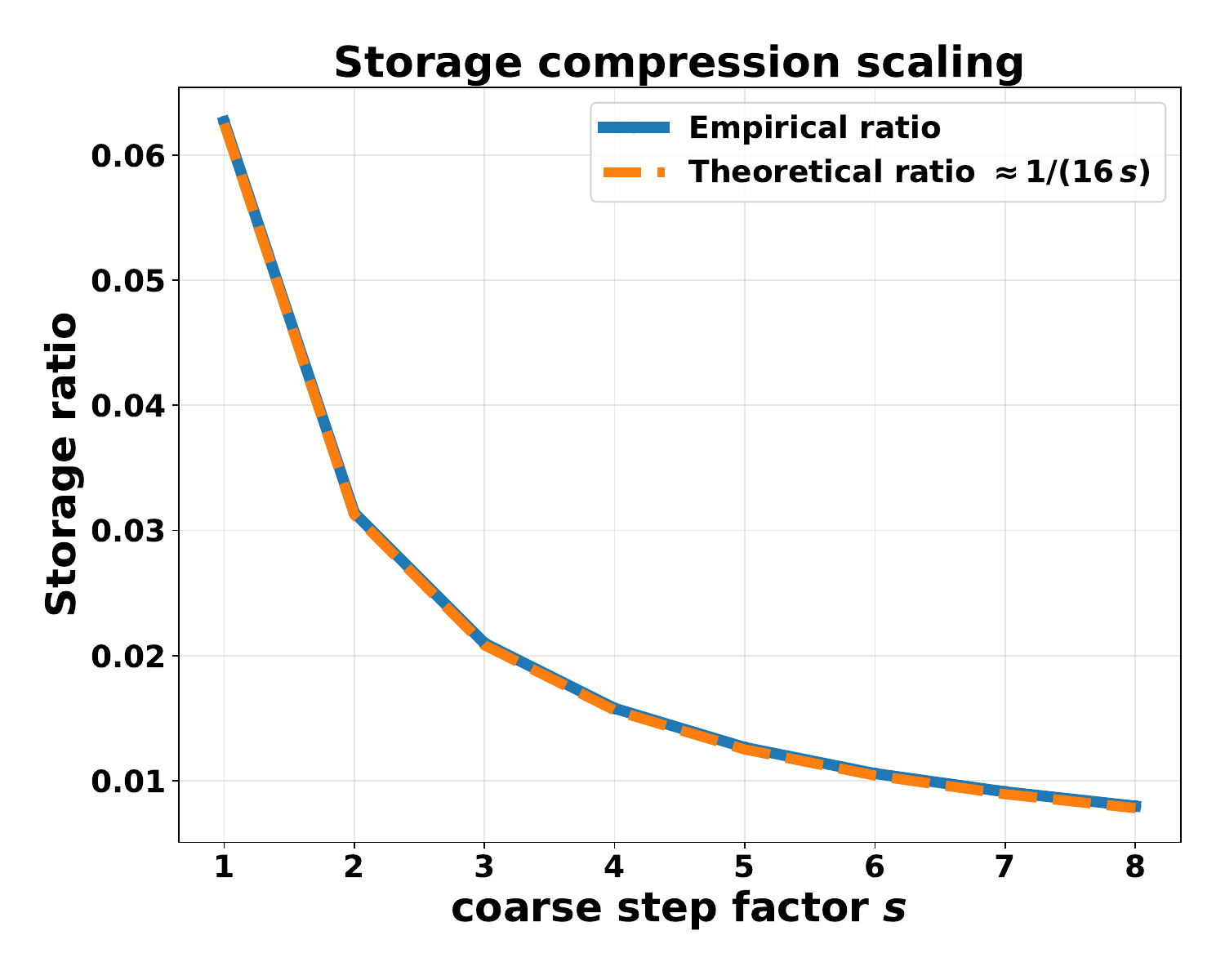}
  \end{minipage}\hfill
  \begin{minipage}[b]{0.32\textwidth}
    \includegraphics[width=\linewidth]{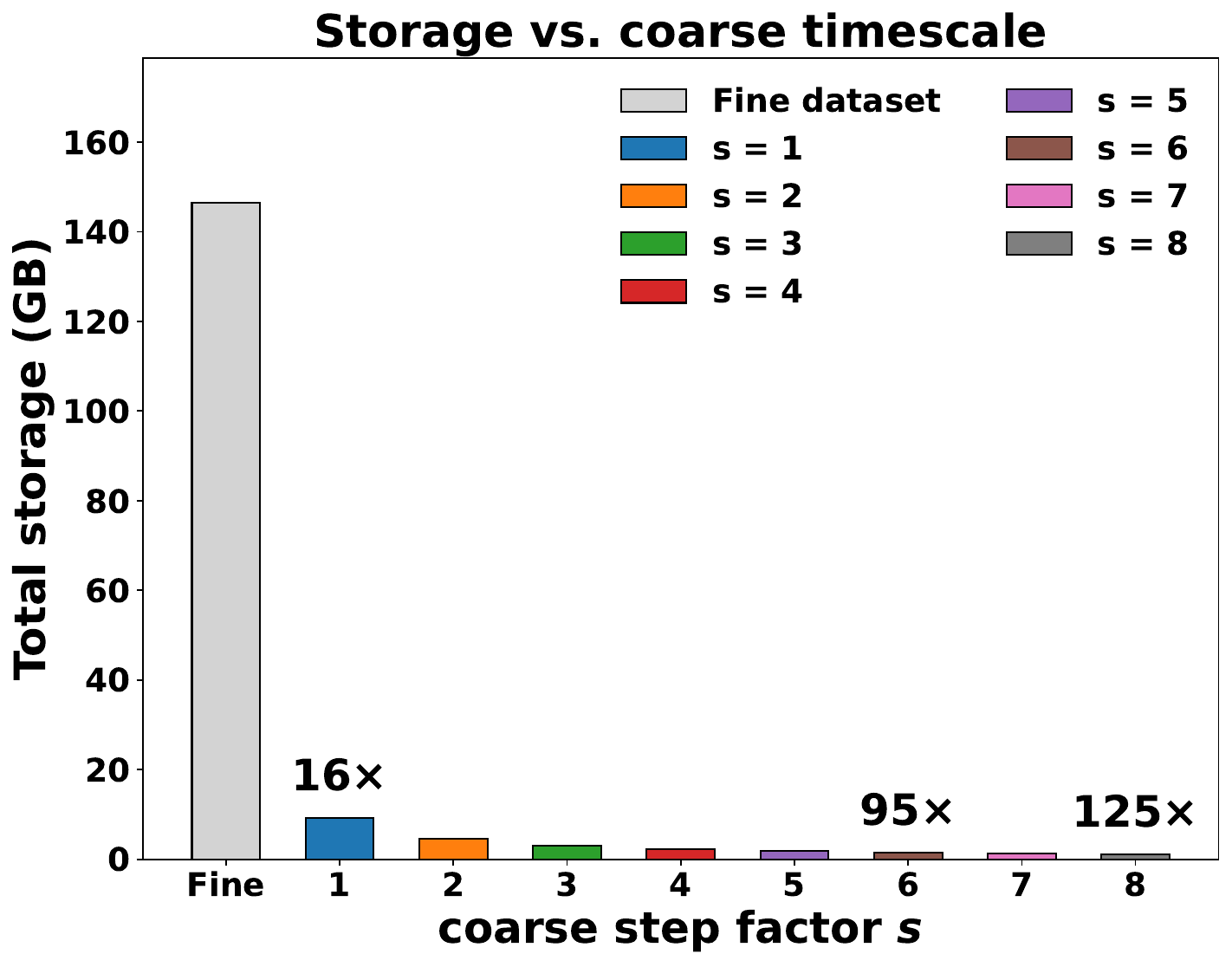}
  \end{minipage}\hfill
  \begin{minipage}[b]{0.32\textwidth}
    \includegraphics[width=\linewidth]{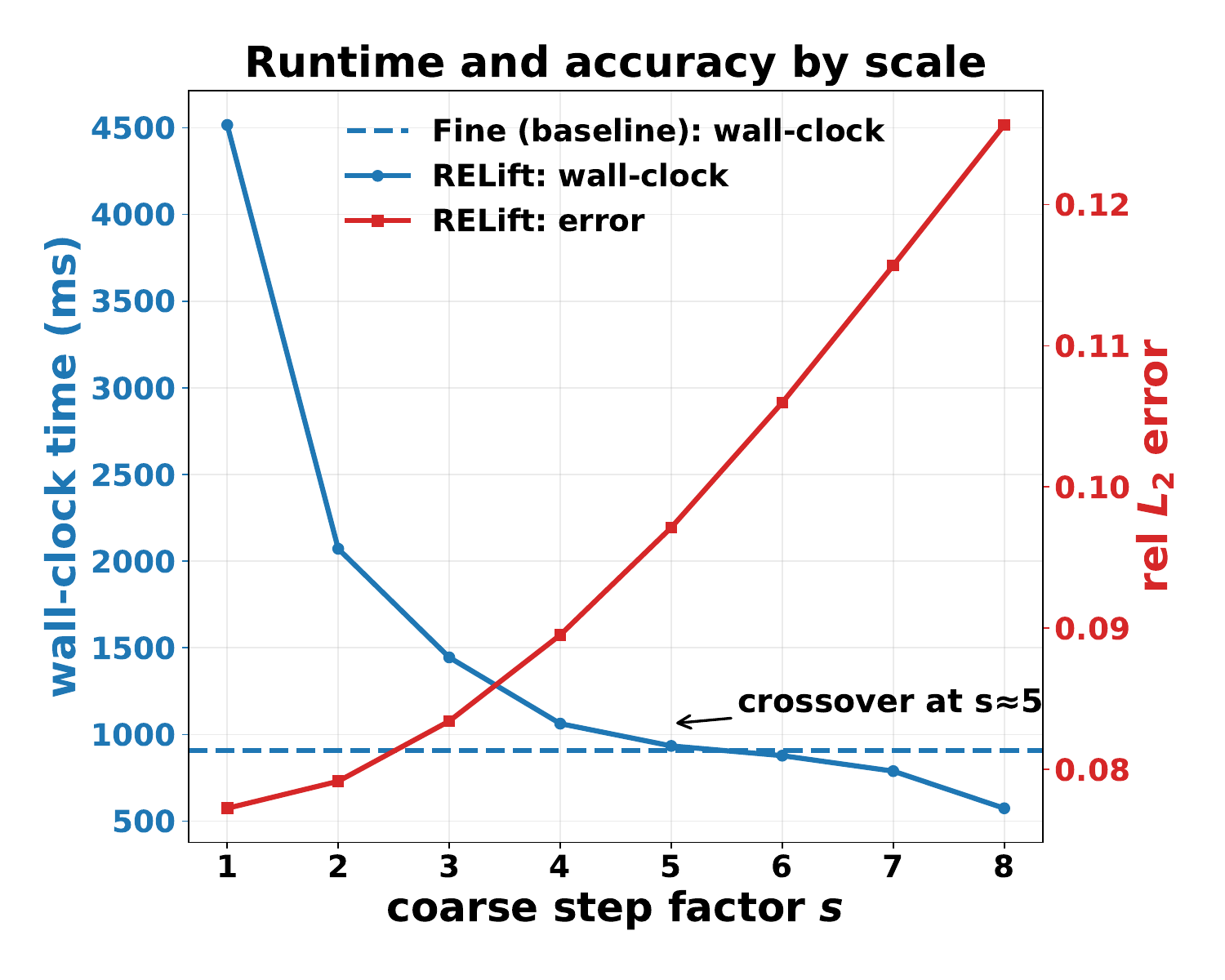}
  \end{minipage}
  \caption{\justifying
\textbf{Storage and runtime benefits for NSE.}  
All results use the FUnet architecture within RELift.  
\textbf{(Left)} Storage compression as a function of the coarse time-scale factor \(s\); the empirical trend closely matches the theoretical \(1/(16s)\) scaling implied by \(4\times\) spatial downsampling. 
\textbf{(Middle)} Total storage versus \(s\) (same trained model and number of trajectories).  
Values above each bar indicate the corresponding computational savings factor.  
\textbf{(Right)} Efficiency–accuracy tradeoff: total wall-clock runtime (left axis) and time-averaged relative \(L_2\) error (right axis).
}
  \label{fig:storage_runtime_benefits}
\end{figure*}

As shown in Figure~\ref{fig:storage_runtime_benefits}, RELift overtakes the fine solver at $s\approx 5$ and yields concrete end-to-end speedups of \text{1.06}$\times$ ($s{=}5$), \text{1.27}$\times$ ($s{=}6$), \text{1.50}$\times$ ($s{=}7$), and \text{1.71}$\times$ ($s{=}8$). Note that all results here use the FUnet as the neural network model. As expected, the accuracy trade–off increases exponentially with $s$, providing a transparent compute–fidelity knob. From the test, we see that taking larger coarse-grid time steps $\Delta t$ can reduce total runtime at acceptable error. 


To address the dataset level storage tradeoff, we compare storing \(1000\) full fine trajectories on disk against storing a single RELift trained FUnet plus \(1000\) coarse trajectories on disk at the same physical time window, shown in the left and center panels of Figure~\ref{fig:storage_runtime_benefits}. Even at \(s{=}1\), storing just the weights of the model and the coarse dataset option is already \({\approx} \, 16\times\) smaller than storing the fine dataset (9.2 GB vs 145.6 GB). As \(s\) increases, total storage decreases further, yielding \({\approx}\,79\times\) at \(s{=}5\) and \({\approx}\,125\times\) at \(s{=}8\). Practically, this enables saving orders of magnitude less data while retaining the ability to reconstruct fine-resolution fields on demand at evaluation time.

In summary, the numerical experiments demonstrate the effectiveness of the RELift framework in PDE dynamics reconstruction and prediction. When integrated with modern super-resolution architectures, RELift enables the reconstruction of fine-scale numerical solutions of PDEs from the coarse-scale solutions. By design, the composition operator serves as a natural ansatz for constructing an effective propagator to predict future fine-scale dynamics. In the NSE example, we find that RELift substantially outperforms baseline approaches such as FNO-AR and naive upsampling.

\section{Error Analysis}\label{sec:phase2-error}

In this section, we present an error analysis of the RELift-learned effective fine-grid propagator $\hat{\mathcal{P}}_f(\Delta t)$ from Phase~2 for predicting the numerical solutions of PDEs. In standard autoregressive forecasting, one learns a single-step map
$\Phi_\theta: \hat{u}^n \mapsto \hat{u}^{n+1}$ directly, so that both spatial
refinement and temporal evolution must be inferred by the network.
By contrast, our construction uses the compositional operator
$\hat{\mathcal{P}}_f(\Delta t; \theta)
 = \mathcal{N}_\theta \circ \mathcal{P}_c(\Delta t) \circ P$,
which explicitly decomposes the task into a \emph{spatial super-resolution step}
and a \emph{temporal dynamics propagation step}, with each stage introducing its
own approximation error.
Accordingly, our analysis focuses on two complementary aspects:

\begin{enumerate}[leftmargin=*]
\item \textbf{Deterministic (trajectory-wise) error bound.}  
For the numerical solution of a well-posed initial value problem (IVP) of a PDE with a prescribed initial condition selected \emph{from the training dataset}, we derive explicit bounds for the one-step \emph{local approximation error} and the $n$-step \emph{global accumulation error}.  
For each case, we systematically compare the contributions of the \emph{super-resolution approximation error} and the \emph{numerical integration error}, analyzing their relative magnitudes and influence on the overall error accumulation.

\item \textbf{Statistical (distributional) generalization bound.}  
In Appendix \ref{subsec:gen}, for the numerical solution of a well-posed IVP of a PDE with an initial condition $u_0$ \emph{randomly sampled} from a distribution $\mathcal{D}$, we derive a bound for the \emph{expected one-step error} in terms of its empirical expectation plus a standard\cite{10.5555/200548,10.5555/2621980} $O(m^{-1/2})$ capacity term, where $m$ denotes the number of training samples.  
This result quantifies the generalization ability of the RELift-learned effective fine-grid propagator $\hat{\mathcal{P}}_f(\Delta t)$ when applied to new initial conditions $u_0 \sim \mathcal{D}$.
\end{enumerate}

\subsection{Preliminaries }\label{subsec:assump}
In the following analysis, we compare the discrete numerical solutions of PDEs generated by two propagators: a sufficiently accurate fine-grid propagator $\mathcal{P}_f(\Delta t): \mathbb{R}^{N_f} \rightarrow \mathbb{R}^{N_f}$; and its RELift-learned approximation $\hat{\mathcal{P}}_f(\Delta t): \mathbb{R}^{N_f} \rightarrow \mathbb{R}^{N_f}$. This formulation allows us to bound the error entirely within a finite-dimensional setting using the $\ell_p$ norm. For linear PDEs, this connection can be formalized using classical consistency-stability arguments (e.g., Lax–Richtmyer-type results). For nonlinear PDEs, analogous conclusions typically require problem-dependent stability or Lipschitz/contractivity assumptions on the time-$\Delta t$ solution operator (or on the discrete propagator) to control the propagation of perturbations.

We define the fine-scale spatial discretization of the  exact PDE solution by $u_f^n=[u(x_j, t_n)]_{j=1}^{N_f} \in \mathbb{R}^{N_f}$, where $N_f$ is the total number of fine grid points. The corresponding coarse-scale representation is written as $u_c^n=[u(x_j, t_n)]_{j=1}^{N_c} \in \mathbb{R}^{N_c}$ with $N_c < N_f$. By definition, we have $Pu_f^n=u_c^n$. Unless otherwise specified, all norms considered below are averaged $\ell_1$ norms on $\mathbb{R}^{N_f}$ or $\mathbb{R}^{N_c}$.

To facilitate the error analysis, we introduce the following assumptions and auxiliary results :
\begin{enumerate}[label=(A\arabic*)]
\item \textbf{Coarse propagator regularity.}
For fixed $\Delta t=t_{n+1}-t_n$, the coarse-scale numerical integrator
$\mathcal{P}_c(\Delta t): \mathbb{R}^{N_c} \rightarrow \mathbb{R}^{N_c}$
is \emph{locally Lipschitz} on a bounded set that contains the relevant coarse-scale states.
Specifically, let $\mathcal{D}\subset \mathbb{R}^{N_c}$ denote a bounded set (e.g., a forward-invariant set
or a ball containing all coarse trajectories under consideration). Then, there exists a constant $L_c(\mathcal{D})>0$
such that for all $w_1,w_2\in\mathcal{D}$,
\begin{align}\label{eq:Lc}
\|\mathcal{P}_c(\Delta t) w_1 - \mathcal{P}_c(\Delta t) w_2\|
\le L_c(\mathcal{D})\,\|w_1-w_2\|.
\end{align}

\item \textbf{Neural lift regularity.} 
The {\em trained and fixed} neural lift operator  $\N_{\theta}:\R^{N_c}\rightarrow\R^{N_f}$ is also Lipschitz continuous and satisfies: 
\begin{equation}\label{eq:LN}
\bigl\|\mathcal N_\theta[w_1]-\mathcal N_\theta[w_2]\bigr\|
\le L_N\,\|w_1-w_2\|,
\end{equation}
for some \(L_N>0\) and all $w_1,w_2\in \R^{N_c}$.

\item \textbf{Error of the coarse numerical integrator.}  
Let the exact coarse-scale one-step propagator be denoted by 
\(\mathcal{P}_c^{\mathrm{exact}}(\Delta t): \R^{N_c} \to \R^{N_c}\).
Assume \(\mathcal{P}_c(\Delta t)\) is a \(q\)-th order time integrator for the coarse-grid ODE.
Then, for any fixed finite step horizon \(M\), there exists a constant \(C=C(M)>0\) such that
\begin{equation}\label{eq:coarse-consistency}
\bigl\| \mathcal{P}_c^{\mathrm{exact}}(\Delta t) u_c^n - \mathcal{P}_c(\Delta t) u_c^n \bigr\|
\le C(M)\, \Delta t^{\,q+1},
\end{equation}
where \( 0\le n \le M\). Here, \(C(M)\) depends on bounds of the derivatives of the coarse-grid ODE up to a chosen order, the numerical scheme, and the initial condition \(u_c^0\) (or equivalently, on the trajectory segment
\(\{u_c^n\}_{n=0}^M\)).
\end{enumerate}

\noindent\textbf{Remark.}
Although \(C(M)\) can in fact be arbitrarily large, for any prescribed tolerance \(\varepsilon>0\) and fixed \(M\),
the bound in~\eqref{eq:coarse-consistency} can, in turn, be made arbitrarily small by choosing
\[
\Delta t \le \left(\frac{\varepsilon}{C(M)}\right)^{\!1/(q+1)}.
\]

For (A1)–(A3), whether the coarse propagator has Lipschitz continuity (A1) can be determined from the
coarse-grid ODE system. In particular, as \(\Delta t \to 0\), the Lipschitz constant \(L_c\) can be related to
the Lipschitz constant of the underlying ODE vector field on the relevant bounded set.
Assumption (A2) holds since standard neural networks are compositions of affine maps and Lipschitz nonlinearities.
Finally, (A3) is the standard local truncation error estimate for a \(q\)-th order integrator; the constant
\(C(M)\) depends on the ODE derivative bounds, the scheme, and the initial condition \(u_c^0\); and the bound
can be made arbitrarily small by taking \(\Delta t\) sufficiently small for any fixed finite \(M\).

\begin{figure*}[thbp!]
  \centering
  \setlength{\abovecaptionskip}{6pt}
  \setlength{\belowcaptionskip}{0pt}
  \setlength{\floatsep}{6pt}
  \setlength{\textfloatsep}{8pt}

  \begin{subfigure}{0.32\linewidth}
    \centering
    \includegraphics[width=\linewidth]{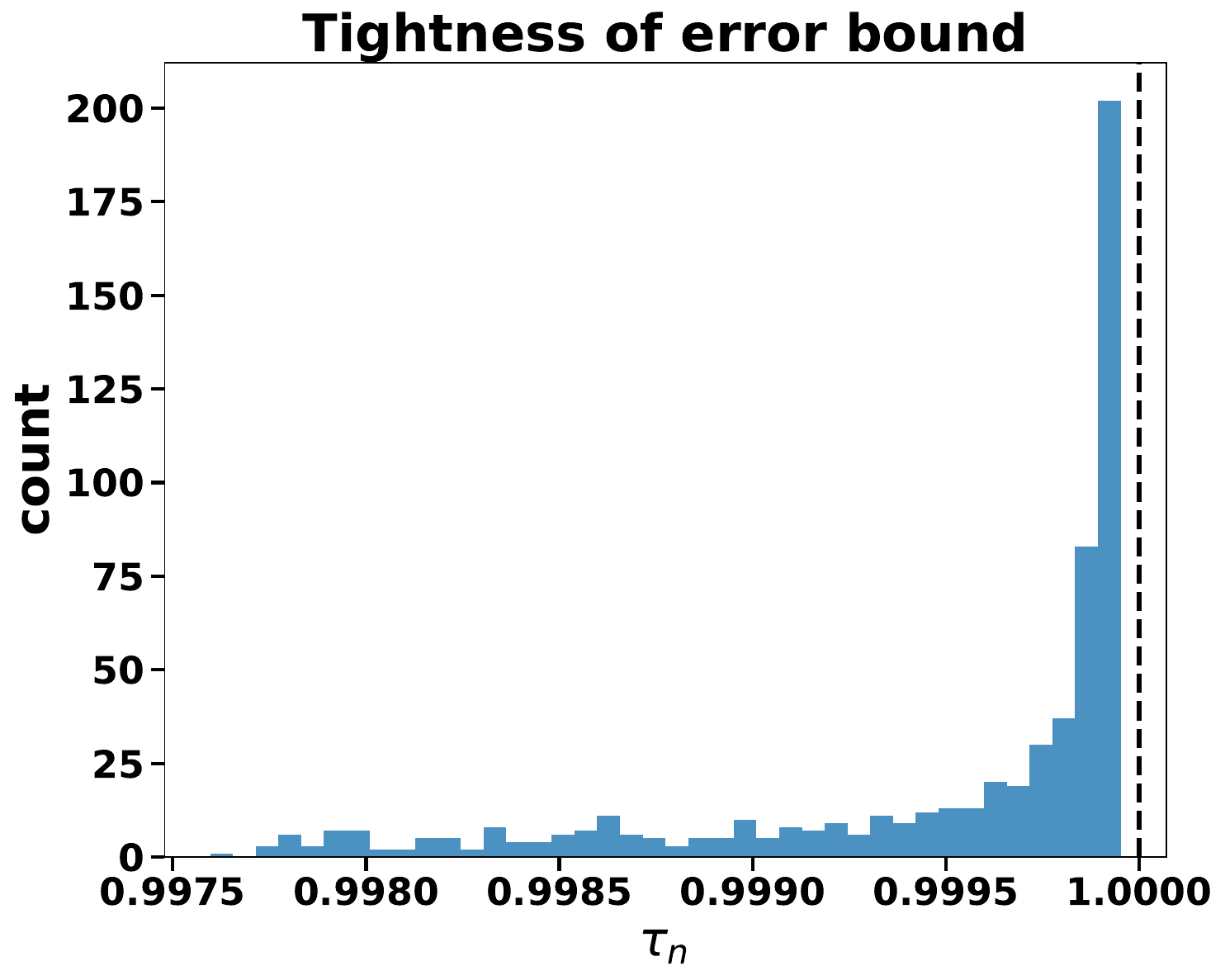}
    \subcaption{Histogram of $\tau_n=e_n/(r_n{+}b_n)$.}
    \label{fig:perstep-a}
  \end{subfigure}\hfill
  \begin{subfigure}{0.32\linewidth}
    \centering
    \includegraphics[width=\linewidth]{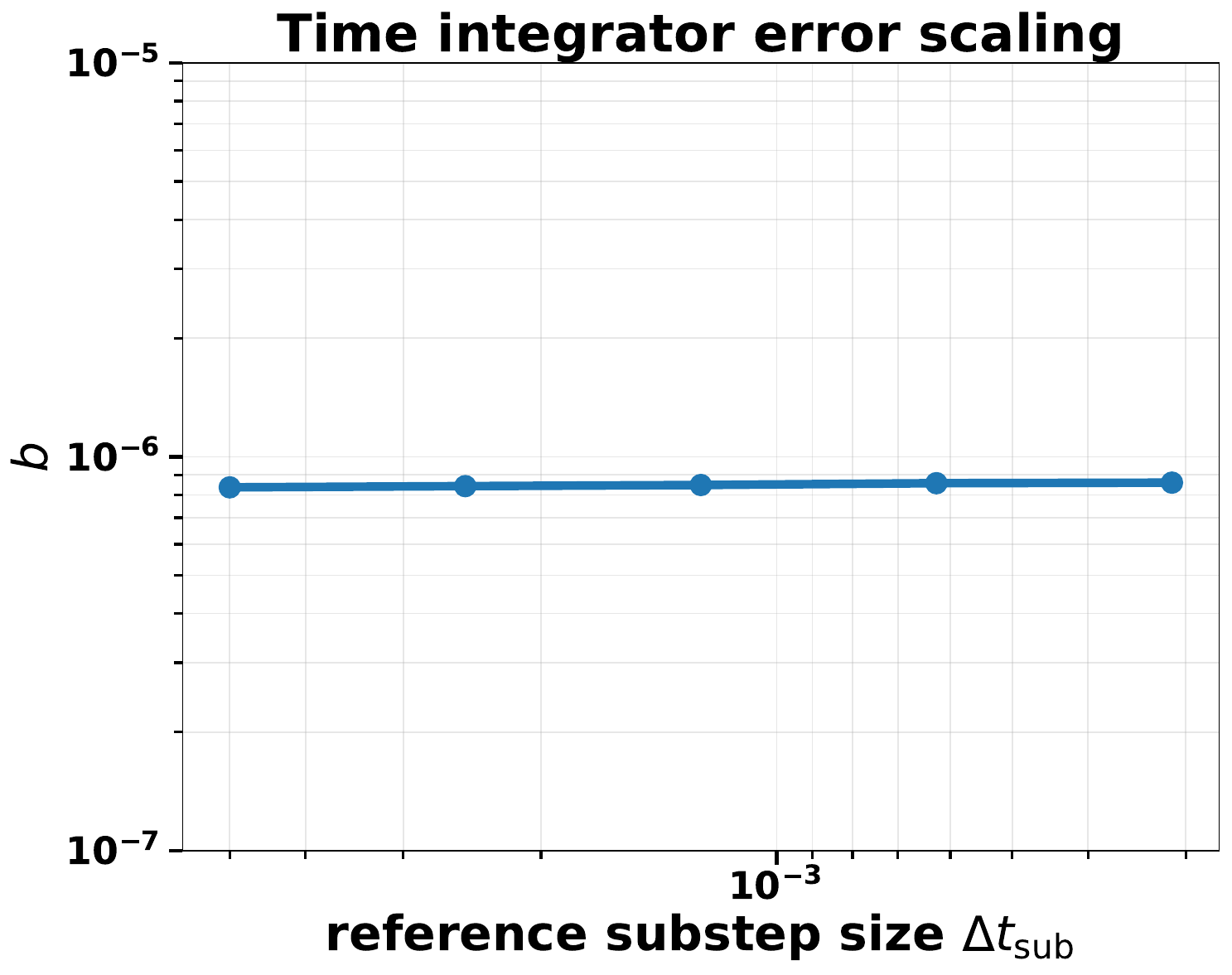}
    \subcaption{$b$ vs.\ substep size $\Delta t_{\mathrm{sub}}$.}
    \label{fig:perstep-b}
  \end{subfigure}\hfill
  \begin{subfigure}{0.32\linewidth}
    \centering
    \includegraphics[width=\linewidth]{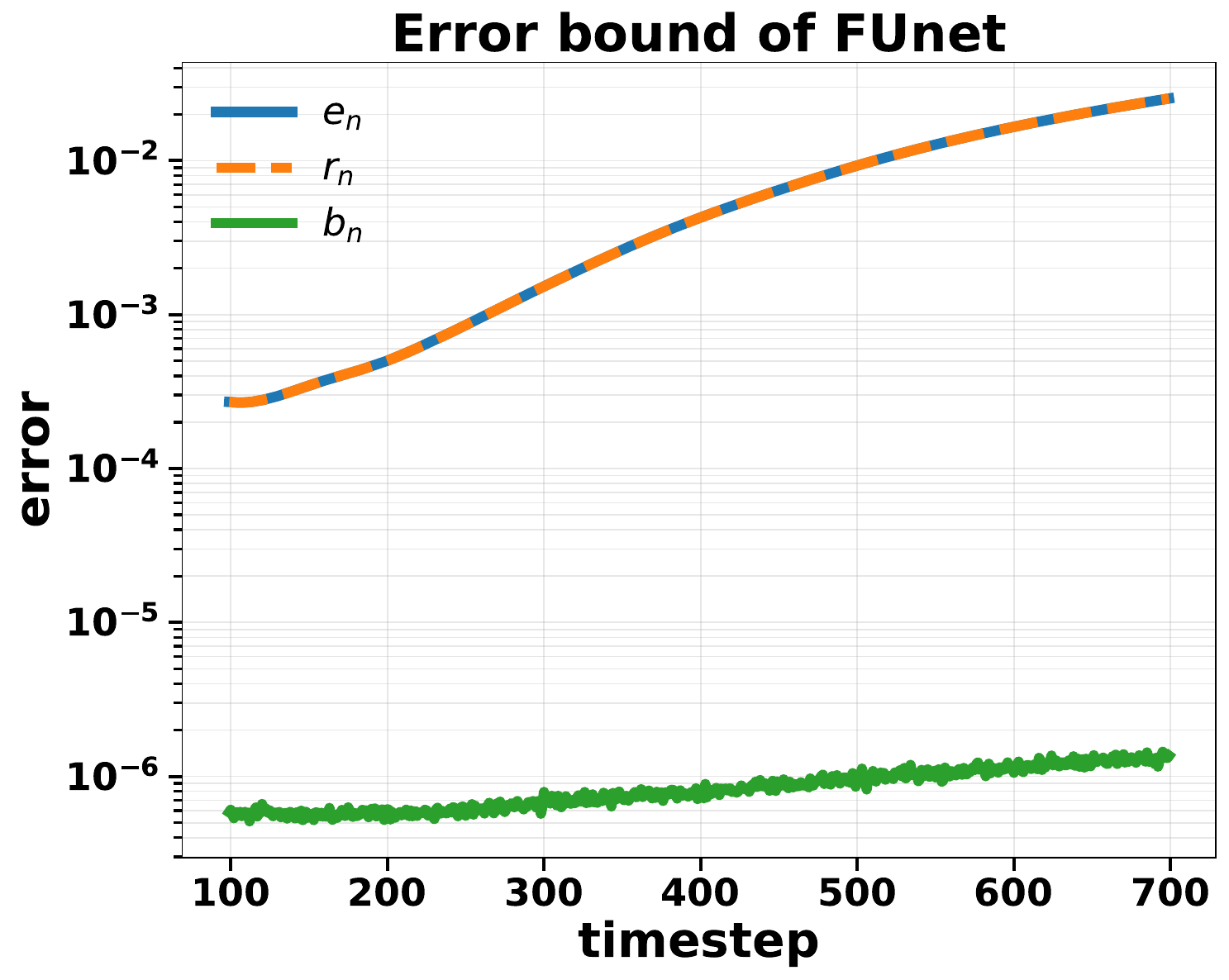}
    \subcaption{Local approximation error.}
    \label{fig:perstep-c}
  \end{subfigure}
\caption{\justifying
\textbf{Local approximation error on a test NSE trajectory.}  
\textbf{(Left)} Histogram of the ratio \(\tau_n = e_n / (r_n + b_n)\) for timesteps \(100 \le n \le 700\).  
\textbf{(Middle)} Mean numerical integration error as a function of timestep refinement.  
\textbf{(Right)} Comparison of \(e_n\), \(r_n\), and \(b_n\) along a long-time trajectory.
}
  \label{fig:perstep-figs}
\end{figure*}
\subsection{Theoretical local and global error bounds}\label{subsec:onestep}
We first consider the trajectory-wise local approximation error defined by:
\begin{equation}\label{eq:step-error-def}
\varepsilon_{\mathrm{step}}(u_f^n;\theta)
:=
\|\mathcal P_f(\Delta t) u_f^n
-  \hat\P_f(\Delta t;
\theta) u_f^n\| ,
\end{equation}
from which it follows that
\begin{equation}\label{eq:step-error-def}
\begin{aligned}
\varepsilon_{\mathrm{step}}(u_f^n; \theta)
&\le
\underbrace{\bigl\| \bigl( \mathcal{P}_f(\Delta t)
- \mathcal{N}_{\theta} \circ \mathcal{P}_c^{\mathrm{exact}}(\Delta t) \circ P \bigr) u_f^n \bigr\|}_{\mathsf{A}(u_f^n; \theta)}
\\
&\quad + \underbrace{
\bigl\|
\underbrace{\mathcal{N}_{\theta}\bigl( \mathcal{P}_c^{\mathrm{exact}}(\Delta t) P u_f^n \bigr)}_{\mathsf{B}_1(u_f^n;\theta)}
-
\underbrace{\mathcal{N}_{\theta}\bigl( \mathcal{P}_c(\Delta t) P u_f^n \bigr)}_{\mathsf{B}_2(u_f^n;\theta)}
\bigr\|
}_{\mathsf{B}(u_f^n; \theta)}.
\end{aligned}
\end{equation}

Using (A1)–(A3), we can derive a uniform bound for $\mathsf{B}(u_f^n)$ as $
\mathsf{B}(u_f^n) \le C L_N \, \Delta t^{\,q+1}
$, where the constant $C$ is independent of $u_c^n$.  
Since $\mathcal{P}_c^{\mathrm{exact}}(\Delta t)$ represents the exact one-step propagator of the coarse dynamics, the two terms can be interpreted as follows:  
$\mathsf{A}(u_f^n; \theta)$ corresponds to the \emph{super-resolution approximation error}; while $\mathsf{B}(u_f^n)$ represents the \emph{numerical integration error}. From the uniform upper bound $\mathsf{B}(u_f^n) \le C L_N \, \Delta t^{\,q+1}$, we observe that as $\Delta t \to 0$, the contribution of $\mathsf{B}(u_f^n)$ to the total error becomes asymptotically negligible and hence the super-resolution error dominates. It is important to note that this directly justifies the objective function ~\eqref{eqn:L_propagator} we chose for Phase 2 training of the neural network, where $\mathsf{A}(u_f^n; \theta)$ is minimized across the training window.

\paragraph{\textbf{Numerical evaluation of \eqref{eq:step-error-def}} } Although exact evaluation of the two error terms is generally intractable, they can be well approximated using a highly accurate reference solution that approximates the action of $\mathcal{P}_c^{\mathrm{exact}}(\Delta t)$:
\[
z_{\mathrm{ref}}(u_f^n)
:= \bigl[\mathcal{P}_c\!\bigl(\tfrac{\Delta t}{\nu}\bigr)\bigr]^{\nu} P u_f^{n}
\approx \mathcal{P}_c^{\mathrm{exact}}(\Delta t) P u_f^n,
\]
where $\nu \ge 2$ is a temporal refinement factor.  Explicitly, for our purposes, we treat this as a sufficiently accurate surrogate for the exact propagator. At each time step $n$, the total error and its two contributing components can be tracked by their corresponding numerical approximations:
\[
\varepsilon_{\mathrm{step}}(u_f^n; \theta) \approx e_n, 
\quad \mathsf{A}(u_f^n; \theta) \approx r_n, 
\quad \mathsf{B}(u_f^n) \approx b_n,
\]
where
\begin{equation}\label{eqn:e_n,r_n,b_n}
\begin{aligned}
e_n 
&:= \bigl\| \mathcal{P}_f(\Delta t)\tilde{u}_f^{n}
      - \mathcal{N}_{\theta}\mathcal{P}_c(\Delta t)P\tilde{u}_f^n \bigr\|, \\
r_n 
&:= \bigl\| \mathcal{P}_f(\Delta t)\tilde{u}_f^{n}
      - \mathcal{N}_{\theta}\!\bigl[z_{\mathrm{ref}}(\tilde{u}_f^n)\bigr] \bigr\|, \\
b_n
&:= \bigl\| \mathcal{N}_{\theta}\!\bigl[z_{\mathrm{ref}}(\tilde{u}_f^n)\bigr]
      - \mathcal{N}_{\theta}\bigl[\mathcal{P}_c(\Delta t)P\tilde{u}_f^n\bigr] \bigr\|.
\end{aligned}
\end{equation}
Here $\tilde{u}_f^n := [\mathcal{P}_f(\Delta t)]^n u_f^0 \approx u_f^n$ denotes the approximate fine-scale numerical solution of the PDE.  
This approximation is valid when $\mathcal{P}_f(\Delta t)$ provides an accurate representation of the exact propagator. All the quantities above are computed using the $\ell_1$-norm on $\mathbb{R}^{N_f}$, i.e., evaluated on the fine grid.

Next, we consider the global approximation error of the effective fine-grid propagator $\hat{\mathcal{P}}_f(\Delta t)$.  
For a neural operator $\Phi_\theta: \hat{u}^n \mapsto \hat{u}^{n+1}$ that learns the discrete-time dynamics of a PDE (e.g., an FNO-type model), classical autoregressive (AR) analysis~\cite{mccabe2023stability} provides a global error bound in terms of its Lipschitz constant 
\(L_\theta := \mathrm{Lip}(\Phi_\theta)\):
\begin{equation}\label{eq:ar-kstep}
\|u_f^{n} - \hat{u}_f^{n}\|
\lesssim
\sum_{j=0}^{n-1} (L_\theta)^{\,n-1-j}\, e_j.
\end{equation}
For our composite effective map 
\(\hat{\mathcal{P}}_f := \mathcal{N}_\theta \circ \mathcal{P}_c \circ P\),
the submultiplicativity of the Lipschitz constant yields
\begin{align}\label{lip_P_f}
\mathrm{Lip}(\hat{\mathcal{P}}_f)
\le L_N\,L_c\,C_P,
\end{align}
where $L_N$ and $L_c$ are the Lipschitz constants defined in~\eqref{eq:Lc}–\eqref{eq:LN}, and $C_P$ denotes the induced operator norm of the projection operator $P$.  
Combining~\eqref{eq:ar-kstep}–\eqref{lip_P_f}, we obtain the following global approximation error bound:
\begin{equation}\label{eq:global_err}
\begin{aligned}
\|u_f^{n} - \hat{u}_f^{n}\|
&\lesssim
\sum_{j=0}^{n-1} (L_N L_c C_P)^{\,n-1-j} \bigl(r_j + b_j\bigr)\\
&+ \mathcal{O}(\Delta t^{\,p+1} + h_f^{\,l}),
\end{aligned}
\end{equation}
where $u_f^{n}$ denotes the exact solution of the PDE evaluated at time $t_n = n\Delta t$ on the fine grid, and 
$\hat{u}_f^{n} = [\hat{\mathcal{P}}_f(\Delta t)]^n u_f^0$ denotes the approximate solution obtained via the effective fine-grid propagator $\hat{\mathcal{P}}_f(\Delta t)$. The additional term $\mathcal{O}(\Delta t^{\,p+1} + h_f^{\,l})$ represents the discretization error associated with a $p$-th order in time and $l$-th order in space fine-scale finite difference scheme. For a convergent fine-scale numerical scheme, this error can be made arbitrarily small as $\Delta t, h_f \to 0$.  Consequently, the dominant contribution to the global approximation error should also arise from the accumulated super-resolution error $r_j, 1\leq j\leq n$. 

\paragraph{\textbf{Numerical evaluation of \eqref{eq:global_err}} } Since the one-step error bound can be approximated by $r_n$ and $b_n$ via their definitions in~\eqref{eqn:e_n,r_n,b_n}, numerically approximating the global error bound~\eqref{eq:global_err} reduces to estimating the Lipschitz constants $L_c$, $L_N$, and $C_p$. By definition, this can be done by  randomly sampling inputs $x, y$ within their respective domains and evaluating:
\begin{equation}\label{eqn:proxies}
\begin{aligned}
L_c &= \sup_{x, y \in \mathbb{R}^{N_c}} 
\frac{\| \mathcal{P}_c(\Delta t)x - \mathcal{P}_c(\Delta t)y \|_{\ell_1(\mathbb{R}^{N_c})}}
     {\|x - y\|_{\ell_1(\mathbb{R}^{N_c})}}, \\[4pt]
L_N &= \sup_{x, y \in \mathbb{R}^{N_c}} 
\frac{\| \mathcal{N}_{\theta}x - \mathcal{N}_{\theta}y \|_{\ell_1(\mathbb{R}^{N_f})}}
     {\|x - y\|_{\ell_1(\mathbb{R}^{N_c})}}, \\[4pt]
C_p &= \sup_{x \in \mathbb{R}^{N_f}} 
\frac{\| P x \|_{\ell_1(\mathbb{R}^{N_c})}}
     {\|x\|_{\ell_1(\mathbb{R}^{N_f})}}.
\end{aligned}
\end{equation}
For AR models, the Lipschitz constant
\(
L_\theta := \mathrm{Lip}(\Phi_\theta)
\)
associated with the single-step map
\(
\Phi_\theta: \hat{u}^n \mapsto \hat{u}^{n+1}
\)
can be estimated in a similar manner. Combining these two, the error bound analysis can help to explain the observed numerical stability and instability behaviors of these two operator-learning strategies.  A central challenge of the autoregressive operator-learning approach is that small one-step approximation errors can rapidly accumulate over time, leading to drift and eventual instability in the predicted trajectories~\cite{chattopadhyay2023long,mccabe2023stability}.  
Such instability is primarily induced by the large Lipschitz constant \(L_\theta = \mathrm{Lip}(\Phi_\theta)\) of the learned operator.  
In contrast, as we shall demonstrate in the next section, the effective Lipschitz constant \(L_N L_c C_P\) in the RELift framework is normally much smaller, resulting in a numerically more stable scheme for long-term dynamics predictions.


\subsection{Numerical evaluation of local approximation errors}\label{sec:phase2-implications}
We now use the NSE dataset as an illustrative example to numerically evaluate the local error bounds derived in Section~\ref{subsec:onestep} . 
Specifically, the one-step approximation error of $\hat{\mathcal{P}}_f(\Delta t)$ is evaluated at \emph{each} time step by comparing the predicted solution $\{\tilde u_f^n, \hat{\mathcal{P}}_f(\Delta t) \tilde u_f^n\}$ against the reference one $\{\tilde u_f^n, \mathcal{P}_f(\Delta t)\tilde u_f^n\}$. Recall that we assume the fine-scale solution well-approximates the exact solution, i.e., $
\tilde{u}_f^n := [\mathcal{P}_f(\Delta t)]^n u_f^0 \approx u_f^n$. The results presented here correspond to a representative test trajectory and the test for other trajectories yields similar results.  
Our main findings are summarized as follows:

\paragraph{\textbf{Bound validity and tightness.}}  
To examine the tightness of the local approximation error bound ~\eqref{eq:step-error-def},  
we evaluate the ratio between the measured error and the theoretical bound:
\(\tau_n =  e_n/(r_n + b_n)\)
over all $100 \le n \le 700$ extrapolated time steps, and we analyze the corresponding histogram of $r_n$. As shown in Figure~\ref{fig:perstep-a}, $\tau_n\approx 1$, which confirms that the error decomposition in~\eqref{eq:step-error-def} provides a tight bound for the local approximation~error.

\paragraph{\textbf{Negligible numerical integration error.}} 
In Figure~\ref{fig:perstep-b}, we plot the mean numerical integration error 
\(b = \frac{1}{N} \sum_{n} b_n\)
as a function of different refinements of the timestep: 
$\Delta t_{\mathrm{sub}} = \Delta t / n_{\mathrm{sub}}$ 
for $n_{\mathrm{sub}} \in \{2, 4, 8, 16, 32\}$, obtained using the same coarse-grid solver. Namely, for each $n_{\mathrm{sub}}$, we replace the coarse propagator $\mathcal P_c(\Delta t)$ inside $\hat{\mathcal P}_f(\Delta t)$ by the substepped composition $[\mathcal P_c(\Delta t/n_{\mathrm{sub}})]^{n_{\mathrm{sub}}}$.

If the numerical integration error were a dominant contributor to the one-step approximation error, then refining $\Delta t_{\mathrm{sub}}$ would noticeably decrease both $b$ and, in turn, the mean error 
$e = \frac{1}{N} \sum_{n} e_n$.  
For the selected test trajectory, however, we find that $b \approx 8.49 \times 10^{-7}$, which is four to five orders of magnitude smaller than the observed mean error 
$e \approx 7.34 \times 10^{-3}$. The growth of $e_n$, $r_n$, and $b_{n}$ in a long extrapolated trajectory is shown more clearly from Figure~\ref{fig:perstep-c}, which further confirms the claim we made by pure theoretical analysis that 
the numerical integration error introduced through the evaluation of $\mathcal{P}_c(\Delta t)$ in $\hat{\mathcal{P}}_f(\Delta t)$ is negligible compared to the super-resolution approximation error. We also repeated this experiment using Richardson extrapolation on the same coarse integrator and obtained nearly identical values of $b$, leading to the same conclusion that time-discretization error is negligible compared to the super-resolution approximation error.


It is important to emphasize that this observation does not imply  
(i) that using the coarse grid would resolve the underlying physics or  
(ii) that all numerical artifacts are negligible.  
Rather, it establishes that, for the chosen coarse discretization, the numerical integration error is insignificant relative to the super-resolution approximation error. The observed difference between the upsampled baseline and the RELift results for the NS equation can thus be attributed to the fact that the learned super-resolution operator $\mathcal{N}_{\theta}$ more accurately recovers the missing subgrid information and mitigates the projection–lifting mismatch.



\begin{figure}[t]
  \centering
  \includegraphics[width=\linewidth]{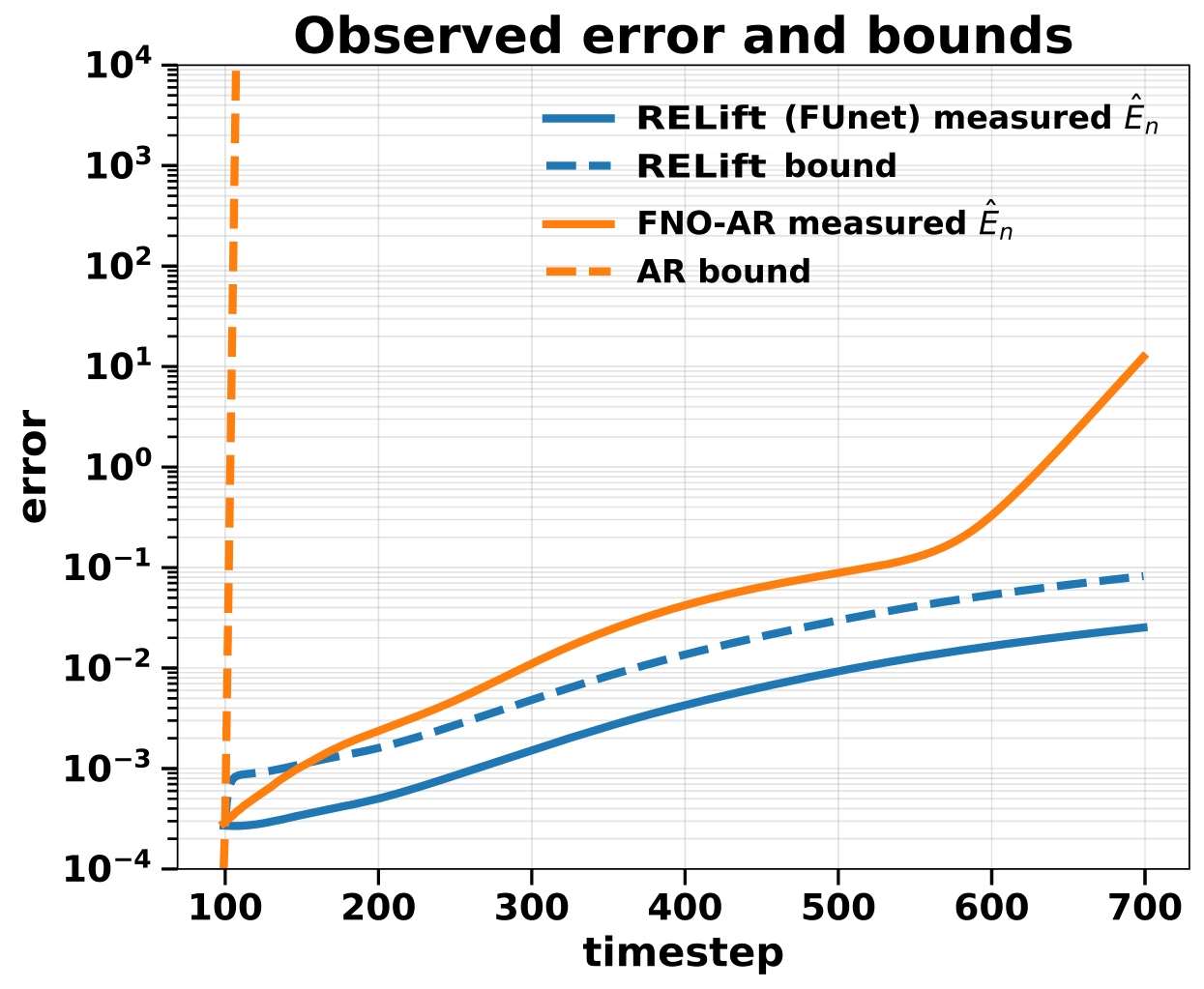}
  \caption{\justifying
\textbf{  Measured global approximation errors and their respective bounds for RELift and the AR baseline.}
 }
  \label{fig:ar-vs-srafte-envelope}
\end{figure}

\subsection{Numerical evaluation of global approximation errors }
The global approximation error can be evaluated by comparing the \emph{entire} predicted trajectory 
\(\{\hat{u}_f^n\}_{n=0}^{N}\) with the corresponding reference trajectory \(\{\tilde{u}_f^n\}_{n=0}^{N}\), 
where 
\(\hat{u}_f^n = [\hat{\mathcal{P}}_f(\Delta t)]^n u_f^0\) 
and 
\(\tilde{u}_f^n = [\mathcal{P}_f(\Delta t)]^n u_f^0\). 
Define
\[
\hat{E}_n = \|\tilde{u}_f^n - \hat{u}_f^n\|
\]
as the global approximation error at step \(n\).  
The error bound~\eqref{eq:global_err} then implies
\begin{align}\label{eqn:E_n_bound}
\hat{E}_n
\;\lesssim\;
\sum_{j=0}^{n-1}
(L_N L_c C_P)^{\,n-1-j}\, e_j,
\end{align}
where \(L_N L_c C_P\) represents the Lipschitz constant of the surrogate propagation operator, which can be approximated by Monte Carlo sampling and evaluating \eqref{eqn:proxies}. Similarly, for a neural operator 
\(\Phi_\theta : \hat{u}^n \mapsto \hat{u}^{n+1}\) 
with empirical Lipschitz constant
\(L_{\theta} := \mathrm{Lip}(\Phi_\theta)\),
the corresponding global approximation error bound~\eqref{eq:ar-kstep} can be expressed as 
\begin{equation}\label{eq:ar-bound}
\hat{E}_n
\;\lesssim\;
\sum_{j=0}^{n-1}
L_{\theta}^{\,n-1-j}\, e_j,
\end{equation}
which can be evaluated numerically in the same way. For a test trajectory of the NS equation, 
Figure~\ref{fig:ar-vs-srafte-envelope} compares the observed global approximation error 
with the theoretical error bounds~\eqref{eqn:E_n_bound}–\eqref{eq:ar-bound} 
for the RELift-learned effective fine propagator 
\(\hat{\mathcal{P}}_f(\Delta t)\) 
and the FNO-AR-learned solution operator 
\(\Phi_{\theta}\), respectively. 
We find that the composite Lipschitz constant 
\(\Lambda = L_N L_c C_P \approx 0.68\), 
so that the geometric weights 
\(\Lambda^{\,n-1-j}\) 
effectively attenuate the contribution of the one-step approximation errors. 
As a result, the theoretical bound~\eqref{eqn:E_n_bound} 
closely tracks the actual global error \(\hat{E}_n\). In contrast, for FNO-AR, the empirical Lipschitz constant is 
\(L_{\theta} \approx 
3.11\), 
which causes the corresponding error bound~\eqref{eq:ar-bound} 
to grow rapidly, reflecting the fast error accumulation observed 
in the FNO-AR predictions of the PDE solution. 

Hence, for standard autoregressive operator-learning methods such as FNO-AR, the global approximation error grows significantly faster 
than that of the RELift-learned dynamical propagator. 
This behavior is consistently observed in our numerical simulations, 
and the corresponding global error bound analysis 
provides a coherent theoretical explanation 
for the observed discrepancy in error growth rates.

\begin{figure*}[htbp]
    \centering

    \begin{subfigure}[t]{0.92\textwidth}
        \centering
        \includegraphics[width=\linewidth]{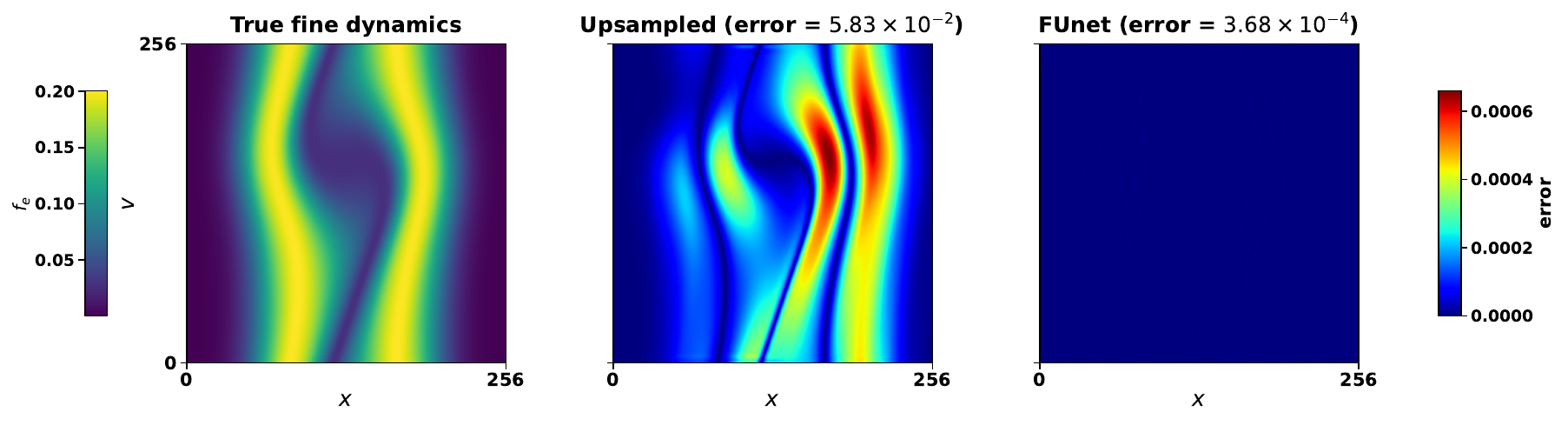}
        \caption{Two-stream instability electron species \(f_e\) at \(t=22.5\).}
        \label{fig:vp-fe-twostream-t05}
    \end{subfigure}

    \vspace{0.6em}

       \begin{subfigure}[t]{0.92\textwidth}
        \centering
        \includegraphics[width=\linewidth]{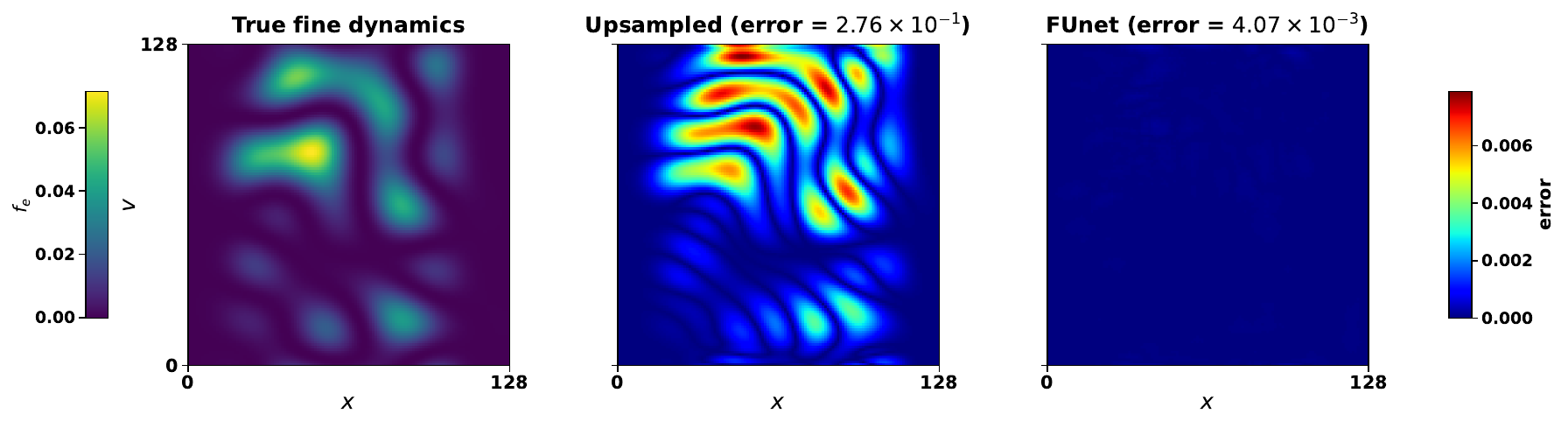}
        \caption{Randomized electrons \(f_e\) at \(t=0.5\).}
        \label{fig:vp-fe-rand-t05}
    \end{subfigure}

    \caption{\justifying{\textbf{Phase 1 test set results for 1D1V Vlasov--Poisson (two-stream instability and randomized datasets)  electron distribution \(f_e\) at $t=22.5$ and $t=0.5$, respectively}. We switch to error heatmaps to make spatially localized discrepancies explicit and to better separate methods whose predicted fields appear visually similar. The titles indicate the relative $L_2$ errors.}}
    \label{fig:vp-phase1-fe-twostream-random}
\end{figure*}

\section{Electron dynamics application}\label{sec:eldyapp}
Electron-scale phenomena are notoriously costly to simulate because often times the physics of interest only occur on very fine space-time grids \cite{burch2016electron, li2020real, ye2021reducing, xie2024coarse, xie2024ab,weinan2011principles,weinan2003heterognous}. Under uniform refinement in $d$ dimensions, halving the mesh width increases the number of degrees of freedom by approximately $2^d$, and it typically tightens stability constraints on the time step. Fast electron time scales and sub-Debye spatial structure further amplify this cost, making long-horizon studies and parameter sweeps impractical for conventional solvers. These regimes therefore provide a stringent and practically motivated test suite for RELift.  Accordingly, we next consider the classic Vlasov–Poisson system as a representative model for electron dynamics within a plasma.
 
\subsection{Vlasov--Poisson equation}

The Vlasov--Poisson (VP) equation governs the dynamics of collisionless plasmas. The term \textit{collisionless} means that the particle's  motion is primarily governed by electromagnetic fields or gravitational forces, rather than collisions with other particles; consequently, the particle's behavior is not significantly influenced by direct, physical impacts with other particles. As such, equations such as the VP equation are of great interest, as many applications in nuclear fusion rely on accurate quantification of plasma dynamics~\cite{ye2022quantum, baumann2021landau, gibbon2020introduction}.

In full, the collisionless Vlasov--Poisson equation under the effect of electric fields is given as,
\begin{equation}\label{eqn:vpe}
    \frac{\partial f_s}{\partial t} + \mathbf{v}_s \cdot \nabla_\mathbf{x}f_s + \frac{q_s}{m_s}\mathbf{E}\cdot \nabla_{\mathbf{v}}f_s=0.
\end{equation}
Here, $f_s$ represents the time-dependent distribution function of a particle species $s$, i.e., electrons or ions, with real-space coordinates $\mathbf{x}$ and velocity-space components~$\mathbf{v}$. The species mass is denoted by $m_s$ and the species charge is denoted by $q_s$. The prescribed electric field $\mathbf{E}$ acts in spatial coordinates $\mathbf{x}$ and is computed via Poisson's equation:
\begin{align}
  \nabla^2\phi  & =  -\frac{1}{\epsilon_0}\sum_s q_s 
     \int d\mathbf{v}\,f_s(\mathbf{x},\mathbf{v},t), \\  
  \mathbf{E}       & =  -\nabla\phi.
\end{align}
Here, $\phi$ is an electric field and $\epsilon_0$ is known as the \textit{permittivity of free space}; this is a fundamental constant that quantifies how much electric field can fill a vacuum, relating the electric charge density to the resulting electric field~\cite{gibbon2020introduction}. 

In its entirety, (\ref{eqn:vpe}) is a (6 + 1)-dimensional PDE, as the real and velocity space coordinates are typically 3-dimensional. The nonlinearity in (\ref{eqn:vpe}) arises from the definition of the electric field. In the following results, we assume that the ions are fixed in time (see Appendix~\ref{sec:ions}), and we focus on modeling the electron species. Furthermore, for this paper, we will be focused on a 1D1V (2 + 1)-dimensional version of the problem, aligning this setting with our previous results. Full implementation details of the solver used can be found in Appendix~\ref{sec:datagen-vp}.

\begin{figure*}
    \centering

    \begin{subfigure}[t]{0.92\textwidth}
        \centering
        \includegraphics[width=\linewidth]{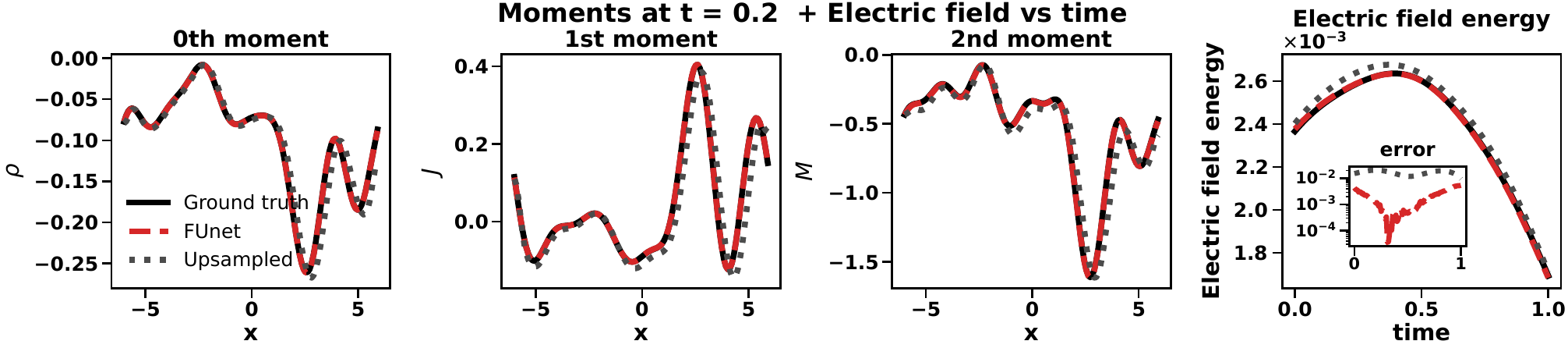}
        \caption{Randomized dataset observables and electric field energy density at $t=0.21$.}
        \label{fig:vp-moments-efielderr-t21}
    \end{subfigure}

    \vspace{0.6em}

    \begin{subfigure}[t]{0.92\textwidth}
        \centering
        \includegraphics[width=\linewidth]{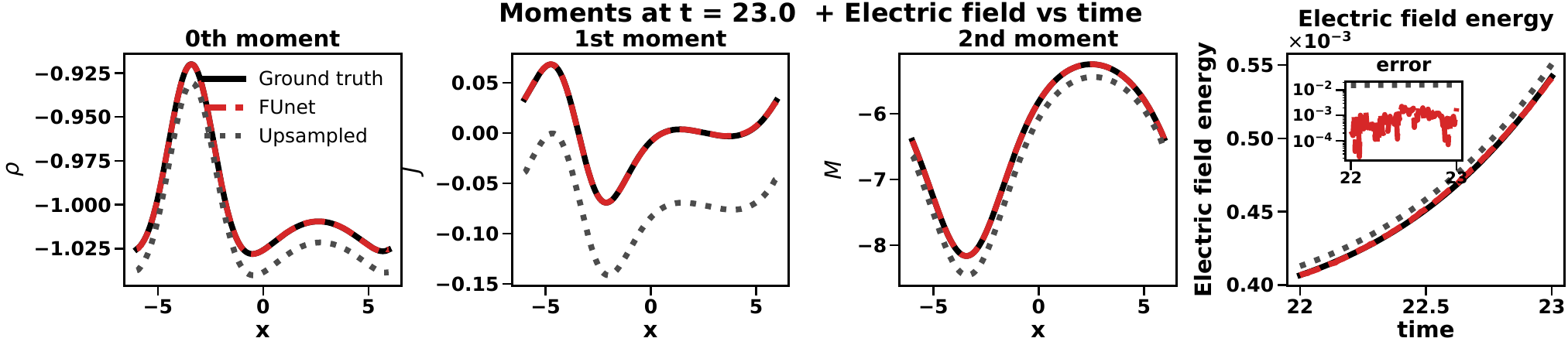}
        \caption{Two-stream instability observables and electric field energy density at $t=23.0$.}
        \label{fig:vp-moments-efielderr-t100}
    \end{subfigure}

    \caption{\justifying{\textbf{Prediction of moments at a fixed point in time for Phase~1 and error in electric field energy density for both the two-stream instability and randomized datasets}.}}
    \label{fig:vp-moments-efield-error}
\end{figure*}

\subsection{Physics-based metrics}
For the 1D1V Vlasov–Poisson test cases, we also track system‑specific, physics‑based diagnostics.
In addition to the metrics reported in Appendix \ref{sec:metrics}, we employ the following physics based metrics for the 1D1V VP system:  the negative-energy penalty, the entropy error, and the electric field error. First, we note that the VP electron distribution function is a probability density, and must remain non-negative. Any predicted negative values are nonphysical, and thus will degrade moment accuracy. The \textit{negative-mass penalty} is thus given by, 
\[
\mathcal{P}(t)=
\iint
\bigl|\,\min\bigl(\tilde f_e(x,v,t),\,0\bigr)\bigr|\,dv\,dx ,
\]
where
$\mathcal{P}(t)$ measures the total negative mass in phase space. Ideally, this quantity is as close to 0 as possible, as it is unphysical to have a negative mass.
Second, for the \textit{entropy error}, we note that for an exact solution, the entropy
\[
  S(t)=\iint f_e(x,v,t)\,\ln f_e(x,v,t)\,\,\mathrm dv\,\mathrm dx
\]
should be conserved in time.  Given a prediction \(\tilde f_e\) we define
\[
  \tilde S(t)=\iint \tilde f_e\,\ln\tilde f_e\,\mathrm dv\,\mathrm dx,
  \qquad
  \mathcal E_{S}(t)=
  \frac{\bigl|\tilde S(t)-S(t)\bigr|}{\bigl|S(t)\bigr|}.
\]
A small \(\mathcal E_{S}\) signals that the prediction has not introduced spurious velocity-space dissipation. Finally, we can measure the \textit{electric-field error}. 
From the predicted charge density, \(\tilde\rho(x,t)=1-\int\tilde f_e\,\mathrm dv\), we solve
\(\partial_{xx}\tilde\phi=\tilde\rho\) and set
\(\tilde E(x,t)=-\partial_x\tilde\phi\).
Comparing with the ground-truth field \(E=-\partial_x\phi\) gives
\[
  \mathcal E_{E}(t)=
  \frac{\bigl\lVert\tilde E-E\bigr\rVert_{2}}
       {\lVert E\rVert_{2}}.
\]
A low \(\mathcal E_{E}\) indicates that the predicted charge distribution produces an electrostatic force consistent with Gauss’ law, ensuring accurate particle acceleration and energy exchange.

 Furthermore, for Vlasov--Poisson systems, physical observables in the form of velocity moments are sometimes the quantities of interest for plasma physicists~\cite{gibbon2020introduction}. Many physically relevant quantities pertaining to plasmas are encoded within the velocity moments,  including the charge density, the current density, and the kinetic energy density. These moments and other higher-order moments provide a compact macroscopic fingerprint of the plasma's underlying dynamics in the form of quantities that physicists can measure. Concretely, given a species distribution \(f_s(\mathbf{x},\mathbf{v},t)\), we define its velocity moments by
\begin{equation}
  M^{(k)}_s(\mathbf{x},t)
  := \int_{\mathbb{R}^{d_v}} \mathbf{v}^{\otimes k}\, f_s(\mathbf{x},\mathbf{v},t)\, d\mathbf{v},
  \qquad k=0,1,2,\dots,
\end{equation}
where \(d_v\) is the velocity-space dimension and \(\mathbf{v}^{\otimes k}\) denotes the \(k\)-fold outer product. In particular, the zeroth moment is the number density,
\begin{equation}
  \rho_s(\mathbf{x},t) := \int_{\mathbb{R}^{d_v}} f_s(\mathbf{x},\mathbf{v},t)\, d\mathbf{v}
  = M^{(0)}_s(\mathbf{x},t),
\end{equation}
and the charge density appearing in Poisson's equation is the species-weighted sum
\begin{equation}
  \rho(\mathbf{x},t) := \sum_s q_s\,\rho_s(\mathbf{x},t)
  = \sum_s q_s \int_{\mathbb{R}^{d_v}} f_s(\mathbf{x},\mathbf{v},t)\, d\mathbf{v}.
\end{equation}
The first moment gives the current density,
\begin{equation}
  \mathbf{J}(\mathbf{x},t)
  := \sum_s q_s \int_{\mathbb{R}^{d_v}} \mathbf{v}\, f_s(\mathbf{x},\mathbf{v},t)\, d\mathbf{v}
  = \sum_s q_s\, M^{(1)}_s(\mathbf{x},t),
\end{equation}
and the second moment we report in our results is
\begin{equation}
  \mathbf{M}_s(\mathbf{x},t)
  := \int_{\mathbb{R}^{d_v}} \mathbf{v}\otimes\mathbf{v}\, f_s(\mathbf{x},\mathbf{v},t)\, d\mathbf{v}
  = M^{(2)}_s(\mathbf{x},t).
\end{equation}


\squeezetable
\begin{table*}[ht!]
\setlength{\tabcolsep}{3pt}        
\renewcommand{\arraystretch}{1.05} 
\centering
\footnotesize                 

\newcolumntype{C}{>{\centering\arraybackslash}X}
\newcolumntype{b}{>{\columncolor{gray!8}}c}

\begin{tabular*}{\linewidth}{@{\extracolsep{\fill}} l b c b c @{}}
\toprule
& \multicolumn{2}{c}{\cellcolor{gray!6}\textbf{Two-stream}}
& \multicolumn{2}{c}{\cellcolor{gray!6}\textbf{Randomized ($m_x=m_v=4$)}} \\
\cmidrule(lr){2-3}\cmidrule(lr){4-5}
& \cellcolor{gray!15}\textbf{Bicubic} & \textbf{FUnet}
& \cellcolor{gray!15}\textbf{Bicubic} & \textbf{FUnet} \\
\midrule

\multicolumn{5}{l}{\textbf{(a) Flow diagnostics}} \\ \midrule
$f_e$ &
$5.799\times10^{-2}$ &
$\mathbf{4.378\times10^{-4}}\pm 1.7\times10^{-5}$ &
$2.468\times10^{-1}$ &
$\mathbf{6.160\times10^{-3}}\pm 4.4\times10^{-4}$ \\
\midrule

\multicolumn{5}{l}{\textbf{(b) Moment errors}} \\ \midrule
$\rho$ &
$1.242\times10^{-2}$ &
$\mathbf{7.466\times10^{-5}}\pm 4.3\times10^{-6}$ &
$1.229\times10^{-1}$ &
$\mathbf{3.439\times10^{-3}}\pm 3.1\times10^{-4}$ \\[0.3em]

$J$ &
$1.746\times10^{0}$ &
$\mathbf{4.547\times10^{-3}}\pm 3.3\times10^{-4}$ &
$2.895\times10^{-1}$ &
$\mathbf{9.912\times10^{-3}}\pm 1.3\times10^{-3}$ \\[0.3em]

$M_2$ &
$3.758\times10^{-2}$ &
$\mathbf{1.458\times10^{-4}}\pm 8.8\times10^{-6}$ &
$1.547\times10^{-1}$ &
$\mathbf{7.575\times10^{-3}}\pm 9.1\times10^{-4}$ \\
\bottomrule
\end{tabular*}

\caption{\justifying{\textbf{Phase 1 test set results for the Vlasov--Poisson benchmarks.} Errors are reported in relative $L_2$. Columns group the two-stream instability and randomized distribution ($m_x=m_v=4$) test cases. Bold indicates the best \emph{mean} within each dataset for a given metric. Bicubic interpolation is reported as means only since it is deterministic. }}
\label{tab:phase1_vp_combined}
\end{table*}

\subsection{Phase 1 test results}
We begin by demonstrating the test set results of Phase 1 for the 1D1V VP system for both the two-stream instability and randomized electron datasets. For testing each different system, we draw new initial conditions from the same random distributions that were used for training, and we propagate both the coarse-grain and fine-grain solvers until the same fixed final physical time. For the randomized dataset, we use $T=1$ and for the two-stream instability dataset, we use $T=24$. For the two-stream instability data set, we learn a $64\times64 \to 256\times256$ mapping, and for the randomized dataset, we learn a $32\times32\to 128\times128$ mapping. Finally, since FUnet consistently performed best across the preceding case studies, we report the results below only for this choice of neural operator. 

Figures~\ref{fig:vp-phase1-fe-twostream-random} and \ref{fig:vp-moments-efield-error} and Table~\ref{tab:phase1_vp_combined} summarize the Phase~1 test set reconstruction performance for the 1D1V VP benchmarks. Figure~\ref{fig:vp-phase1-fe-twostream-random} shows representative fine-grid snapshots of the electron distribution \(f_e(x,v,t)\) for an unseen two-stream trajectory (at \(t=22.5\)) and an unseen randomized trajectory (at \(t=0.5\)), comparing bicubic upsampling against the FUnet prediction. Here, we observe error heatmaps because they localize where the super-resolution mapping succeeds or fails in phase space. In contrast to the previous cases, we found it more difficult to discern differences purely visually based on the fields themselves. In both test cases, the bicubic upsampling baseline produces an overly smooth reconstruction with systematic magnitude mismatch relative to the ground truth, while the FUnet more faithfully matches both the amplitude and the localized fine-scale variability. This is consistent with the previous results from Section~\ref{sec:results}.

Quantitatively, Table~\ref{tab:phase1_vp_combined} reports mean relative \(L_2\) errors over \(20\) new test trajectories for both the field-level reconstruction of \(f_e\) and the derived velocity-moment observables \(\rho\), \(J\), and \(M\). For the two-stream case, the relative $L_2$ error in \(f_e\) drops from \(5.79\times 10^{-2}\) (bicubic) to \(4.37\times 10^{-4}\) (FUnet), and the moment errors improve by multiple orders of magnitude (e.g., \(\rho\) from \(1.24\times 10^{-2}\) to \(7.46\times 10^{-5}\)). For the randomized dataset, FUnet similarly reduces the reconstruction error in \(f_e\) from \(2.46\times 10^{-1}\) to \(6.16\times 10^{-3}\), with consistent reductions in \(\rho\), \(J\), and \(M\). Figure~\ref{fig:vp-moments-efield-error} further illustrates that these gains persist at the observable level: the predicted moment profiles closely track the ground truth at fixed times, and the resulting electric-field energy error remains small, indicating that the learned lift improves not only pixel-wise agreement, but also physically relevant aggregates driven by phase-space structure. Finally, for the sake of clarity, we add Figures~\ref{fig:vlasov_phase1_randomized_full} and ~\ref{fig:vlasov_phase1_two_stream_full}, which shows the evolution in time of the ground truth dynamics, baseline, and our Phase~1 super-resolution predictions.

An interesting question one can investigate is the transferability of high frequency training to a lower frequency test set with no fine-tuning. This zero-shot generalization phenomenon has also been noticed in a variety of different scientific contexts, including weather and fluid turbulence modeling \cite{subel2023explaining, darman2025fourier, chattopadhyay2023long}, although strong negative results have recently been established \cite{FalsePromizeZeroShot_TR}. In Appendix ~\ref{sec:supresults}, we explore this phenomenon in our electron dynamics context. The main findings indicate the ability to transfer high frequency training to zero-shot (i.e., no fine-tuning) low frequency generalization, but not vice-versa (Figure~\ref{fig:vlasov-mx-sweeps} and Table~\ref{tab:phase1_moments}). Furthermore, mathematical and empirical analysis of the model weights demonstrate an interpretable analysis that shows that high frequency training yields a more expressive model, thus is capable of generating fine-scale features that low frequency training misses (Figure~\ref{fig:spec-norms} and Figure~\ref{fig:ntk-spec-eig}).

In summary, for the Vlasov--Poisson benchmarks, the Phase~1 learned lift \(\N_\theta\) reconstructs the fine-grid phase-space state from coarse inputs with high fidelity; in representative test cases, the resulting field-level agreement is roughly two orders of magnitude better than bicubic upsampling, and this improvement propagates to the derived observables \(\rho\), \(J\), and \(M\), where we observe up to three orders of magnitude in terms of error improvement.

\begin{figure*}[!t]
  \centering
  \includegraphics[width=\linewidth]{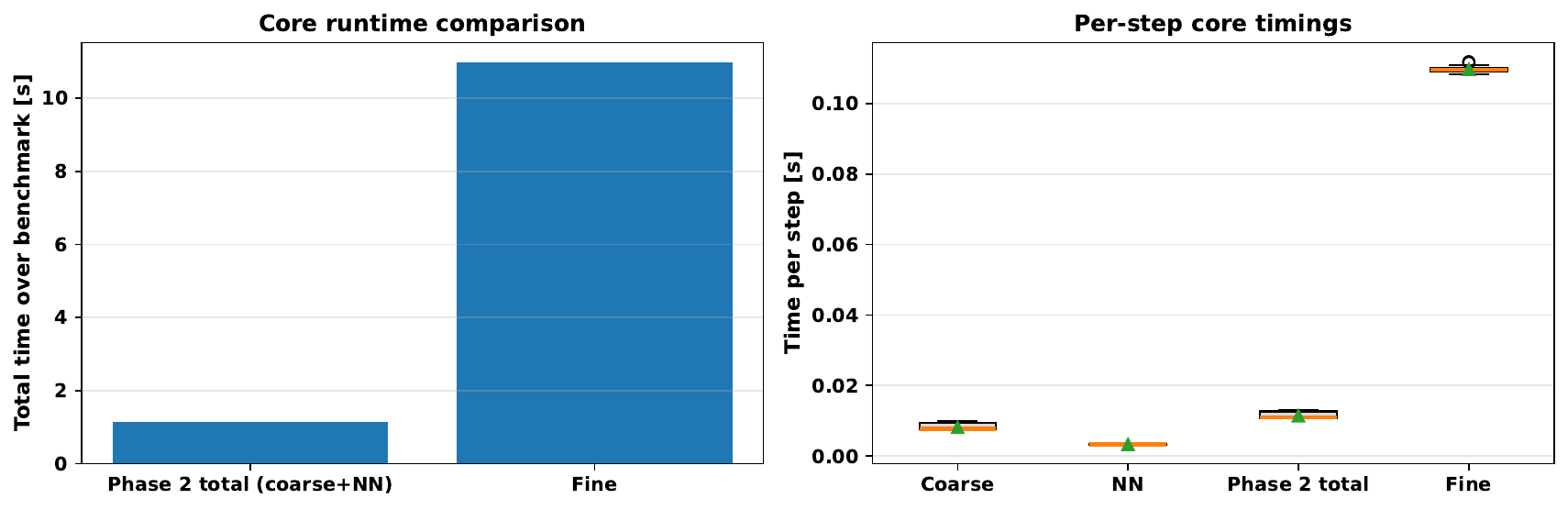}
  \caption{\justifying{\textbf{Core runtime comparison and per-step timings.} We compare one fine-grid Vlasov–Poisson step on a $256\times256$ phase-space grid against our Phase~2 extrapolation, which consists of the projection into the coarse grid, the time stepping on the coarse grid, and a single neural super-resolution query to the fine-resolution. \textbf{Left}: total time accumulated over a fixed number of benchmark steps. \textbf{Right}: decomposed per-step distributions. The Phase~2 pathway achieves a consistently lower runtime envelope on our node, with approximately $11\times$ reduction.}}
  \label{fig:timings}
\end{figure*}

\begin{table*}[ht!]
\setlength{\tabcolsep}{2pt}        
\renewcommand{\arraystretch}{0.95} 
\centering
\scriptsize

\newcolumntype{L}{l}
\newcolumntype{B}{>{\columncolor{gray!10}}c}
\newcolumntype{S}{c}
\begin{tabular*}{\textwidth}{@{\extracolsep{\fill}} L B S @{}}
\toprule
& \multicolumn{1}{c}{\cellcolor{gray!20}\textbf{Baseline}} & \multicolumn{1}{c}{\cellcolor{gray!8}\textbf{RELift model}} \\
\cmidrule(lr){2-2}\cmidrule(lr){3-3}
\textbf{Metric} & \cellcolor{gray!20}\textbf{Upsampled} & \textbf{FUnet} \\
\midrule

\multicolumn{3}{l}{\textbf{Two-stream}} \\ \midrule
\multicolumn{3}{l}{\textbf{(a) Flow / physics diagnostics}} \\ \midrule
rel.\,$L_2$ ($\downarrow$)
  & $5.18\!\times\!10^{-2}$
  & $\mathbf{7.86\!\times\!10^{-3}}\pm 3.70\!\times\!10^{-4}$ \\[0.15em]
SSIM ($\uparrow$)
  & $9.85\!\times\!10^{-1}$
  & $\mathbf{9.99\!\times\!10^{-1}}\pm 5.60\!\times\!10^{-5}$ \\[0.15em]
Spect.\ MSE ($\downarrow$)
  & $2.07$
  & $\mathbf{5.89\!\times\!10^{-2}}\pm 4.50\!\times\!10^{-3}$ \\[0.15em]
PSNR ($\uparrow$)
  & $31.00$
  & $\mathbf{49.10}\pm 6.60$ \\[0.15em]
Neg-mass pen.\ ($\downarrow$)
  & $\mathbf{0.00}$
  & $1.95\!\times\!10^{-5}\pm 3.00\!\times\!10^{-7}$ \\[0.15em]
Entropy err.\ ($\downarrow$)
  & $7.09\!\times\!10^{-4}$
  & $\mathbf{5.66\!\times\!10^{-4}}\pm 6.70\!\times\!10^{-5}$ \\[0.15em]
$E$-field err.\ ($\downarrow$)
  & $8.08\!\times\!10^{-2}$
  & $\mathbf{4.12\!\times\!10^{-2}}\pm 1.40\!\times\!10^{-3}$ \\
\midrule
\multicolumn{3}{l}{\textbf{(b) Velocity-moment errors (relative $L_2$)}} \\ \midrule
$\rho$
  & $3.38\!\times\!10^{-3}$
  & $\mathbf{1.08\!\times\!10^{-3}}\pm 5.60\!\times\!10^{-5}$ \\[0.15em]
$J$
  & $1.56$
  & $\mathbf{7.27\!\times\!10^{-2}}\pm 2.80\!\times\!10^{-3}$ \\[0.15em]
$M$
  & $8.75\!\times\!10^{-3}$
  & $\mathbf{5.31\!\times\!10^{-3}}\pm 2.50\!\times\!10^{-4}$ \\
\midrule

\multicolumn{3}{l}{\textbf{Randomized ($m_x=m_v=4$)}} \\ \midrule
\multicolumn{3}{l}{\textbf{(a) Flow / physics diagnostics}} \\ \midrule
rel.\,$L_2$ ($\downarrow$)
  & $1.66\!\times\!10^{-1}$
  & $\mathbf{9.32\!\times\!10^{-2}}\pm 6.30\!\times\!10^{-3}$ \\[0.15em]
SSIM ($\uparrow$)
  & $9.17\!\times\!10^{-1}$
  & $\mathbf{9.63\!\times\!10^{-1}}\pm 3.10\!\times\!10^{-3}$ \\[0.15em]
Spect.\ MSE ($\downarrow$)
  & $1.83\!\times\!10^{-1}$
  & $\mathbf{7.94\!\times\!10^{-2}}\pm 7.70\!\times\!10^{-3}$ \\[0.15em]
PSNR ($\uparrow$)
  & $3.28\!\times\!10^{1}$
  & $\mathbf{4.15\!\times\!10^{1}}\pm 1.00\!\times\!10^{0}$ \\[0.15em]
Neg-mass pen.\ ($\downarrow$)
  & $5.99\!\times\!10^{-1}$
  & $\mathbf{7.89\!\times\!10^{-2}}\pm 4.10\!\times\!10^{-3}$ \\[0.15em]
Entropy err.\ ($\downarrow$)
  & $7.68\!\times\!10^{-2}$
  & $\mathbf{6.60\!\times\!10^{-3}}\pm 6.40\!\times\!10^{-4}$ \\[0.15em]
$E$-field err.\ ($\downarrow$)
  & $\mathbf{8.84\!\times\!10^{-3}}$
  & $8.77\!\times\!10^{-2}\pm 7.10\!\times\!10^{-3}$ \\
\midrule
\multicolumn{3}{l}{\textbf{(b) Velocity-moment errors (relative $L_2$)}} \\ \midrule
$\rho$
  & $6.06\!\times\!10^{-2}$
  & $\mathbf{4.33\!\times\!10^{-2}}\pm 2.50\!\times\!10^{-3}$ \\[0.15em]
$J$
  & $1.34\!\times\!10^{-1}$
  & $\mathbf{8.52\!\times\!10^{-2}}\pm 5.20\!\times\!10^{-3}$ \\[0.15em]
$M$
  & $1.14\!\times\!10^{-1}$
  & $\mathbf{7.99\!\times\!10^{-2}}\pm 4.70\!\times\!10^{-3}$ \\
\bottomrule
\end{tabular*}

\caption{\justifying
\textbf{Vlasov--Poisson Phase 2 rollout metrics (top: two-stream; bottom: randomized distribution).}
Block (a) lists flow/physics diagnostics; block (b) reports relative $L_2$ errors for velocity moments. Within each case, the best mean per row is highlighted in bold.}
\label{tab:vlasov_phase2_combined_twostream_randomized}
\end{table*}

\subsection{Phase 2 test results}
\begin{figure*}[htbp]
    \centering

    \begin{minipage}[b]{\linewidth}
        \centering
        \includegraphics[width=\linewidth]{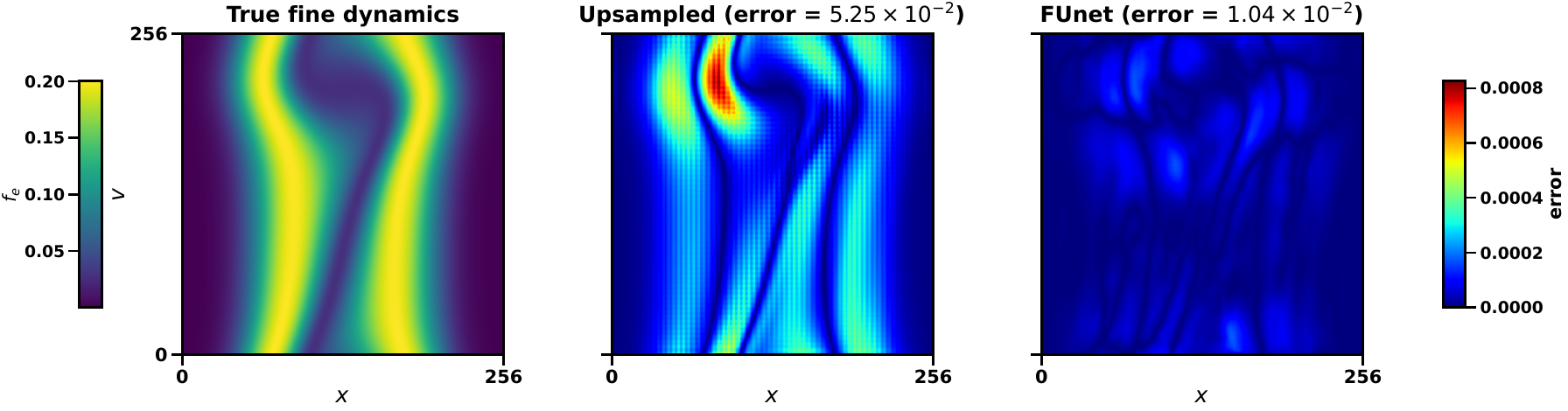}
    \end{minipage}

    \vspace{1em}

        \begin{minipage}[b]{\linewidth}
        \centering
        \includegraphics[width=\linewidth]{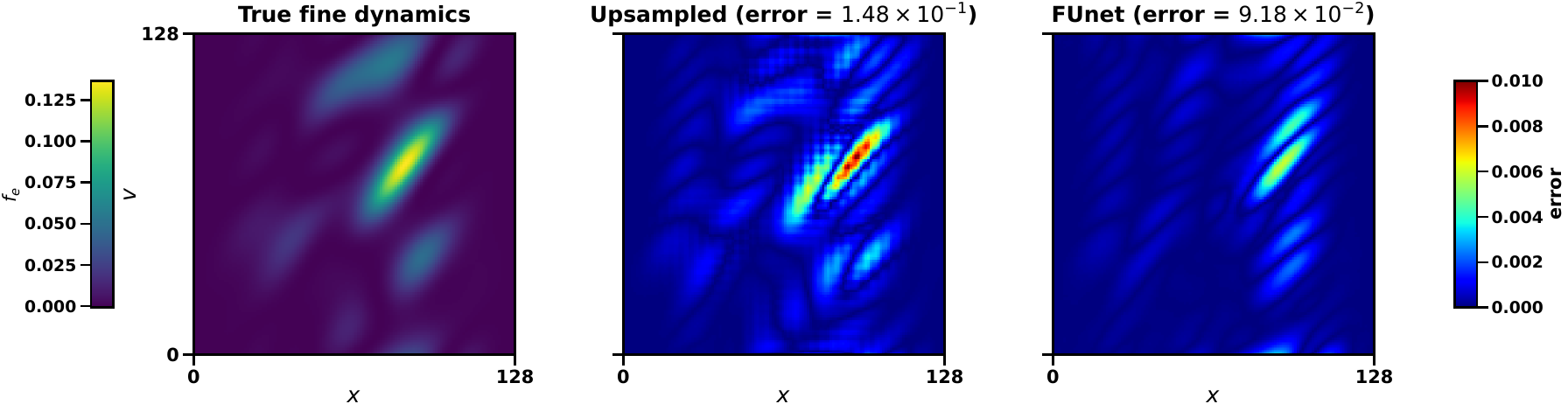}
    \end{minipage}

    \caption{\justifying \textbf{Phase~2 long-horizon Vlasov--Poisson rollouts for the two-stream and randomized datasets: ground truth fine-grid \(f_e\), bicubic upsampling baseline, and Phase~2 RELift (coarse step \(+\) learned lift), with relative \(L_2\) errors reported in each panel.}}
    \label{fig:vlasov_phase2_two_stream+random}
\end{figure*}

These Phase~1 results motivate the Phase~2 evaluation, where we compose the coarse Vlasov--Poisson update with the learned lift to obtain an effective fine-grid propagator and assess both efficiency and rollout accuracy. For the two-stream instability, we roll out to physical time \(T=28\), corresponding to \(400\) time steps beyond the training window (a \(4\times\) extrapolation relative to the training horizon); for the randomized dataset, we roll out to \(T=3\), corresponding to \(200\) steps beyond the training window (a \(2\times\) extrapolation).

We first quantify the core runtime advantage of the Phase~2 pipeline for the VP datasets in Figure~\ref{fig:timings} by comparing a single fine-grid step on a $256\times256$ phase-space mesh against one Phase~2 step consisting of a $64\times64$ coarse physics update followed by a single neural upscaling query on GPU. Both the accumulated time benchmark and the per-step distributions show that the Phase~2 pathway maintains a consistently lower runtime envelope than the fine solver, with the coarse update dominating the Phase~2 per-step cost and the neural lift contributing only a modest overhead. Note that, in contrast to the preceding PDE solvers, we do not additionally coarse-grain the time step here. Due to the greater complexity of the Vlasov--Poisson update in time, substantial wall-clock savings are already realized under the solver's natural time stepping. Finally, we note that the fine-scale dataset requires approximately 80~GB of storage, whereas the coarse-scale dataset together with the neural network weights requires only about 4~GB, yielding roughly a $20\times$ reduction in storage cost.

\begin{figure*}[htbp]
    \centering

    \begin{minipage}[b]{\linewidth}
        \centering
        \includegraphics[width=\linewidth]{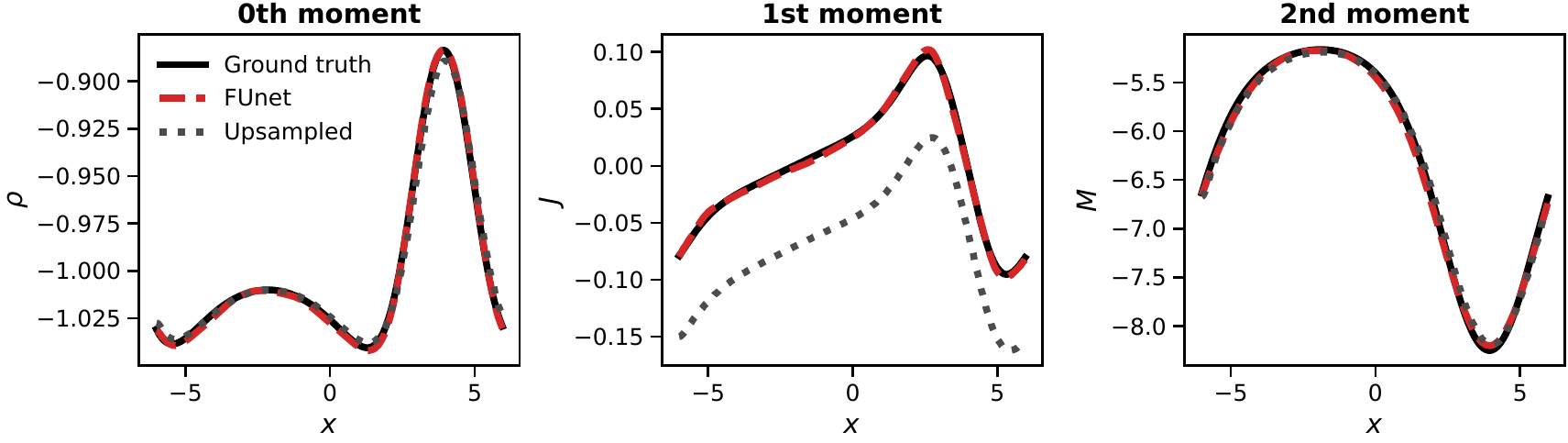}
    \end{minipage}

    \vspace{1em}

        \begin{minipage}[b]{\linewidth}
        \centering
        \includegraphics[width=\linewidth]{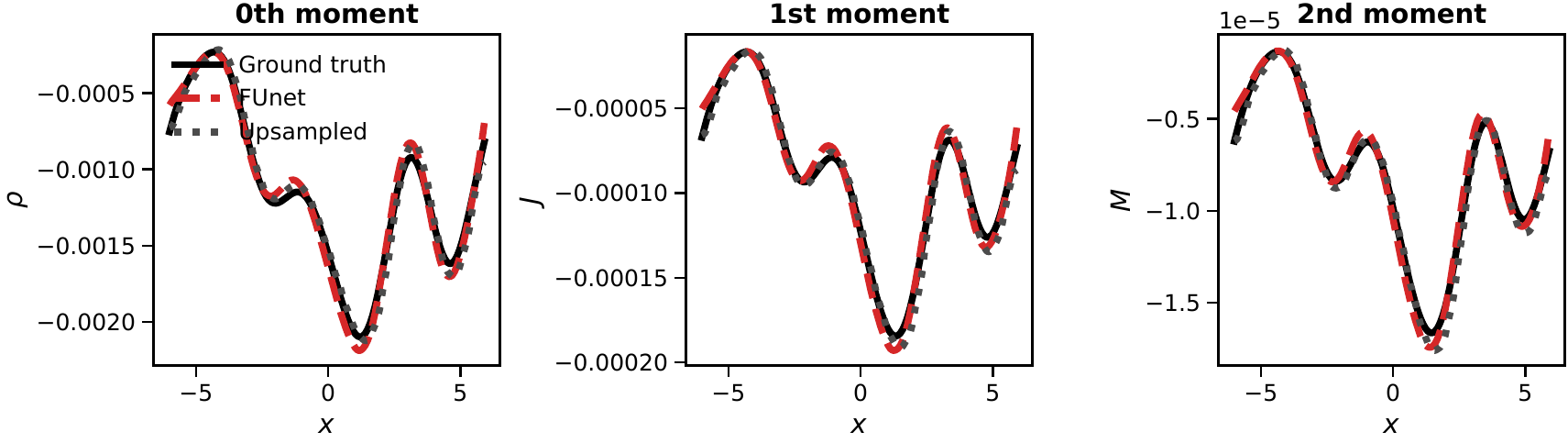}
    \end{minipage}

    \caption{\justifying \textbf{Phase~2 Vlasov--Poisson rollout diagnostics for the two-stream instability \textit{(top) }and randomized \textit{(bottom)} datasets}: velocity moments \(\rho\), \(J\), and \(M\) comparing ground truth, bicubic upsampling, and Phase~2 RELift predictions.}

    \label{fig:observables_vlasov_phase2_two_stream+random}
\end{figure*}

Turning to accuracy, Figure~\ref{fig:vlasov_phase2_two_stream+random} compares representative snapshots in the future time extrapolation for the two-stream and randomized datasets. In both cases, advancing on the coarse grid and lifting via bicubic interpolation yields noticeably larger phase-space mismatch compared to Phase~2 extrapolation with the FUnet. In the shown examples, the reported errors decrease from $5.25\times10^{-2}$ to $1.04\times10^{-2}$ (two-stream) and from $1.48\times10^{-1}$ to $9.18\times10^{-2}$ (randomized). To make the implications for physically meaningful aggregates explicit, Figure~\ref{fig:observables_vlasov_phase2_two_stream+random} plots the 0th, 1st, and 2nd moments $\rho$, $J$, and $M$. Here, we see that the learned Phase~2 rollouts track the ground-truth profiles more closely than bicubic lifting, most prominently for the current density $J$.

Finally, Table~\ref{tab:vlasov_phase2_combined_twostream_randomized} summarizes test-set rollout metrics. For the two-stream case, the learned Phase~2 propagator improves field-level errors (e.g., rel.\ $L_2$ from $5.18\times10^{-2}$ to $7.86\times10^{-3}$, spectral MSE from $2.07$ to $5.89\times10^{-2}$, PSNR from $31.0$ to $49.1$) and reduces physics/observable discrepancies (e.g., $J$ error from $1.56$ to $7.27\times10^{-2}$, and electric-field error from $8.08\times10^{-2}$ to $4.12\times10^{-2}$).

For the randomized case, FUnet again improves most metrics, including relative $L_2$ error ($1.66\times10^{-1} \text{ vs. }9.32\times10^{-2}$), spectral MSE ($1.83\times10^{-1} \text{ vs. }7.94\times10^{-2}$), and entropy error ($7.68\times10^{-2} \text{ vs. }6.60\times10^{-3}$), while the relative electric-field metric is the main exception, favoring bicubic lifting ($8.84\times10^{-3}$ vs.\ $8.77\times10^{-2}$). We suspect that the electric field diagnostic in this dataset can favor a smoother lift because \(E\) is obtained from \(\rho\) via Poisson’s equation, which emphasizes large-scale (low-\(k\)) charge modes and effectively filters small-scale phase-space errors. Consequently, bicubic interpolation may achieve a smaller \(E\)-error by damping high-frequency density fluctuations, even when its phase-space and moment-level matches are worse. For the sake of clarity, we add  Figures~\ref{fig:vlasov_phase2_randomized_full}--~\ref{fig:vlasov_phase2_randomized_full_mom} which show the evolution in time of the ground truth dynamics, baseline, and our Phase~2 extrapolation predictions for the electron species of both datasets and its relevant moments.

Finally, as an additional long-horizon stress test, Figure~\ref{fig:vp-long-time} shows that the same Phase~2 construction remains qualitatively stable and physically coherent far beyond the rollout horizons used above, including extrapolations to $T=64$ for the two-stream case and 
$T=20$ for the randomized case, while also remaining more accurate than the baseline. These results further support that the composition-based effective propagator can sustain accurate future-time extrapolation well outside the temporal window seen in Phase~1 training; full details are given in Section~\ref{sec:ablation}.

\section{Conclusion}\label{sec:conclusion}
In this paper, we introduced RELift, a two-phase learning framework that couples a super-resolution neural operator with a coarse-grid numerical integrator to approximate the fine-scale dynamics of PDEs.  
Across different PDEs studied in the work, RELift successfully reconstructs high-resolution fine-grid solutions from coarse inputs (Phase~1) and delivers stable, accurate long-time forecasts (Phase~2) of the underlying dynamics.  
We also provided a systematic error analysis of the RELift-learned effective fine-grid propagator, demonstrating that the derived local and global approximation error bounds closely follow the observed numerical errors in simulations of a representative nonlinear PDE.
While this work focuses on PDE dynamics on regular geometries with uniform discretizations, future extensions may incorporate adaptive meshes, irregular spatial domains, or more complex boundary conditions.  
Overall, our results highlight RELift as a practical and flexible framework for pairing learned neural operators with existing numerical solvers, offering a scalable path toward high-resolution forecasting of spatiotemporal systems.

\section*{Acknowledgments}
This research was sponsored by the Office of Naval Research
and was accomplished under Grant Number W911NF-23-1-0153 (H. Bassi and H.S. Bhat) and the U.S. Department of Energy, Office of Science, Office
of Advanced Scientific Computing Research and Office
of Basic Energy Sciences, Scientific Discovery through
Advanced Computing (SciDAC) program under Contract No. DE-AC02-05CH11231
 (H. Bassi, Y. Zhu, E. Ye, A. Dektor and C. Yang). 
This work is also supported by the Center for Computational Study of Excited-State
Phenomena in Energy Materials (C2SEPEM) at the Lawrence Berkeley National Laboratory, which is funded by the U.S. Department of Energy, Office of Science, Basic Energy Sciences, Materials Sciences and Engineering Division, under Contract No. DE-AC02-05CH11231, as part of the Computational Materials Sciences Program. (Y. Zhu and C. Yang)
This research used
resources of the National Energy Research Scientific Computing Center (NERSC), a U.S. Department of Energy Office of Science User Facility located at Lawrence Berkeley
National Laboratory, operated under Contract No. DE-AC02-
05CH11231, using NERSC awards ASCR-m4577 and ASCR-m1027.
The views and conclusions contained in this document
are those of the authors and should not be interpreted as representing the official policies, either expressed or implied, of
the Army Research Office or the U.S. Government. The U.S.
Government is authorized to reproduce and distribute reprints
for Government purposes notwithstanding any copyright notation herein.
M.W. Mahoney acknowledges DARPA, NSF, the DOE Competitive Portfolios grant, and the DOE SciGPT grant.
H. Bassi also gratefully acknowledges the additional fruitful discussions and invaluable insights provided to him by mentors and friends, including Changho Kim, Matthew Blomquist, and Pinchen Xie.


\section*{DATA AVAILABILITY}
Full details of each dataset are provided in Appendix \ref{app:data_generation}. The neural network models and weights, training, evaluation, and data generation scripts are provided at  \hyperlink{https://github.com/hbassi/RELift}{https://github.com/hbassi/RELift}.

\section*{AUTHOR DECLARATIONS}

\subsection*{Conflict of Interest}
The authors have no conflicts to disclose. 

\bibliographystyle{aipnum4-2}
\bibliography{references} 

\appendix

\setlength{\tabcolsep}{5pt}
\renewcommand{\arraystretch}{1.1}
\newcommand{\NRow}[4]{#1 & #3\if\relax\detokenize{#4}\relax\else\ \textit{(#4)}\fi\\}

\begin{table*}
\small
\caption{Core notation used in the main text.}
\label{tab:notation-core}
\begin{tabular*}{\textwidth}{@{\extracolsep{\fill}} l l @{}}
\hline
\textbf{Symbol} & \textbf{Meaning} \\
\hline
\NRow{$\Omega$}{}{Physical spatial domain (continuous)}{}
\NRow{$\Omega_c,\ \Omega_f$}{}{Coarse and fine computational grids/index sets discretizing $\Omega$}{}
\NRow{$X$}{}{State/function space for fields on $\Omega$ (e.g., $X=L^2(\Omega)$ or $L^1(\Omega)$)}{}

\NRow{$u_c(t),\ u_f(t)$}{}{Coarse and fine states at (continuous) time $t$}{}
\NRow{$u_c^n,\ u_f^n$}{}{Coarse and fine states at discrete step $n$ with step size $\Delta t$}{}
\NRow{$P$}{}{Restriction/projection from fine to coarse (spatial decimation via point sampling)}{}
\NRow{$\mathcal N_\theta$}{}{Learned lift / super-resolution operator (Phase~1 network)}{}

\NRow{$\mathcal P_f(\Delta t)$}{}{Reference fine one-step propagator (ground truth solver)}{}
\NRow{$\mathcal P_c(\Delta t)$}{}{Coarse numerical one-step propagator}{}
\NRow{$\widehat{\mathcal P}_f(\Delta t)$}{}{Phase~2 surrogate fine one-step propagator}{$\widehat{\mathcal P}_f(\Delta t):=\mathcal N_\theta\circ \mathcal P_c(\Delta t)\circ P$}

\NRow{$\Phi_\theta$}{}{Autoregressive (AR) baseline one-step map on the fine grid}{}

\NRow{$\Delta t$}{}{Timestep}{}
\NRow{$T$}{}{Number of time steps used for evaluation}{}
\NRow{$N_c,\ N_f$}{}{Number of grid points on coarse and fine discretizations}{}
\NRow{$K$}{}{Number of supervised training pairs}{}


\NRow{$N_{\text{traj}}$}{}{Number of training trajectories}{}
\NRow{$\mathrm{AR}$, $\mathrm{SR}$}{}{Shorthand for autoregressive and super-resolution settings}{}
\hline
\end{tabular*}
\end{table*}

\section{Intuition on approximate space-time separability and memory}\label{sec:stmem-details}

RELift couples a temporal update on the coarse grid with an instantaneous learned spatial lift to the fine grid. This architecture is most effective when a short-memory condition holds: conditioned on the present coarse field $u_c(t)$, the subgrid content is well approximated by a function of $u_c(t)$ alone,
\[
u_f(t)\approx\Phi\big(u_c(t)\big),
\]
so the one-step surrogate propagator
$\hat{\mathcal P}_f(\Delta t;\theta)=\mathcal N_\theta\circ \mathcal P_c(\Delta t)\circ P$
is a reasonable ansatz for the fine solver. In that sense, RELift exploits an approximate separation of temporal and spatial correlation: the coarse solver $\mathcal P_c$ carries the time advance while $\mathcal N_\theta$ provides an instantaneous spatial correction.

Related assumptions are also made intuitive in the SHRED~\cite{kutz2024shallow} method, which reconstructs the full field from a short temporal history via a shallow recurrent encoder and a nonlinear decoder. In our notation, its closure has the form \(u_f(t)\approx \Phi_h\!\big(u_c(t-\tau{:}t)\big)\), i.e., a history–dependent (non-Markovian) map. By contrast, RELift delegates temporal evolution to the coarse integrator \(\mathcal P_c\) and applies a time-local lift \(\mathcal N_\theta\), corresponding to a Markovian closure \(u_f(t)\approx \mathcal N_\theta\!\big(u_c(t)\big)\). 
From a separation of variables perspective, SHRED learns a nonlinear space–time factorization with an explicit recurrent temporal component, whereas RELift composes the physics-based temporal factor \(\mathcal P_c\) with an instantaneous spatial factor \(\mathcal N_\theta\).

\paragraph{\textbf{Linear intuition (heat / wave).}}
For linear, constant-coefficient PDEs with periodic boundary conditions (e.g., $\partial_t u=\nu\Delta u+\alpha f$ or $\partial_{tt}u=c^2\Delta u$), a Fourier projection onto low ($|k|\le k_c$) and high ($|k|>k_c$) modes commutes with the generator. Each wavenumber evolves independently in time:
\begin{equation}
\begin{aligned}
\widehat{u}_k(t) &= e^{\lambda_k t}\,\widehat{u}_k(0) + \text{(forced response)},\\
\lambda_k &=
\begin{cases}
-\nu|k|^2 & \text{(heat)},\\
\pm i\,c|k| & \text{(wave)}.
\end{cases}
\end{aligned}
\end{equation}

Thus, the \textit{temporal} evolution is multiplicative and modewise, while the \textit{spatial} content is set by which modes are present. For the heat equation, high modes decay on a time scale $\tau_k\sim(\nu|k|^2)^{-1}$; when $k_c$ is chosen so that $\tau_{k_c}$ is significantly smaller than the resolved time scale, the unresolved content is quickly chained to the resolved part, justifying an instantaneous lift. For the wave equation, although the dynamics are nondissipative, the $(k,\omega)$ spectrum concentrates along the dispersion ridge $\omega=c|k|$, so phase is primarily governed by $\mathcal P_c$ and the lift mainly corrects spatial aliasing.

\paragraph{\textbf{Nonlinear intuition (2D NSE)}}
Writing $\omega=\omega_<+\omega_>$ for low/high modes, we take $\Pi_<$ and $\Pi_>$ to be orthogonal Fourier projectors at a cutoff $k_c$:
for any solution on a grid $f(x)=\sum_{k\in\mathbb{Z}^2}\hat f_ke^{i kx}$,
\[
\Pi_< f := \sum_{|k|\le k_c}\hat f_ke^{i k\cdot x},\qquad
\Pi_> f := f-\Pi_< f,
\]
so that $\Pi_<^2=\Pi_<$, $\Pi_>^2=\Pi_>$, and $\Pi_<\Pi_>=0$. Both of these operators reduce a fine grid solution to its coarse counterpart, but the earlier spatial decimation projector $P$ is a practical grid-level downsampler (not necessarily orthogonal or spectral) used to form $u_c = Pu_f$, whereas the present projectors $\Pi_<$ and $\Pi_>$ are the ideal orthogonal Fourier projector at a cutoff $k_c$ (these are used only in this subsection for the sake of analysis). This can coincide with $P$ only when $P$ implements spectral truncation in conjunction with spatial decimation. Using these operators yields
\begin{equation*}
\begin{aligned}
\partial_t \omega_> &= \nu\,\Delta\omega_> + \Pi_>\,B(\omega_<,\omega_<) + \Pi_>\,B(\omega_<,\omega_>) \\
&\quad + \Pi_>\,B(\omega_>,\omega_>) ,
\end{aligned}
\end{equation*}
where $B$ is the advection operator in 2D vorticity form, and the corresponding Mori--Zwanzig identity \cite{zhu2018estimation,zhu2019mori,zhu2021effective} expresses the resolved dynamics as
\[
\partial_t \omega_< = \Pi_<\mathcal L\,\omega_< + \int_0^t K(t,s;\omega_<(s))\,ds + \eta(t),
\]
with a memory kernel $K$ and noise $\eta$. When viscosity and the cutoff $k_c$ yield a spectral gap, high modes decay on a fast time $\tau_{\mathrm{fast}}\sim(\nu k_c^2)^{-1}$ and the memory integral is well-approximated by a Markovian functional,
\(
\int_0^t K(\cdot)\,ds \approx \Phi(\omega_<(t)).
\)
RELift's Phase~1 learns an analogous operator to $\Phi$ from the dataset $(u_c,u_f)$; and Phase~2 wraps it around the coarse integrator so that the \textit{time} advance is done by $\mathcal P_c$ while the spatial  correction is refreshed each step by $\mathcal N_\theta$.

Thus, the short-memory condition is more likely to be satisfied for (i) diffusive systems and (ii) linear dispersive systems with constant coefficients, as well as (iii) moderately nonlinear, viscous flows where high modes relax quickly relative to resolved dynamics. The short-memory condition can degrade (and Phase~2 may drift) when (a) the Reynolds number is large enough such that subgrid time scales are not fast compared to resolved scales, (b) strong nonlocal interactions create long memory, or (c) the downsampling ratio is so aggressive that crucial interaction bands move entirely into the unresolved set; and, in this case, the NeurDE approach \cite{neurde_TR} may be appropriate.

\section{Model summaries}\label{sec:model-summaries}
Unless otherwise stated, models operate on stacks of coarse/fine grid solutions in time treated as channels  and use \texttt{GELU} activations. Specifically, all models except FNO-AR take in coarse grid inputs and and all models output fine grid outputs.

\smallskip
\subsection*{Common building blocks}
\paragraph{Spectral convolution (2D).}
All Fourier layers compute a real-space FFT over spatial dimensions $(H,W)$, retain only a low-frequency rectangle of size $(m_1{\times}m_2)$, apply a learned complex linear mixing on the retained modes, zero out all others, and return to the spatial domain by inverse FFT. Concretely, if $x\in\mathbb{C}^{B\times C_{\text{in}}\times H\times (W/2+1)}$ denotes the FFT, the retained block and the learned mixing is given by weights
$W\in\mathbb{C}^{C_{\text{in}}\times C_{\text{out}}\times m_1\times m_2}$, with
\[
\widehat{y}_{b,o,:,:}=\sum_{i=1}^{C_{\text{in}}} \widehat{x}_{b,i,:,:}\odot W_{i,o,:,:}, \qquad
y=\mathrm{iFFT}(\widehat{y}).
\]
A parallel $1{\times}1$ convolutional path is often summed with the spectral branch before a nonlinearity.

\paragraph{PixelShuffle upsampling.}
For scale factor $s\in\{2,4,8\}$ we map $C\mapsto s^2C$ channels by a $3{\times}3$ convolution followed by \texttt{PixelShuffle}$(s)$ and a GELU.

\subsection*{FUnet specification}
\begin{itemize}[leftmargin=1.2em]
\item \textbf{Input and output.}
The network maps a coarse input
\(x \in \mathbb{R}^{B \times C_{\mathrm{in}} \times H_c \times W_c}\)
to a fine output
\(y \in \mathbb{R}^{B \times C_{\mathrm{out}} \times H_f \times W_f}\),
with \(C_{\mathrm{in}} = C_{\mathrm{out}} = 100\).
The spatial resolutions satisfy \((H_f, W_f) = s (H_c, W_c)\) for an upsampling factor \(s\) (e.g., \(s=4\) in the \(32{\times}32 \to 128{\times}128\) experiments).

\item \textbf{Lifting layer.}
A \(1{\times}1\) convolution maps the 100 input channels to \(C_{\mathrm{lift}} = 128\) channels.

\item \textbf{Encoder.}
Two encoder stages operate at 128 channels.
Each stage is a \texttt{ConvBlock} consisting of two consecutive \(3{\times}3\) convolutions with GELU activations.
After the second \texttt{ConvBlock}, a \(2{\times}2\) max-pooling layer with stride 2 halves the spatial resolution, and we store skip connections from both encoder stages.

\item \textbf{Bottleneck.}
The bottleneck first applies a \texttt{ConvBlock} at 256 channels, followed by a Fourier layer that performs a spectral convolution with retained modes \((m_1,m_2) = (32,32)\) at 256 channels, and a final GELU activation.

\item \textbf{Decoder.}
The decoder consists of three upsampling stages based on pixel-shuffle operators with skip connections:
starting from the bottleneck representation (256 channels), successive stages
(1) reduce channels from 256 to 128 after combining with the deeper encoder skip,
(2) reduce channels from 128 to 64 after combining with the shallower encoder skip, and
(3) reduce channels from 64 to 32 on the finest grid.
The product of the pixel-shuffle factors over the three stages equals the overall upsampling factor \(s\).

\item \textbf{Projection (head).}
A \(3{\times}3\) convolution maps 32 channels to 32 channels, followed by GELU and a final \(3{\times}3\) convolution mapping 32 channels to \(C_{\mathrm{out}} = 100\).

\item \textbf{Nonlinearity and normalization.}
GELU activations are used throughout; no normalization layers are applied.
\end{itemize}

\subsection*{U-Net~\cite{ronneberger2015u} specification}
\begin{itemize}[leftmargin=1.2em]
\item \textbf{Input and output.}
The U-Net maps
\(x \in \mathbb{R}^{B \times C_{\mathrm{in}} \times H_c \times W_c}\)
to
\(y \in \mathbb{R}^{B \times C_{\mathrm{out}} \times H_f \times W_f}\),
with \(C_{\mathrm{in}} = C_{\mathrm{out}} = 100\) and \((H_f, W_f) = s (H_c, W_c)\).

\item \textbf{Base width.}
The base number of channels is \(C_{\mathrm{base}} = 64\).

\item \textbf{Encoder.}
There are two encoder levels.
At level~1, two \texttt{conv\_block}s operate at 64 channels, followed by a \(2{\times}2\) max-pooling layer (stride 2).
At level~2, two \texttt{conv\_block}s operate at 128 channels, again followed by max-pooling.
Each \texttt{conv\_block} consists of a \(3{\times}3\) convolution, batch normalization, and GELU activation.
Skip connections are taken from the outputs of the encoder levels before pooling.

\item \textbf{Bottleneck.}
Two \texttt{conv\_block}s operate at 256 channels on the most downsampled resolution.

\item \textbf{Decoder.}
Each decoder level upsamples via a \texttt{ConvTranspose2d} layer (stride 2), concatenates the corresponding encoder skip feature map, and applies two \texttt{conv\_block}s.
Channel counts mirror the encoder in reverse (256\(\to\)128, then 128\(\to\)64).

\item \textbf{Projection (head).}
A final \(1{\times}1\) convolution maps 64 channels to \(C_{\mathrm{out}} = 100\).
The output is then upsampled to the fine resolution \((H_f, W_f)\) by bilinear interpolation when needed to match the desired scale factor \(s\).

\item \textbf{Nonlinearity and normalization.}
All convolutions inside \texttt{conv\_block}s use GELU activations and batch normalization.
\end{itemize}

\subsection*{EDSR~\cite{lim2017enhanced} specification}
\begin{itemize}[leftmargin=1.2em]
\item \textbf{Input and output.}
The EDSR network maps
\(x \in \mathbb{R}^{B \times C_{\mathrm{in}} \times H_c \times W_c}\)
to
\(y \in \mathbb{R}^{B \times C_{\mathrm{out}} \times H_f \times W_f}\),
with \(C_{\mathrm{in}} = C_{\mathrm{out}} = 100\).

\item \textbf{Width and depth.}
We use \texttt{n\_feats} \(= 128\) feature channels and \texttt{n\_res\_blocks} \(= 16\) residual blocks.

\item \textbf{Head.}
A \(3{\times}3\) convolution maps \(C_{\mathrm{in}} = 100\) channels to 128 channels.

\item \textbf{Body.}
The body consists of 16 \texttt{ResBlock}s in series.
Each \texttt{ResBlock} applies a \(3{\times}3\) convolution, ReLU activation, and a second \(3{\times}3\) convolution, with a residual connection scaled by \texttt{res\_scale} \(= 0.1\).
A final \(3{\times}3\) convolution at 128 channels is applied at the end of the body, and the head output is added as a skip connection (standard EDSR formulation).

\item \textbf{Tail (upsampling and projection).}
An \texttt{Upsampler} implements the overall spatial upsampling factor \(s\) (for example, two pixel-shuffle layers with factor 2 each when \(s = 4\)).
A concluding \(3{\times}3\) convolution maps 128 channels back to \(C_{\mathrm{out}} = 100\).

\item \textbf{Input normalization.}
Before the head, we apply channelwise affine normalization using stored mean and standard deviation; the inverse transform is applied after the tail.

\item \textbf{Nonlinearity and normalization.}
ReLU activations are used in the residual blocks; no batch normalization is used inside residual blocks (as in the original EDSR).
\end{itemize}

\subsection*{FNO-SR~\cite{li2021fourier} specification}
\begin{itemize}[leftmargin=1.2em]
\item \textbf{Input and output.}
FNO-SR maps
\(x \in \mathbb{R}^{B \times C_{\mathrm{in}} \times H_c \times W_c}\)
to
\(y \in \mathbb{R}^{B \times C_{\mathrm{out}} \times H_f \times W_f}\),
with \(C_{\mathrm{in}} = C_{\mathrm{out}} = 100\).

\item \textbf{Width and depth.}
We use \texttt{width} \(= 64\) channels and stack 3 spectral blocks.

\item \textbf{Lifting layer.}
A \(1{\times}1\) convolution maps 100 input channels to 64 channels before the spectral blocks.

\item \textbf{Spectral blocks.}
Each block applies a \texttt{SpectralConv2d} operator with input and output width 64 and Fourier modes \((m_1,m_2) = (16,16)\), in parallel with a \(1{\times}1\) convolution in physical space.
The two outputs are summed and passed through a GELU activation.
The channel dimension remains 64 across all blocks.

\item \textbf{Upsampling and projection.}
After the spectral blocks, we apply bilinear interpolation from \((H_c, W_c)\) to \((H_f, W_f)\) with scale factor \(s\), followed by a \(1{\times}1\) convolution mapping 64 channels to \(C_{\mathrm{out}} = 100\).

\item \textbf{Nonlinearity and normalization.}
GELU activations are used after each spectral block; no normalization layers are applied.
\end{itemize}

\subsection*{FNO-AR specification}
\begin{itemize}[leftmargin=1.2em]
\item \textbf{Input and output.}
The autoregressive FNO operates on the fine grid.
Given a temporal context of \(K\) past fine snapshots stacked as channels,
\(x \in \mathbb{R}^{B \times K \times H_f \times W_f}\),
it outputs \(K\) channels of the same spatial shape.
In our experiments we take \(K = 100\); the predicted next snapshot is read from the output channels in the usual autoregressive manner.

\item \textbf{Width and depth.}
We use \texttt{width} \(= 64\) channels and 4 spectral blocks.

\item \textbf{Lifting layer.}
At each spatial location, a per-pixel linear map \(\mathbb{R}^K \to \mathbb{R}^{64}\) (implemented as a \(1{\times}1\) convolution) lifts the \(K\) input channels to 64 channels.

\item \textbf{Spectral blocks.}
Each block applies a \texttt{SpectralConv2dAR} operator with 64 input and output channels and retained modes \((m_1,m_2) = (16,16)\), in parallel with a \(1{\times}1\) convolution.
The outputs are summed and, in the first three blocks, passed through a GELU activation.
The channel dimension remains 64.

\item \textbf{Head.}
A per-pixel multi-layer perceptron (implemented as two linear layers with a GELU activation) maps \(\mathbb{R}^{64} \to \mathbb{R}^{K}\), producing the \(K\) output channels at each spatial location.

\item \textbf{Nonlinearity and normalization.}
GELU activations are used in the spectral blocks and in the head MLP; no normalization layers are applied.
\end{itemize}

\smallskip

\textbf{\text{Model training.} } We selected all configurations via standard hyperparameter sweeps and report, for each method, the best validation configuration. All models used Adam~\cite{kingma2014adam} with cosine–annealed learning rates; the default peak rate was $\eta\in\{1{\times}10^{-3},\,5{\times}10^{-4},\,3{\times}10^{-4}\}$ and weight decay $\in\{0,\,10^{-4}\}$. Phase~1 was trained for 5000 epochs and Phase~2 for 2500 epochs and the checkpoint with the lowest validation error was kept for Phase~1 and fine-tuned for Phase~2. For convolutional SR baselines (U-Net, EDSR) we swept depth (encoder/decoder blocks) $\in\{2,4,6\}$, base channels $\in\{32,64,128\}$, residual blocks per stage $\in\{4,8,16\}$, kernel size $\in\{3,5\}$, upsampler (pixel-shuffle vs nearest neighbor), normalization (none vs. batch), activation (ReLU vs. GELU), and dropout $\in\{0,0.05,0.1\}$. For Fourier models (FUnet, FNO-SR, FNO-AR) we swept layer count $\in\{4,6,8\}$, hidden width $\in\{32,64,128\}$, Fourier truncation per axis $k_x,k_y\in\{8,16,32\}$, and activation (ReLU vs. GELU). Once a set of parameters is chosen, training RELift models in Phase 1 and Phase 2 with all models tested took approximately 3 GPU hours (5000 + 2500 epochs on 1000 samples) on a NERSC Perlmutter GPU node with 1 80GB Nvidia A100.

\section{Dataset generation}\label{app:data_generation}
We generate paired fine–coarse datasets for three canonical PDEs on periodic domains. Both resolutions are advanced with the same time step $\Delta t$.

\paragraph{\textbf{Navier–Stokes (2D vorticity).}}
For each trajectory, the fine–grid initial vorticity $\omega_0$ is drawn i.i.d.\ from $\mathcal N(0,1)$, transformed to Fourier space, sharply low–pass filtered to modes $\{|k_x|,|k_y|\le 16\}$, real-part taken, and inverse–transformed to yield a band–limited field. The dynamics
\[
\partial_t \omega = -\,u\cdot\nabla\omega + \nu\,\Delta\omega + f,\qquad
u=\nabla^\perp\psi, -\Delta\psi=\omega,
\]
are advanced with a pseudo–spectral forward–Euler scheme under periodic boundary conditions. Parameters: viscosity $\nu=10^{-4}$, step size $\Delta t=0.01$, and forcing
$f(x,y)=0.025\bigl[\sin\bigl(2\pi(x+y)\bigr)+\cos\bigl(2\pi(x+y)\bigr)\bigr]$.
The identical integrator is applied on both spatial resolutions.

\paragraph{\textbf{Heat equation.}}
The fine-grid initial temperature $u_0(x,y)\sim \mathrm U(0,1)$; the coarse initial field is obtained by point decimation of $u_0$. Each Fourier mode $(k,\ell)$ is updated with the exact integrating–factor step using diffusivity $\nu=10^{-3}$, forcing amplitude $\alpha=10^{-2}$, and $\Delta t=0.01$. The update is diagonal in Fourier space, so the same routine is used at both resolutions without additional dealiasing.

\paragraph{\textbf{Wave equation.}}
The initial velocity is $g_0\equiv 0$ and the initial displacement is
\[
f_0(x,y)=\sum_{m=1}^{8} A_m \sin\!\bigl(2\pi(k_{x,m}x+k_{y,m}y)+\varphi_m\bigr),
\]
with amplitudes $A_m\sim\mathrm U[0.5,1.0]$, wavenumbers $(k_{x,m},k_{y,m})\in\{1,2,3\}^2$, and phases $\varphi_m\sim\mathrm U[0,2\pi)$. Each trajectory is linearly rescaled to $[0,1]$. Time integration uses an explicit leap–frog scheme with the finite–difference Laplacian $\Delta_h$ and periodic boundaries; parameters are wave speed $c=0.5$ and $\Delta t=0.01$, yielding $\mathrm{CFL}=c\,\Delta t\sqrt{2}<1$.

\section{\textbf{Metrics}}\label{sec:metrics}

To assess the performance of all models and baselines, we adopt the following metrics: 

\textbf{Relative $L_2$ error.} For a predicted \(\tilde u\) and the corresponding true \(u\) on the same \(N_f\times N_f\) grid,  the \textit{relative $L_2$ error} is
\begin{equation}\label{eqn:rel-l2}
 \mathcal{E}(t)
  =
  \frac{\|\,\tilde {u}(\cdot,\cdot,t)-u(\cdot,\cdot,t)\|_{2}}
       {\|\,u(\cdot,\cdot,t)\|_{2}},
\end{equation}
where the Euclidean norm is taken over all spatial grid points.
In the following results, unless otherwise indicated, we report the mean over every time snapshot \(t_m\).

\textbf{Structural similarity index measure (SSIM).} The single-scale SSIM between a prediction and the reference is
\[
\text{SSIM}(t)
=\frac{(2\mu_u\mu_{\tilde u}+C_1)\,(2\sigma_{u\tilde u}+C_2)}
        {(\mu_u^{2}+\mu_{\tilde u}^{2}+C_1)\,(\sigma_u^{2}+\sigma_{\tilde u}^{2}+C_2)},
\]
where \(\mu_u=\mathbb E[u]\) and \(\mu_{\tilde u}=\mathbb E[\tilde u]\) are local pixel-sample means, \(\sigma_u^{2}=\mathbb E[(u-\mu_u)^2]\) and \(\sigma_{\tilde u}^{2}=\mathbb E[(\tilde u-\mu_{\tilde u})^2]\) are local pixel-sample variances, \(\sigma_{u\tilde u}=\mathbb E[(u-\mu_u)(\tilde u-\mu_{\tilde u})]\) is the local pixel-sample covariance, and \(C_1, C_2>0\) are small constants to stabilize division. 

\textbf{Spectral mean-squared error (Spectral MSE).} Denote the 2D Fourier amplitudes by
\(
\widehat{u}_{k,\ell}
=\mathrm{DFT}[u]_{k,\ell},
\widehat{\tilde u}_{k,\ell}
=\mathrm{DFT}[\tilde u]_{k,\ell}.
\)
The spectral MSE is $\text{SpecMSE}(t)
=\frac{1}{N_f^{2}}\sum_{k,\ell}
\bigl|\,
|\widehat{\tilde u}_{k,\ell}|
-|\widehat{u}_{k,\ell}|
\bigr|^{2}$. Physically, this is the mean-squared difference between spectral magnitudes. It penalizes mismatches in the distribution of energy over wavenumbers.

\textbf{Pearson linear correlation \(r\).} If we flatten each field into a vector, the Pearson coefficient \(r\) is given as
\[
r(t)
=\frac{\sum_{i}(u_i-\bar u)\,(\tilde u_i-\bar{\tilde u})}
       {\sqrt{\sum_{i}(u_i-\bar u)^{2}}
        \sqrt{\sum_{i}(\tilde u_i-\bar{\tilde u})^{2}}},
\] where bars denote sample means. Higher values indicate higher linear correlation.

\section{One–step test error generalization}\label{subsec:gen}
Let \(\mathcal D\) denote the distribution of fine snapshots \(u\) used in Phase~2. Define the population risk
\begin{equation}\label{eq:pop-risk}
\mathcal R(\theta):=\mathbb E_{u\sim\mathcal D}\,
\bigl\|\mathcal P_f(\Delta t)\,u-\mathcal N_\theta\!\bigl[\mathcal P_c(\Delta t)\,Pu\bigr]\bigr\| ,
\end{equation}
and let \(\widehat{\mathcal R}_m(\theta)\) be the empirical counterpart computed on $m$ approximately independent per-step pairs.

Write the loss on one instance as
\(f_\theta(u):=\bigl\|\mathcal P_f(\Delta t)u-\mathcal N_\theta\!\bigl[\mathcal P_c(\Delta t)Pu\bigr]\bigr\|\in[0,B]\)
and define the function class \(\mathcal F:=\{f_\theta:\theta\in\Theta\}\). Here $\Theta$ is the space of all feasible parameter vectors $\theta$.
Given a sample \(S=(u^{(1)},\dots,u^{(m)})\) drawn i.i.d.\ from \(\mathcal D\), denote the empirical mean
\(\widehat{\mathbb E}_S[f]=\tfrac1m\sum_{i=1}^m f(u^{(i)})\) and the population mean \(\mathbb E[f]=\mathbb E_{u\sim\mathcal D}[f(u)]\).
Define the empirical Rademacher complexity as
\begin{align*}
\widehat{\mathfrak R}_m(\mathcal F) &:=
\mathbb E_{\boldsymbol\sigma}\!\Biggl[\sup_{f\in\mathcal F}\frac{1}{m}\sum_{i=1}^m \sigma_i f(u^{(i)})\Biggr]\\
\sigma_i &\stackrel{\text{i.i.d.}}{\sim}\mathrm{Unif}\{-1,+1\},
\end{align*}
and, for brevity, write \(\mathfrak R_m:=\widehat{\mathfrak R}_m(\mathcal F)\).  Note that $\mathfrak R_m$ is a data–dependent measure of how well functions in \(\mathcal F_\theta\) can fit random \(\pm1\) noise. While such complexity measures are typically used to derive generalization bounds for i.i.d. classification problems, in the RELift setting---where we learn continuous-valued PDE dynamics from temporally correlated trajectories—--$\mathfrak R_m$ should be viewed only as a qualitative proxy for the capacity of the RELift hypothesis class, rather than as a fully rigorous complexity measure.

Suppose the loss $f_{\theta}$ is bounded by \(B>0\). As we detail below, standard empirical risk minimization (via symmetrization and contraction) implies that with probability at least \(1-\delta\),
\begin{equation}\label{eq:gen-bound}
\mathcal R(\theta)
\leq
\widehat{\mathcal R}_m(\theta)
+2\,\mathfrak R_m
+3B\sqrt{\frac{\ln(2/\delta)}{2m}}
\end{equation}
This bound is stated for independently, identically distributed (i.i.d.) one-step pairs $(u^n,u^{n+1})$. In practice, we subsample sequential steps from trajectories; thus the reported empirical results should be viewed as supportive evidence rather than a strict proof where the data contains temporal dependence. Therefore, the \text{per–step test error} is controlled by the observed Phase~2 training error plus a capacity/finite–sample term decaying as \(m^{-1/2}\). 

\paragraph{\textbf{Derivation.}}
\emph{Uniform symmetrization. }
Consider the deviation functional
\(\Phi(S):=\sup_{f\in\mathcal F}\bigl(\mathbb E[f]-\widehat{\mathbb E}_S[f]\bigr)\).
Introducing a ghost sample \(S'=(u'^{(1)},\dots,u'^{(m)})\) i.i.d.\ as \(S\) and using standard symmetrization yields,
\begin{align*}
\mathbb E_S[\Phi(S)]
&=\mathbb E_{S}\Bigl[\sup_{f\in\mathcal F}\bigl(\mathbb E_{u'}[f(u')]-\widehat{\mathbb E}_S[f]\bigr)\Bigr] \\
&\le \mathbb E_{S,S'}\Bigl[\sup_{f\in\mathcal F}\bigl(\widehat{\mathbb E}_{S'}[f]-\widehat{\mathbb E}_{S}[f]\bigr)\Bigr] \\
&= \mathbb E_{S,S',\boldsymbol\sigma}\Bigl[\sup_{f\in\mathcal F}\frac{1}{m}\sum_{i=1}^m \sigma_i\bigl(f(u'^{(i)})-f(u^{(i)})\bigr)\Bigr] \\
&\le 2\,\mathbb E_{S,\boldsymbol\sigma}\Bigl[\sup_{f\in\mathcal F}\frac{1}{m}\sum_{i=1}^m \sigma_i f(u^{(i)})\Bigr]
=2\,\mathfrak R_m.
\end{align*}

\emph{Bounded differences concerntration. }
By McDiarmid’s inequality, with probability at least \(1-\delta\),
\[
\Phi(S)\le \mathbb E_S[\Phi(S)]+B\sqrt{\frac{\ln(2/\delta)}{2m}}.
\]
Combining with the symmetrization argument yields
\[
\sup_{f\in\mathcal F}\bigl(\mathbb E[f]-\widehat{\mathbb E}_S[f]\bigr)\le
2\,\mathfrak R_m
+ B\sqrt{\frac{\ln(2/\delta)}{2m}}
\quad\text{w.p.\ }\ge 1-\delta.
\]

Applying the same bounded-differences argument to the random functional
\(S\mapsto \mathfrak R_m=\widehat{\mathfrak R}_m(\mathcal F)\)
and a union bound yields, with probability at least \(1-\delta\),
\[
\forall\,\theta\in\Theta:\quad
\mathcal R(\theta)-\widehat{\mathcal R}_m(\theta)
\le
2\,\mathfrak R_m
+3B\sqrt{\frac{\ln(2/\delta)}{2m}},
\]
which is exactly \eqref{eq:gen-bound}.

\begin{figure}[t]
  \centering
  \includegraphics[width=\linewidth]{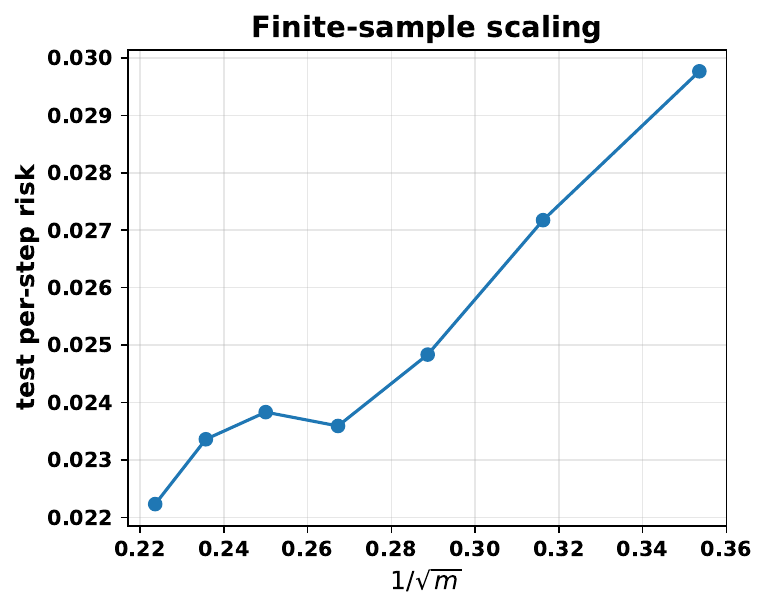}
  \vspace{-0.6em}
  \caption{\justifying \textbf{Finite sample scaling of the Phase~2 per–step test risk.}
  For random subsamples of size $m$ from the test set, we compute the empirical per–step risk
  $\widehat{\mathcal R}_m(\theta)=\tfrac{1}{m}\sum_{i=1}^m
  \|\mathcal P_f(\Delta t)u^{(i)}-\mathcal N_\theta[\mathcal P_c(\Delta t)P u^{(i)}]\|$
  and plot it against $1/\sqrt{m}$. The nearly linear trend is consistent with the generalization bound in~\eqref{eq:gen-bound}, which give a finite sample  decay as $m^{-1/2}$.
  }
  \label{fig:finite-sample-scaling}
\end{figure}

\paragraph{\textbf{Finite sample behavior of the test risk.}}
To empirically validate the generalization prediction in~\eqref{eq:gen-bound} from multiple test trajectories, we randomly subsample $m$ one-step pairs $(u^n, u^{n+1})$, evaluate $\|u^{n+1} - \widehat u^{n+1}\|$, average over the $m$ samples, and plot the mean against $1/\sqrt m$. A near linear increase in Figure~\ref{fig:finite-sample-scaling} versus $1/\sqrt m$ matches the $O(m^{-1/2})$ finite-sample result in the generalization bound \eqref{eq:gen-bound}. Although per-step samples along a trajectory are dependent, the observed near-linear trend versus $1/\sqrt{m}$ aligns with the i.i.d.\ prediction and is robust under random subsampling. As expected, we find that larger $m$ reduces variance in the empirical per-step~risk.

\section{Details for 1D1V VP data  generation}\label{sec:datagen-vp}
For the 1D1V VP time integrator, we use a second-order Lax-Wendroff (LW) scheme. Compared to Runge-Kutta 4 (RK4), we find that the LW scheme better preserves the phase-space structure in the long-time limit. For a scalar advection equation,
\(\partial_t f + a\,\partial_x f = 0\), the LW update is given by,
\[
  f^{\,n+1}_i
  =
  f^{\,n}_i
  - \frac{a\,\Delta t}{2\Delta x}
    \bigl(f^{\,n}_{i+1}-f^{\,n}_{i-1}\bigr)
  + \frac{a^2\Delta t^2}{2\Delta x^2}\bigl(f^{\,n}_{i+1}-2f^{\,n}_i+f^{\,n}_{i-1}\bigr).
\]

In this work, we consider a (2 + 1)-dimensional system with one real spatial component and one velocity spatial component (1D1V). We study two types of systems that arise as a result of changing the initial distribution of electrons: (i) randomized electron distribution; and (ii) the two-stream instability. For the randomized electron distribution,  we fix the  ions (see Appendix~\ref{sec:ions}) with a Maxwell-Boltzmann distribution given by,
\[
f_i(x,v,0)
= \frac{1}{\sqrt{2\pi\,v_{\rm th,i}^2}}
  \exp\!\Bigl(-\frac{v^2}{2\,v_{\rm th,i}^2}\Bigr),
\] and initialize the electrons with a randomized cosine-superposition given by,
\[
f_e(x,v,0)
= \bigl[\,
1 + \sum_{n_x,n_v=0}^{m_x-1,m_v-1}
\cos\bigl(2\pi(n_x\,x + n_v\,v) + \varphi_{n_x,n_v}\bigr)\bigr]^2.
\]
Here, $v_{\rm th,i}$ represents the thermal velocity (average speed of ions moving randomly within the plasma due to their temperature), the phases \(\varphi_{n_x,n_v}\) are drawn i.i.d.\ in \([0,2\pi)\), and the result is renormalized to unit density. In practice, we fix $m_x,m_v = \{(4,4), (8,8)\}$, and generate 1000 trajectories each by randomly sampling \(\varphi_{n_x,n_v}\).  Both datasets used uniform grid spacing where $x\in[-2\pi, 2\pi]$ and $v\in[-6,6]$. We refer to the dataset with $m_x,m_v=4,4$ as the low frequency dataset, and the dataset with $m_x,m_v=8,8$ as the high frequency dataset.  We do not adaptively timestep, as this poses some slight complications for Phase 2 training and future time extrapolation. We instead fix $\Delta t=0.01$, which satisfies the CFL condition, and propagate until physical time $T_{train} = 1.0$.  For the subsequent results, we assume a $32\times32$ grid for the coarse grain data generation and a $128\times128$ grid for the fine grain data generation.

For the two-stream dataset, the ions are taken to be stationary Maxwell-Boltzmann distributions in the same way used for the randomized dataset. The electron distribution is a pair of counter-propagating beams whose
density is perturbed by a small, broadband perturbation:
\begin{multline*}
f_e(x,v,0)
=
\Bigl[\,1+A_1\cos(kx+\phi_1)+A_2\cos(2kx+\phi_2)\Bigr]\,
\frac{1}{2\sqrt{2\pi}}
\\
\cdot
\Biggl[
e^{-\frac{(v-\bar v)^2}{2}}
+
e^{-\frac{(v+\bar v)^2}{2}}
\Biggr].
\end{multline*}

Here, \(k=\beta/\lambda_D\) with \(\lambda_D\) the Debye length,
\(A_1=\varepsilon\) is the fundamental-mode amplitude,
\(A_2=\varepsilon r\) is the second-harmonic amplitude, and
\(\bar v\) is the beam drift speed. Each trajectory is obtained by independently sampling the
parameter set
\begin{multline*}
\varepsilon=\varepsilon_0(1+\delta_\varepsilon),\quad
\delta_\varepsilon\sim\mathrm U[-0.10,0.10],\quad
\varepsilon_0=5\times10^{-3},
\\
\beta=\beta_0+\delta_\beta,\quad
\delta_\beta\sim\mathrm U[-0.05,0.05],\quad
\beta_0=0.20,
\\
\bar v=\bar v_0(1+\delta_v),\quad
\delta_v\sim\mathrm U[-0.05,0.05],\quad
\bar v_0=2.4,
\\
r\sim\mathrm U[0.05,0.10],\qquad
\phi_{1},\phi_{2}\sim\mathrm U[0,2\pi).
\end{multline*}

This ensemble therefore varies the perturbation strength, the effective
wavenumber, the beam drift, the second-harmonic content, and the phases. This creates a dataset that captures a wide range of growth and saturation behaviors. We generate 1000 trajectories for training and evolved each one on a coarse \(64\times64\) grid and a fine \(256\times 256\) grid using a fixed timestep $\Delta t=0.01$ and simulate until physical time $T_{\text{sim}} = 50$. During post-processing, we retain only the one-second slice $22 \leq t < 23 ,$
for both grids, resulting in 101 equally spaced frames. This captures the portion in time evolution when the electron two-stream instability has grown well enough into its nonlinear regime, but has yet to fully mix/form coherent vortices.

\section{Why freezing the ions is valid}\label{sec:ions}
A \emph{stationary-ion} (single-species) Vlasov--Poisson model is routinely adopted whenever the dynamics of interest unfold on electronic time-scales.  Canonical examples of interest from the plasma community include linear and nonlinear Landau damping \cite{baumann2021landau}. Landau damping is the collisionless attenuation of electrostatic waves due to phase-mixing in the \emph{electron} distribution function. This process unfolds on electronic time scales and therefore treats the far heavier ions as effectively immobile. More recently, data-driven studies of Vlasov--Poisson systems using latent-space identification or PINNs \cite{he2025physics,qin2023data} have also employed this approximation to make predictions with respect to the electron dynamics.  

The rationale is the large separation between the electron and ion plasma periods, quantified by the large discrepancies between their masses as 
\[
\frac{T_i}{T_e}=\sqrt{\frac{m_i}{m_e}},
\]
where \(m_e\) and \(m_i\) are the electron and ion masses and \(T_e\) and \(T_i\) their respective plasma periods.  In general, the plasma frequency of a species with number density \(n_s\), charge magnitude \(q_s\) (for electrons \(q_e = -1\)), and mass \(m_s\) is  
\[
\omega_{p,s}
  =\sqrt{\frac{n_s q_s^{2}}{\epsilon_{0} m_s}},
\]
so the period \(T_s = 2\pi/\omega_{p,s}\) scales as \(T_s \propto \sqrt{m_s}\); hence \(T_i/T_e = \sqrt{m_i/m_e}\), showing that heavier ions oscillate intrinsically more slowly than light electrons under the same electrostatic restoring force \cite{gibbon2020introduction}.  Even with a reduced mass ratio \(m_i/m_e = 25\) one obtains \(T_i= 5\,T_e\).  The electron period is \(T_e = 2\pi/\omega_{pe}\) with  
\[
\omega_{pe}
  =\sqrt{\frac{n_e e^{2}}{\epsilon_{0} m_e}},
\]
where \(n_e\) is the electron number density, \(e\) the elementary charge, and \(\epsilon_{0}\) the vacuum permittivity (in our normalized units \(n_e = 1\) and \(e = \epsilon_{0} = 1\), giving \(\omega_{pe}=1\)).  Over simulation windows of only a few electron periods, the majority of the ion velocity changes negligibly, so omitting the ion species from Vlasov--Poisson equation leaves the linear growth, nonlinear saturation, and phase-space structure of the electron species virtually unchanged.  A full two-species treatment becomes essential only when simulations extend well beyond \(T_i\) or deliberately excite coupled electron-ion modes such as the Buneman instability\cite{gibbon2020introduction, baumann2021landau}; such regimes lie outside the methodological scope of this work.

\section{Supplementary results}\label{sec:supresults}

\subsection{Additional ablation studies}\label{sec:ablation}
\squeezetable

\begin{table}[ht!]
\setlength{\tabcolsep}{4pt}
\renewcommand{\arraystretch}{0.95}
\centering
\scriptsize
\newcolumntype{b}{>{\columncolor{gray!8}}c}
\begin{tabular}{@{}l b c@{}}
\toprule
& \multicolumn{1}{c}{\cellcolor{gray!15}\textbf{Baseline}} 
& \multicolumn{1}{c}{\cellcolor{gray!6}\textbf{RELift model}} \\
\cmidrule(lr){2-2}\cmidrule(lr){3-3}
\textbf{Viscosity $\nu$} & \cellcolor{gray!15}\textbf{Upsampled} & \textbf{FUnet} \\
\midrule

$\nu = 10^{-3}$
  & $1.41\times10^{1}$
  & $\mathbf{5.64\times10^{-1}} \pm 6.12\times 10^{-2}$ \\[0.3em]

$\mathbf{\nu = 10^{-4}}$ 
  & $1.47\times10^{1}$
  & $\mathbf{1.25\times10^{-2}}\pm 4.82\times10^{-3}$ \\[0.3em]

$\nu = 10^{-5}$
  & $1.48\times10^{1}$
  & $\mathbf{1.47\times10^{-1}}  \pm 5.07\times 10^{-2}$ \\[0.3em]

\bottomrule
\end{tabular}
\caption{
\justifying
\textbf{ \justifying  Relative $L_2$ errors for Phase~1 Navier--Stokes super-resolution across different viscosity values.}
The shaded column corresponds to the bicubic upsampling baseline, while the remaining column reports the FUnet  based RELift model. The FUnet model here was trained on data using $\nu=10^{-4}$.}
\label{tab:phase1_ns_funet_multi_nu}
\end{table}

\begin{figure}[ht!]
\centering

\begin{subfigure}{0.95\linewidth}
    \centering
    \includegraphics[width=\linewidth]{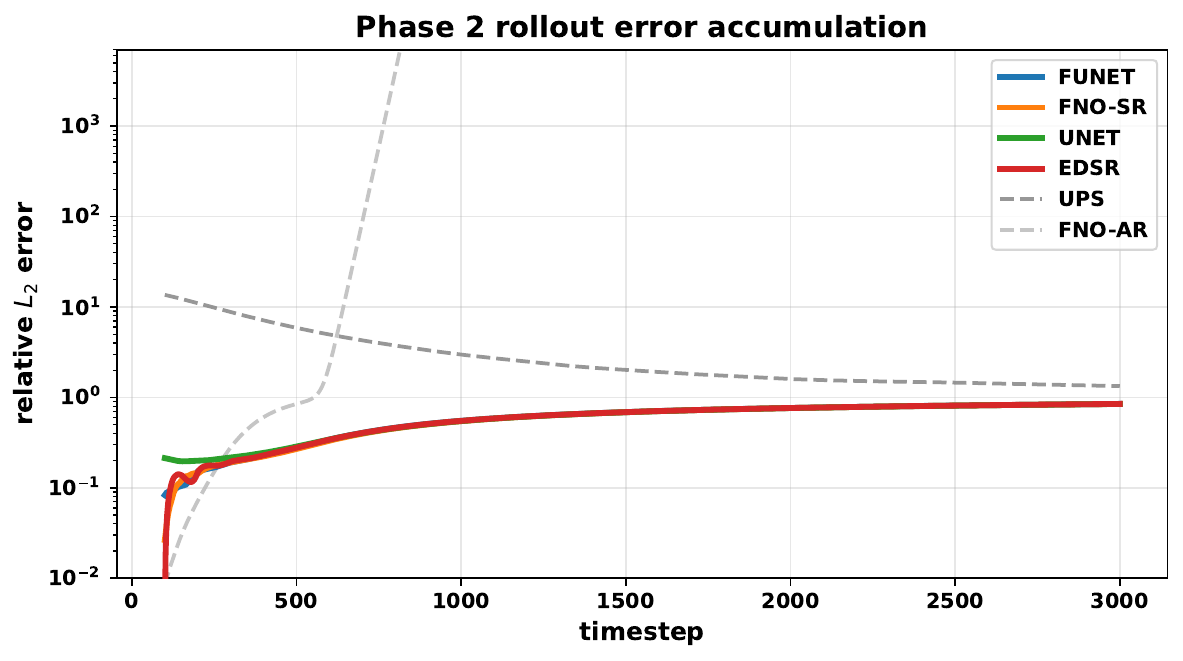}
    \caption{$\nu = 10^{-5}$}
    \label{fig:ns-relL2-time-nu1e-5}
\end{subfigure}

\vspace{0.6em}

\begin{subfigure}{0.95\linewidth}
    \centering
    \includegraphics[width=\linewidth]{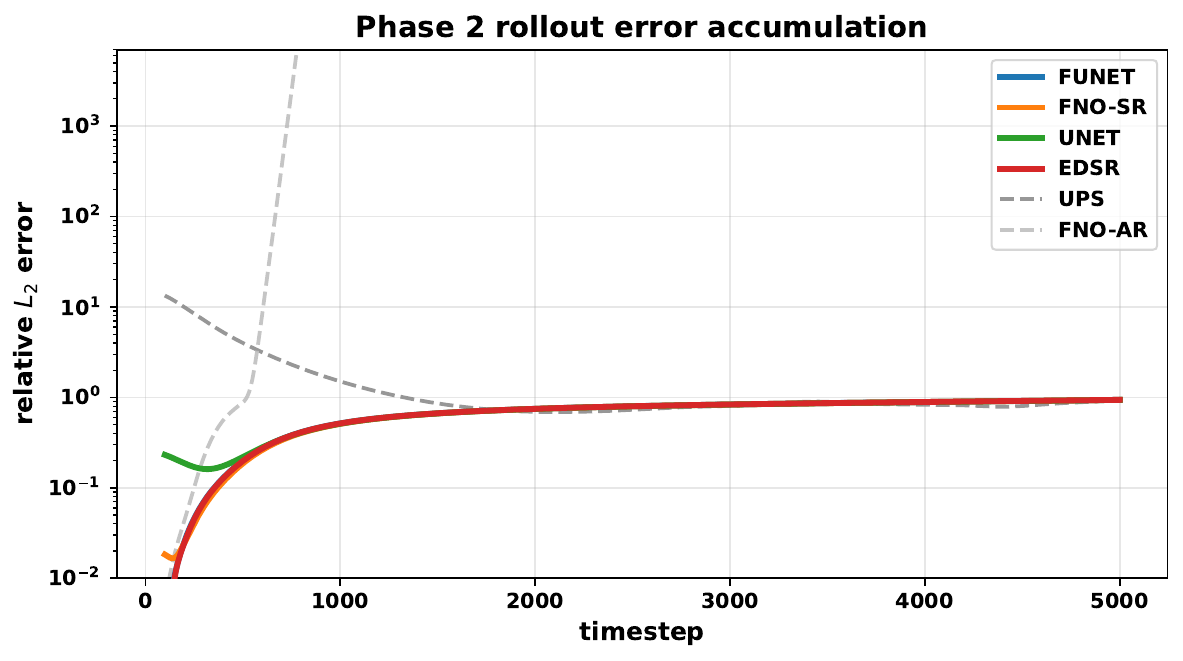}
    \caption{$\nu = 10^{-4}$}
    \label{fig:ns-relL2-time-nu1e-4}
\end{subfigure}

\vspace{0.6em}

\begin{subfigure}{0.95\linewidth}
    \centering
    \includegraphics[width=\linewidth]{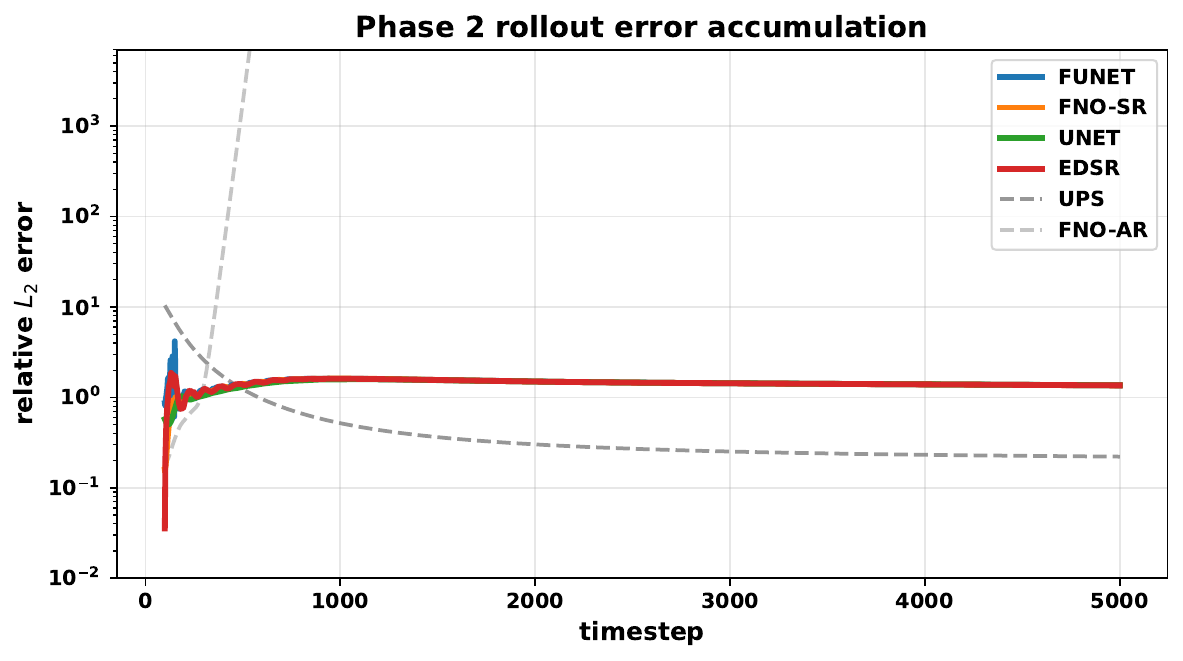}
    \caption{$\nu = 10^{-3}$}
    \label{fig:ns-relL2-time-nu1e-3}
\end{subfigure}

\caption{\justifying 
\textbf{Relative \(L_2\) error for a representative unseen test trajectory.}
At a fixed \(4\times\) downsampling factor, the RELift Phase~2 predictions (colored curves) are compared against the baselines (gray). The model used here was trained using Navier--Stokes data at \(\nu=10^{-4}\) only, and is then evaluated on unseen trajectories at \(\nu=10^{-5}\), \(\nu=10^{-4}\), and \(\nu=10^{-3}\). These results therefore show both performance at the training viscosity and generalization away from the training viscosity, while also illustrating how rollout error evolves beyond the training window of Phase~1.
}
\label{fig:ns-relL2-time-multi-nu}
\end{figure}

\begin{figure*}[ht!]
\centering

\begin{subfigure}{0.95\linewidth}
    \centering
    \includegraphics[width=\linewidth]{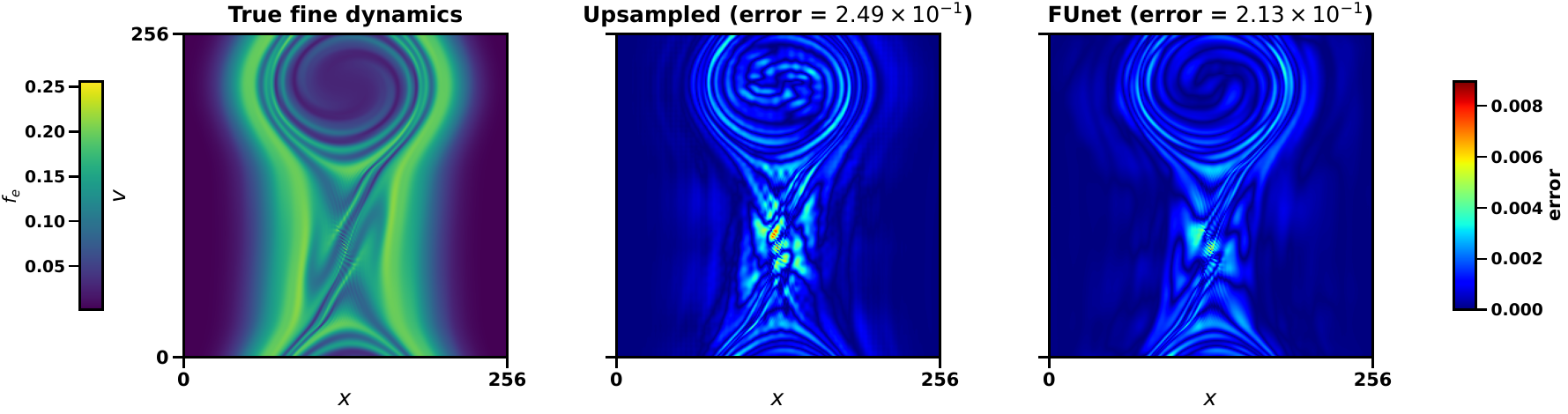}
    \caption{\justifying Two-stream instability at \(T=64\). In this case, the Phase~1 training window spans only \(T\in[22,23]\), so the displayed prediction corresponds to an extrapolation approximately \(4000\) timesteps, or \(41\) physical time units, beyond the end of the training window.}
    \label{fig:vp-two-strema-p2-t40}
\end{subfigure}

\vspace{0.6em}

\begin{subfigure}{0.95\linewidth}
    \centering
    \includegraphics[width=\linewidth]{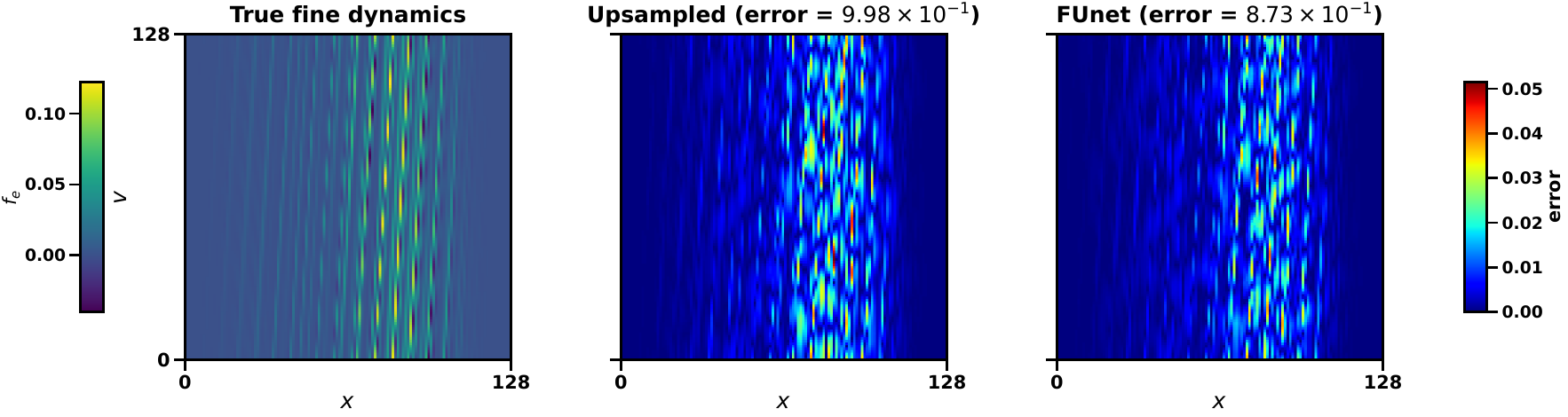}
    \caption{\justifying Randomized initial condition at \(T=20\). Here, the Phase~1 training window covers only \(T\in[0,1]\), so the shown result is an extrapolation \(2000\) timesteps, or \(20\) physical time units, beyond the end of the training interval.}
    \label{fig:vp-random-p2-t20}
\end{subfigure}

\caption{\justifying 
\textbf{Super long-time stability of RELift for future-time extrapolation in the 1D1V Vlasov--Poisson system.}
These examples illustrate that the learned effective fine-grid propagator remains stable and physically coherent far beyond the temporal window used in Phase~1 training. In both cases, the model is rolled out deep into the future, well outside the training regime, demonstrating that RELift can maintain accurate long-time predictions even when extrapolating thousands of timesteps beyond the observed training horizon.
}
\label{fig:vp-long-time}
\end{figure*}

\begin{figure*}[ht!]
    \centering
    \begin{subfigure}[b]{0.48\linewidth}
        \centering
        \includegraphics[width=\linewidth]
            {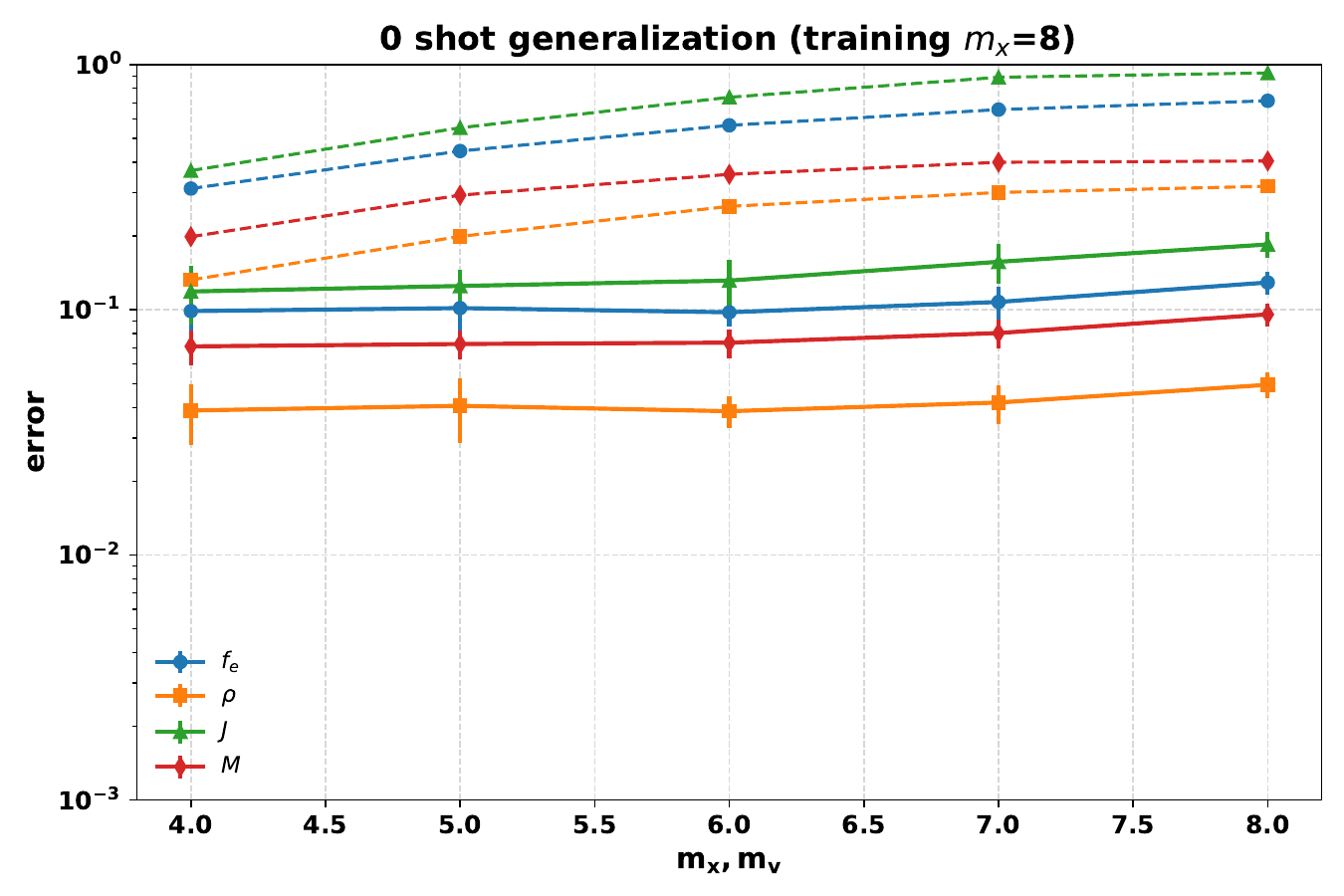}
        \caption{Training on the \emph{high-frequency} dataset 
                 $(m_x,m_v)=(8,8)$ and evaluating on progressively lower
                 test bandwidths.}
    \end{subfigure}
    \hfill
    \begin{subfigure}[b]{0.48\linewidth}
        \centering
        \includegraphics[width=\linewidth]
            {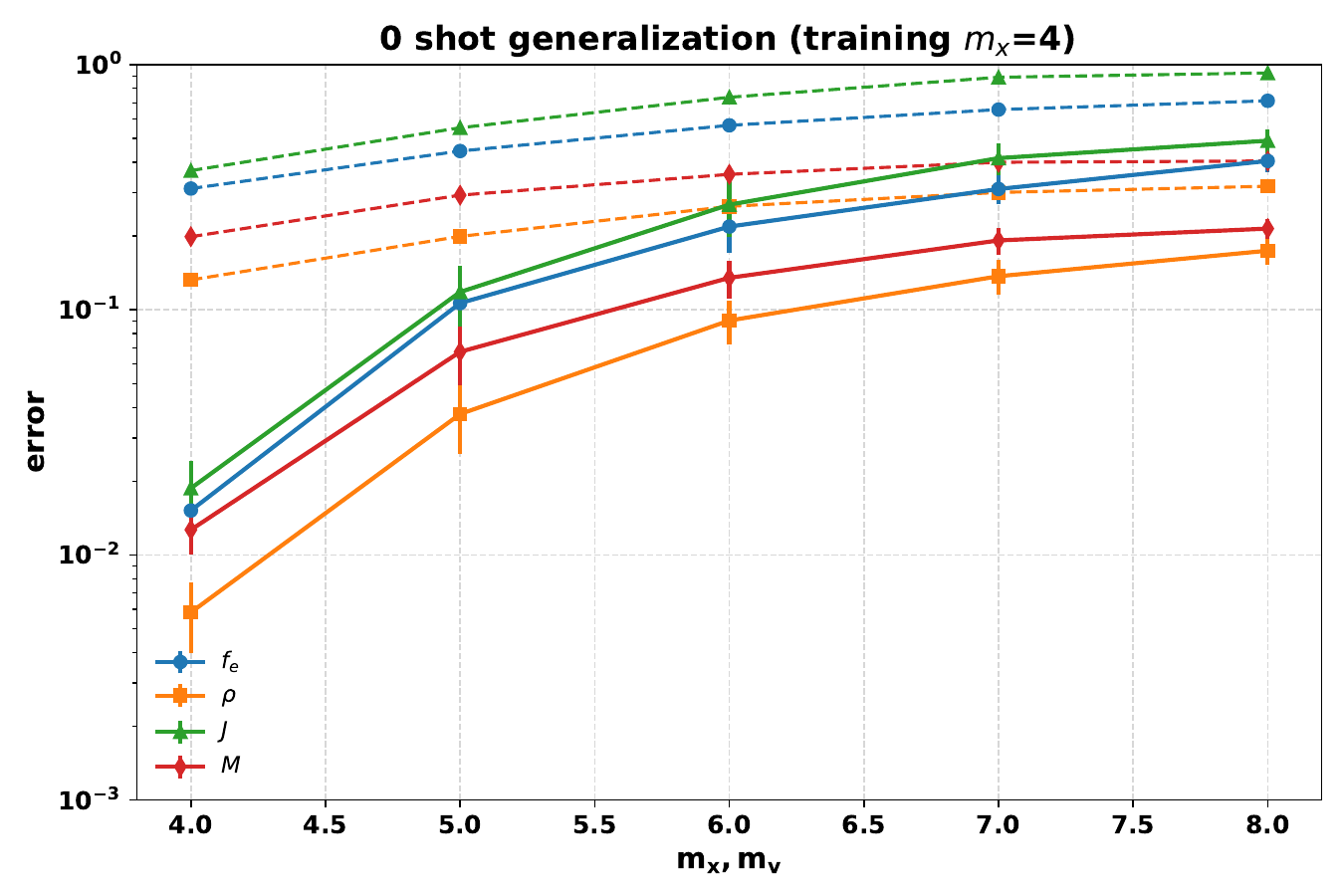}
        \caption{Training on the \emph{low-frequency} dataset 
                 $(m_x,m_v)=(4,4)$ and evaluating on progressively higher
                 test bandwidths.}
    \end{subfigure}
    \caption{\justifying  \textbf{Zero-shot generalization sweeps for the 1D1V Vlasov--Poisson super-resolution task with random cosine superposition initial condition.}  Solid curves show the FUnet prediction
errors while dashed curves indicate the bicubic
upsampling baseline.  
High-frequency training (left) maintains low error when asked to resolve coarser spectra, whereas the converse (right) incurs rapidly growing error when tested with high-frequency data.}
    \label{fig:vlasov-mx-sweeps}
\end{figure*}
\subsection{Training frequency transfer}
The left panel of Figure \ref{fig:vlasov-mx-sweeps} shows a sweep of test set errors where $m_x=m_v\in\{4,\ldots,8\}$ using a model trained \textit{only} on the highest band-limit \((8,8)\). The dashed curves correspond to the bicubic interpolation baseline errors, and the solid curves represent the FUnet predictions errors. Surprisingly, the error appears to decrease as the gap between train- and test-frequencies widen. Nevertheless, we see that for all predictions of the flow and moments, that the FUnet predictions retain anywhere from a $7-10\times$ improvement in error over the baseline,  demonstrating the non-trivial generalization to out of distribution spectra. The right panel shows a model that is trained only on the lowest band-limit \((4,4)\). As expected, the error increases as the frequency band-width increases.

\begin{table*}[ht!]
\setlength{\tabcolsep}{2pt}        
\renewcommand{\arraystretch}{0.95} 
\centering
\scriptsize
\begin{adjustbox}{max width=\textwidth}
\begin{tabular}{@{}l l >{\columncolor{gray!10}}c c@{}}
\toprule
& & \multicolumn{1}{c}{\cellcolor{gray!20}\textbf{Baseline}} & \multicolumn{1}{c}{\cellcolor{gray!8}\textbf{RELift model}} \\
\cmidrule(lr){3-3}\cmidrule(lr){4-4}
\textbf{Moment} & \textbf{$(m_x,m_v)^{\text{train}}\!\to\!(m_x,m_v)^{\text{test}}$}
& \cellcolor{gray!20}\textbf{Bicubic} & \textbf{FUnet} \\
\midrule

\multirow[t]{3}{*}{$\rho$}
  & $(4,4)\!\to\!(4,4)$ & $1.32\!\times\!10^{-1}$ & $\mathbf{6.42\!\times\!10^{-3}}\pm 1.52\!\times\!10^{-3}$ \\[0.15em]
  & $(8,8)\!\to\!(8,8)$ & $3.19\!\times\!10^{-1}$ & $\mathbf{4.81\!\times\!10^{-2}}\pm 5.83\!\times\!10^{-3}$ \\[0.15em]
  & $(8,8)\!\to\!(4,4)$ & $1.32\!\times\!10^{-1}$ & $\mathbf{3.58\!\times\!10^{-2}}\pm 1.01\!\times\!10^{-2}$ \\[0.3em]
\midrule

\multirow[t]{3}{*}{$J$}
  & $(4,4)\!\to\!(4,4)$ & $3.69\!\times\!10^{-1}$ & $\mathbf{1.87\!\times\!10^{-2}}\pm 4.89\!\times\!10^{-3}$ \\[0.15em]
  & $(8,8)\!\to\!(8,8)$ & $9.25\!\times\!10^{-1}$ & $\mathbf{1.77\!\times\!10^{-1}}\pm 1.99\!\times\!10^{-2}$ \\[0.15em]
  & $(8,8)\!\to\!(4,4)$ & $3.69\!\times\!10^{-1}$ & $\mathbf{1.08\!\times\!10^{-1}}\pm 3.00\!\times\!10^{-2}$ \\[0.3em]
\midrule

\multirow[t]{3}{*}{$M$}
  & $(4,4)\!\to\!(4,4)$ & $1.98\!\times\!10^{-1}$ & $\mathbf{1.28\!\times\!10^{-2}}\pm 3.53\!\times\!10^{-3}$ \\[0.15em]
  & $(8,8)\!\to\!(8,8)$ & $4.04\!\times\!10^{-1}$ & $\mathbf{9.17\!\times\!10^{-2}}\pm 9.52\!\times\!10^{-2}$ \\[0.15em]
  & $(8,8)\!\to\!(4,4)$ & $1.98\!\times\!10^{-1}$ & $\mathbf{6.49\!\times\!10^{-2}}\pm 1.08\!\times\!10^{-2}$ \\
\bottomrule
\end{tabular}
\end{adjustbox}

\caption{\justifying
\textbf{1D1V Vlasov Phase 1 test set relative $L_2$ errors.}
Errors are averaged over 20 new initial conditions. FUnet results are averaged over five training runs and reported as mean \(\pm\) standard deviation, while bicubic upsampling is deterministic and reported as a mean only. Lower values are highlighted in bold.}
\label{tab:phase1_moments}
\end{table*}

\begin{figure*}[t!]
  \centering
  \includegraphics[width=\textwidth]{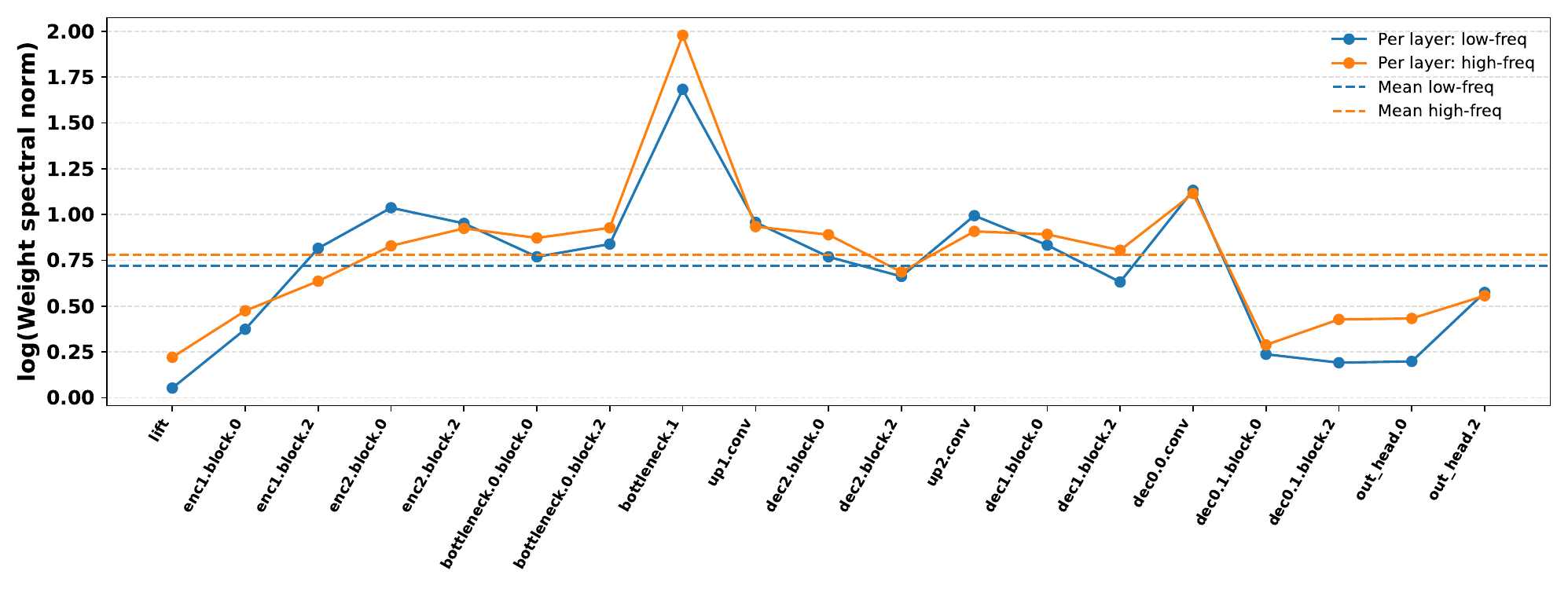}
  \caption{\justifying  \textbf{Per-layer spectral norms.} Each data point is the largest singular value \(\|W\|_2\) of a weight matrix in the model. Larger values therefore signal a greater ability to pass or create high-frequency features. Both models exhibit to clear spikes: one at the Fourier mixing layer \texttt{bottleneck.1} and one at the first decoder convolution \texttt{dec0.0.conv}. These are precisely the layers that must inject or preserve fine-scale modes as the feature map is downsampled and then re-expanded.}
  \label{fig:spec-norms}
\end{figure*}

To understand the transferability, we investigate the weights of the model for insight. Figure \ref{fig:spec-norms} shows that every weight matrix in the high-frequency network (trained on $(m_x, m_v)=(8,8)$ data) has a spectral norm that matches or exceeds its low-frequency analogue in key layers that are responsible for generating fine-scale features. The gap is noticeably the largest in the global Fourier layer
 \texttt{bottleneck.1}, which is primarily responsible for processing the mixed Fourier modes. Thus, we can reason that the high-frequency model retains more coarse directions, whilst additionally activating higher-frequency degrees of freedom that do not appear during $(m_x,m_v)=(4,4)$ training.

Let $f_{\theta}$ be a trained network with parameters
$\theta$ and select a single weight matrix $W_{\ell}\subset\theta$. Given $N$ coarse snapshots $\{x_i\}_{i=1}^{N}$ we form for each $i$, 
\[
  g_i =
  \frac{\partial}{\partial W_{\ell}}
  \Bigl[\,
    \mathcal{L}\bigl(f_{\theta}(x_i)\bigr)
  \Bigr]
  \in\mathbb{R}^{P},
  \qquad
  P = \lvert W_{\ell}\rvert ,
\]
for $\mathcal L$ the loss. If we then stack the parameter gradients row-wise, we obtain the Jacobian, $J=[g_1^{\top}\dot sg_N^{\top}]^{\top}\in\mathbb{R}^{N\times P}$, and the \emph{layer-restricted neural tangent kernel (NTK)} \cite{jacot2018neural} is
\[
  K = J\,J^{\top}\in\mathbb{R}^{N\times N},
  \quad
  K_{ij} = \langle g_i,\,g_j\rangle .
\]
Hence $K$ measures how similar the weight update would be if the network were (re-)trained on snapshots $x_i$ and $x_j$. We will make use of this to analyze two NTK diagnostics \cite{tancik2020fourier}:
\begin{enumerate}
    \item \textbf{ Eigenspectrum.}
The sorted eigenvalues $\lambda_1\geq\cdots\geq\lambda_N$ of $K$ help quantify the effective dimensionality of the subspace generated by the span of the gradients: a rapid decay will indicate that only a few directions in the parameter space are being adjusted during training, whereas a flatter tail indicates a richer set of learnable directions.

\item \textbf{ Kernel heatmap.} Visualizing $K$ itself highlights whether gradients are independent or redundant (i.e., $K$ is diagonal dominated or has strong off-diagonal blocks). A nearly diagonal kernel implies that different coarse snapshots possess distinct parameter updates, thus are likely to be well represented without the need for fine-tuning.
\end{enumerate}

\begin{figure*}[t!]
  \centering
  \begin{subfigure}[b]{0.48\linewidth}
    \includegraphics[width=\linewidth]{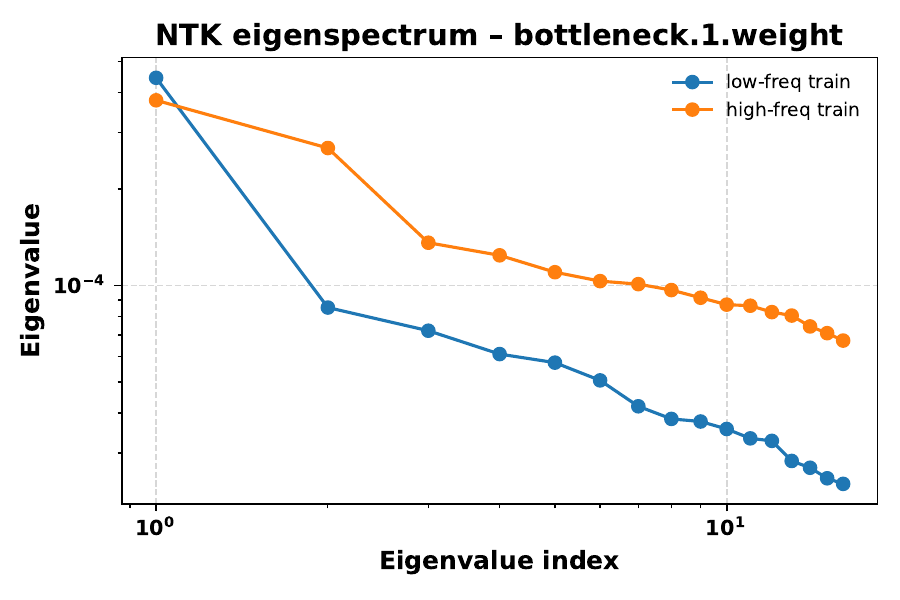}\\[6pt]
    \includegraphics[width=\linewidth]{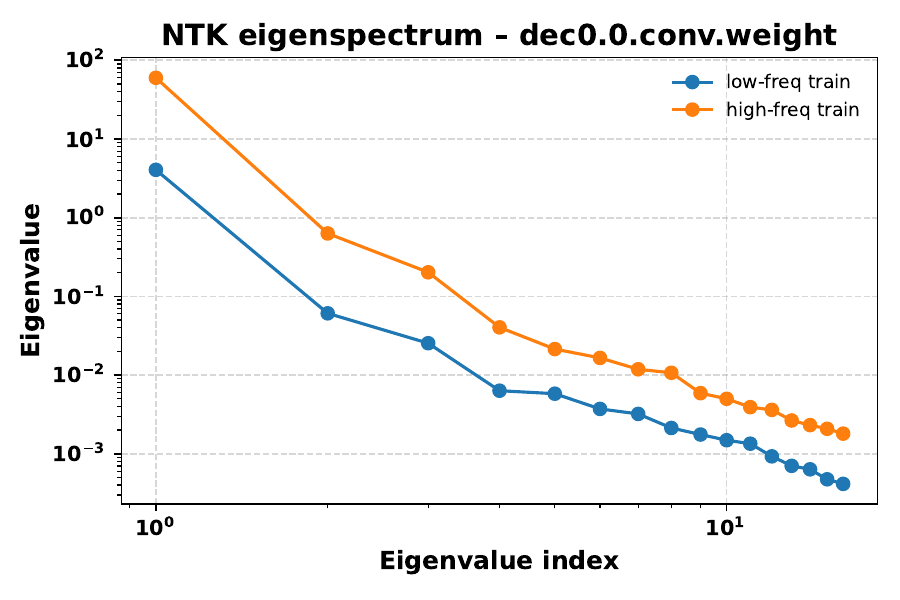}
    \caption{NTK eigenspectra}
  \end{subfigure}\hfill
  \begin{subfigure}[b]{0.48\linewidth}
    \includegraphics[width=\linewidth]{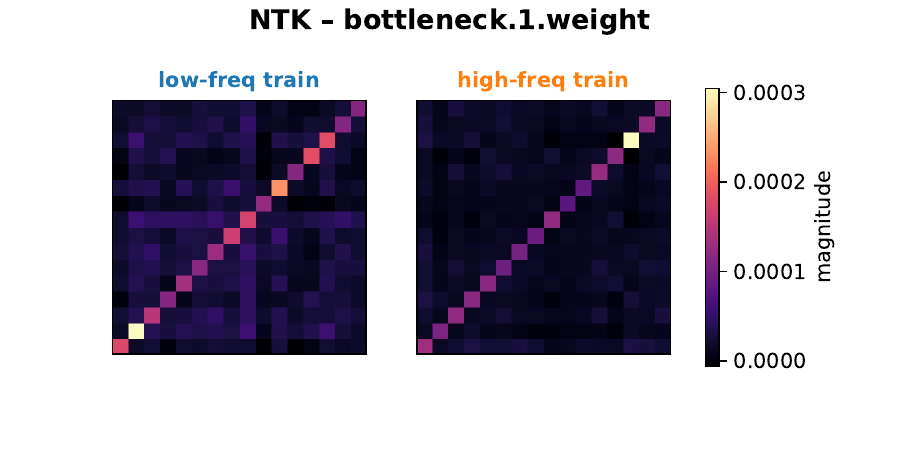}\\[6pt]
    \includegraphics[width=\linewidth]{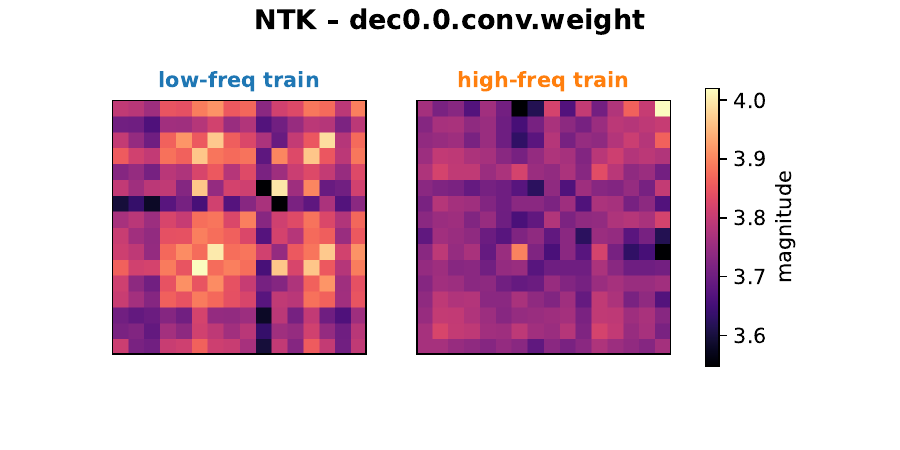}
    \caption{NTK heatmaps}
  \end{subfigure}
  \caption{\justifying  \textbf{Neural tangent kernel statistics}
  \emph{Top row}: Fourier mixing layer \texttt{bottleneck.1};
  \emph{bottom row}: first decoder convolution \texttt{dec0.0.conv}.
  High-frequency training (orange) yields richer eigenspectra and nearly
  diagonal kernels, whereas low-frequency training (blue) collapses to a
  low-rank span with shared gradients.}
  \label{fig:ntk-spec-eig}
\end{figure*}

The NTK diagnostics in Figure \ref{fig:ntk-spec-eig} corroborate the viewpoint that high-frequency training can generalize to low-frequency testing, but not vice-versa. For the bottleneck Fourier layer and the first decoder convolution, the eigenspectrum of the high-frequency training remains order(s) of magnitude larger beyond the first few ranks and its kernel is nearly diagonal. This indicates that for distinct snapshots, there are linearly independent gradients. Conversely, the low-frequency training eigenspectra collapse after rank $k\approx 4$, and the kernels have pronounced off-diagonal blocks. This signals that many input directions share identical or near identical curvature. These kernel-level differences explain the performance sweep in Figure \ref{fig:vlasov-mx-sweeps}: a new low-frequency snapshot will likely lie inside the curvature subspace spanned by the high-frequency NTK and is therefore recovered immediately in evaluation, whereas a high-frequency snapshot falls outside of the low-frequency span, so the model fails to generalize.


\begin{figure*}[htbp]
    \centering

    \begin{minipage}[b]{\linewidth}
        \centering
        \includegraphics[width=\linewidth]{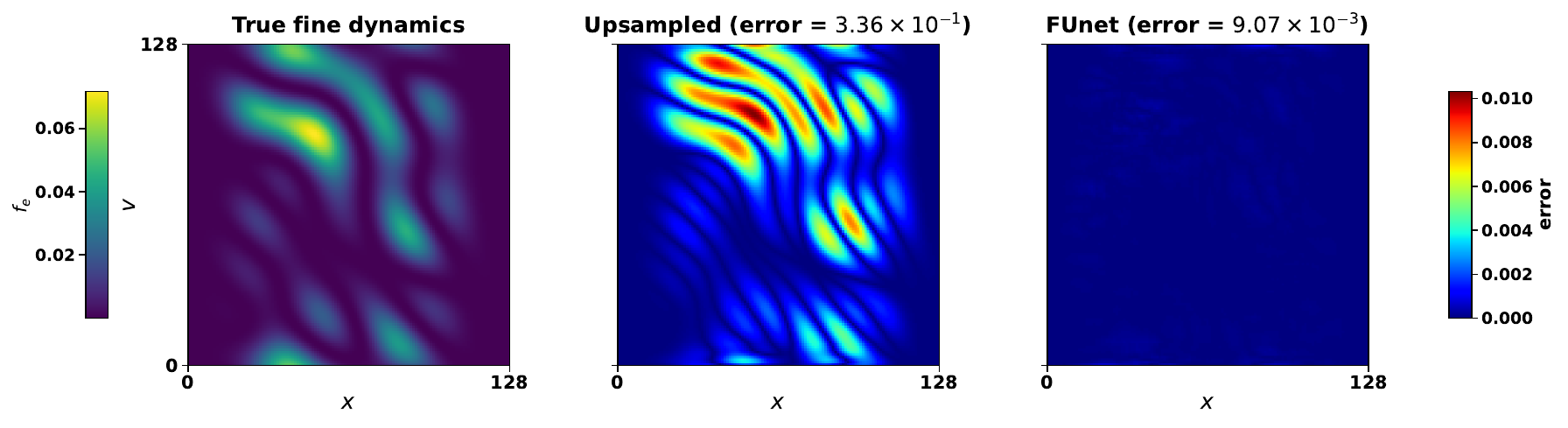}
    \end{minipage}

    \vspace{1em}
    \begin{minipage}[b]{\linewidth}
        \centering
        \includegraphics[width=\linewidth]{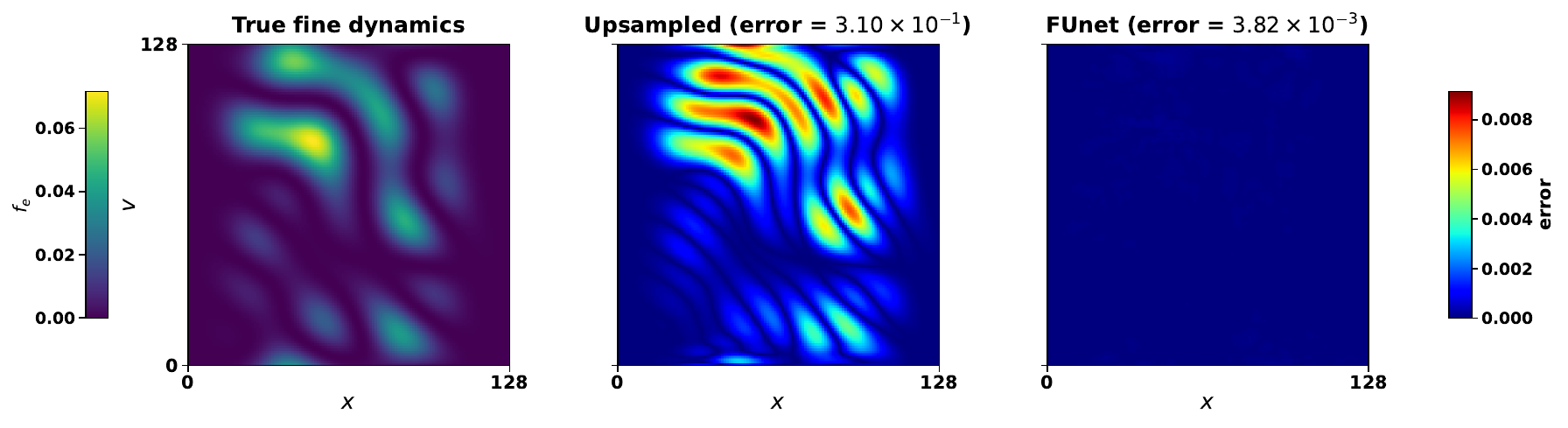}
    \end{minipage}

    \vspace{1em}
        \begin{minipage}[b]{\linewidth}
        \centering
        \includegraphics[width=\linewidth]{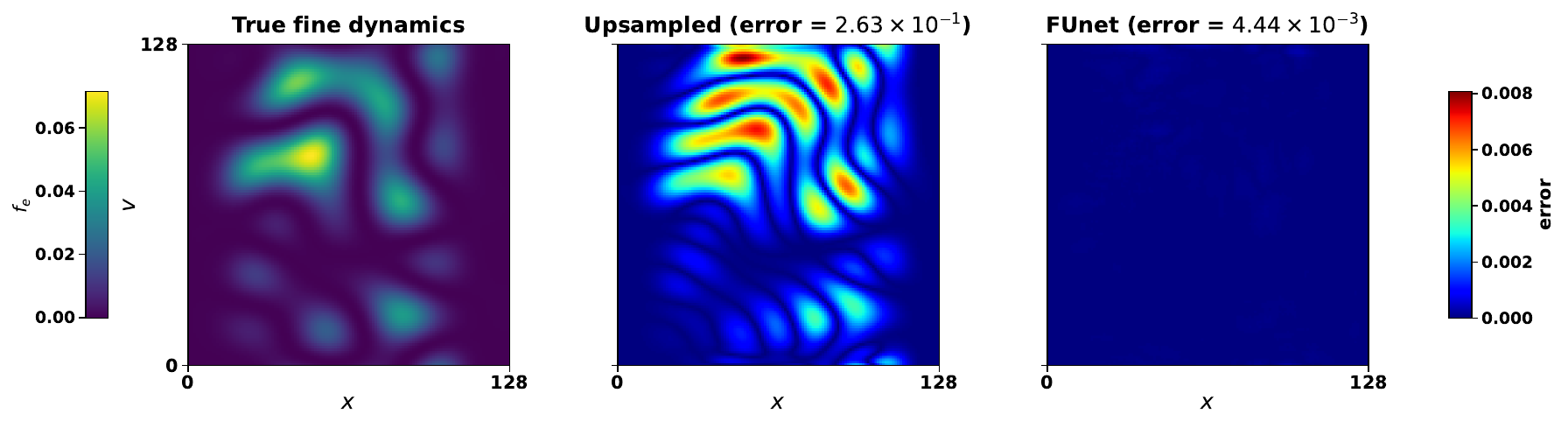}
    \end{minipage}

    \vspace{1em}
        \begin{minipage}[b]{\linewidth}
        \centering
        \includegraphics[width=\linewidth]{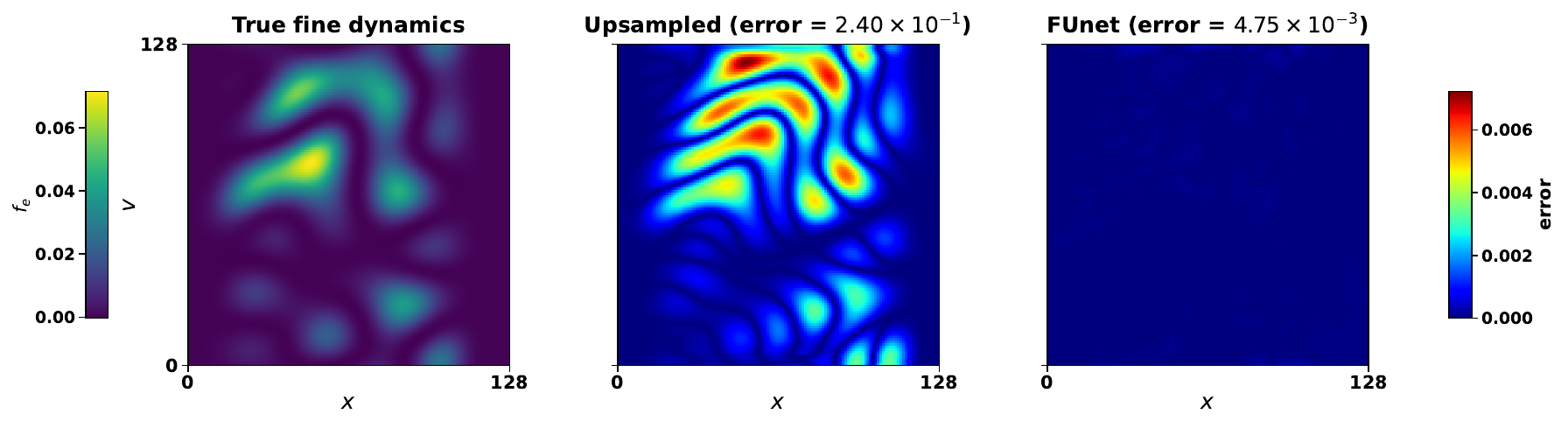}
    \end{minipage}
    \caption{\justifying  \textbf{Total Phase~1 super-resolution for randomized dataset.}}
    \label{fig:vlasov_phase1_randomized_full}
\end{figure*}

\begin{figure*}[htbp]
    \centering

    \begin{minipage}[b]{\linewidth}
        \centering
        \includegraphics[width=\linewidth]{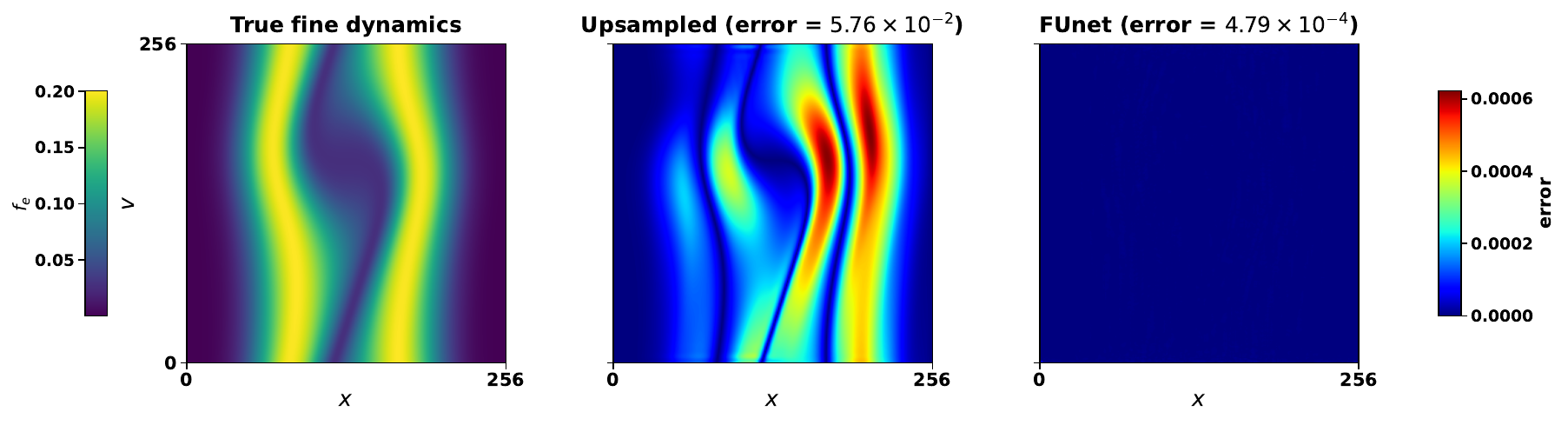}
    \end{minipage}

    \vspace{1em}
    \begin{minipage}[b]{\linewidth}
        \centering
        \includegraphics[width=\linewidth]{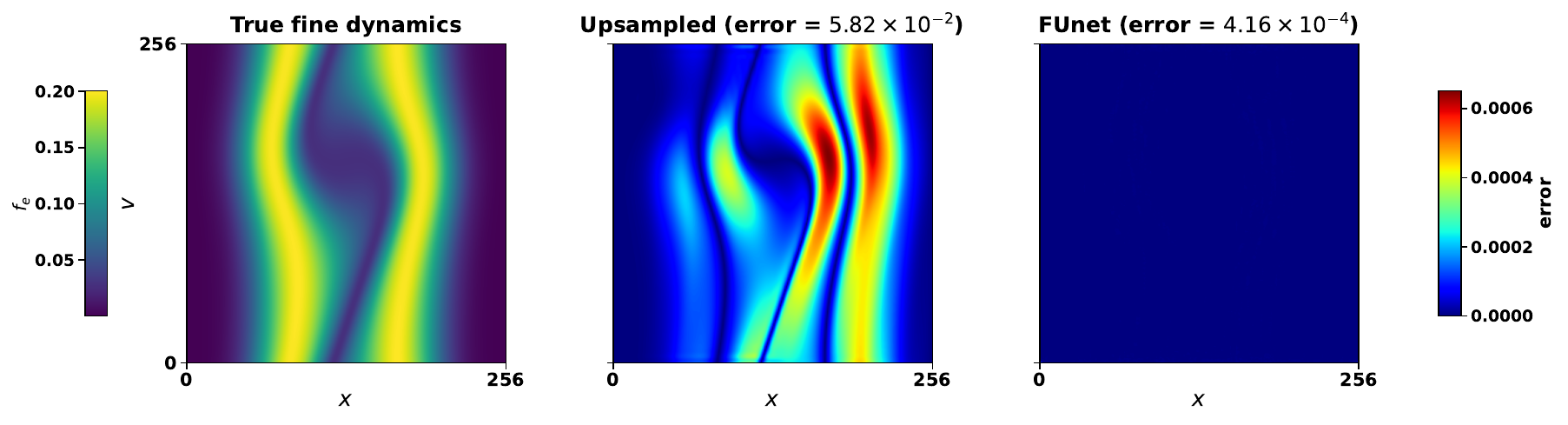}
    \end{minipage}

    \vspace{1em}
        \begin{minipage}[b]{\linewidth}
        \centering
        \includegraphics[width=\linewidth]{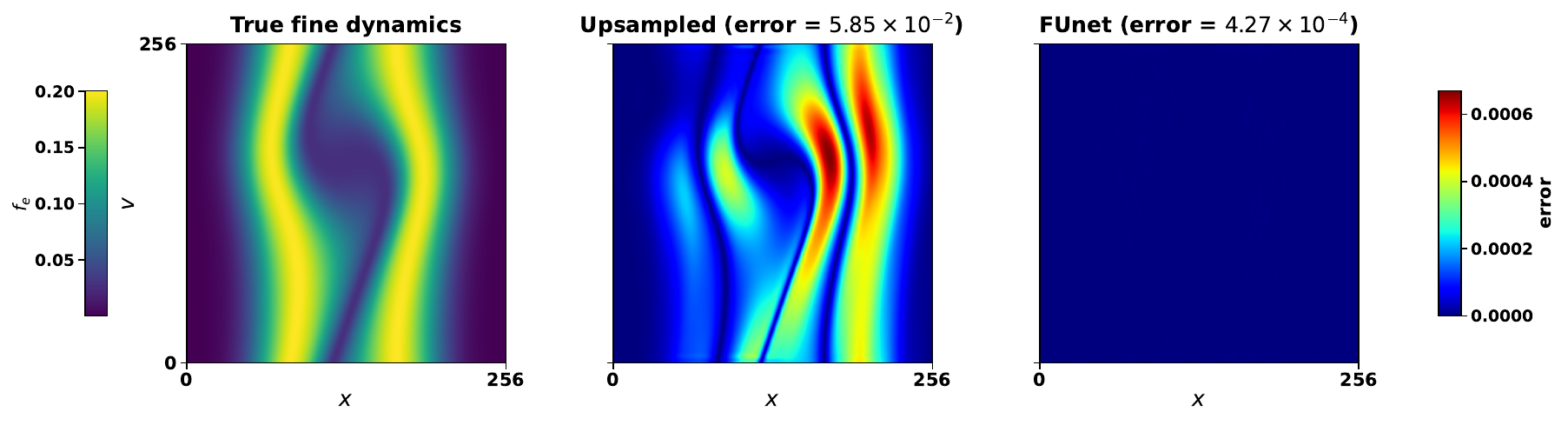}
    \end{minipage}

    \vspace{1em}
        \begin{minipage}[b]{\linewidth}
        \centering
        \includegraphics[width=\linewidth]{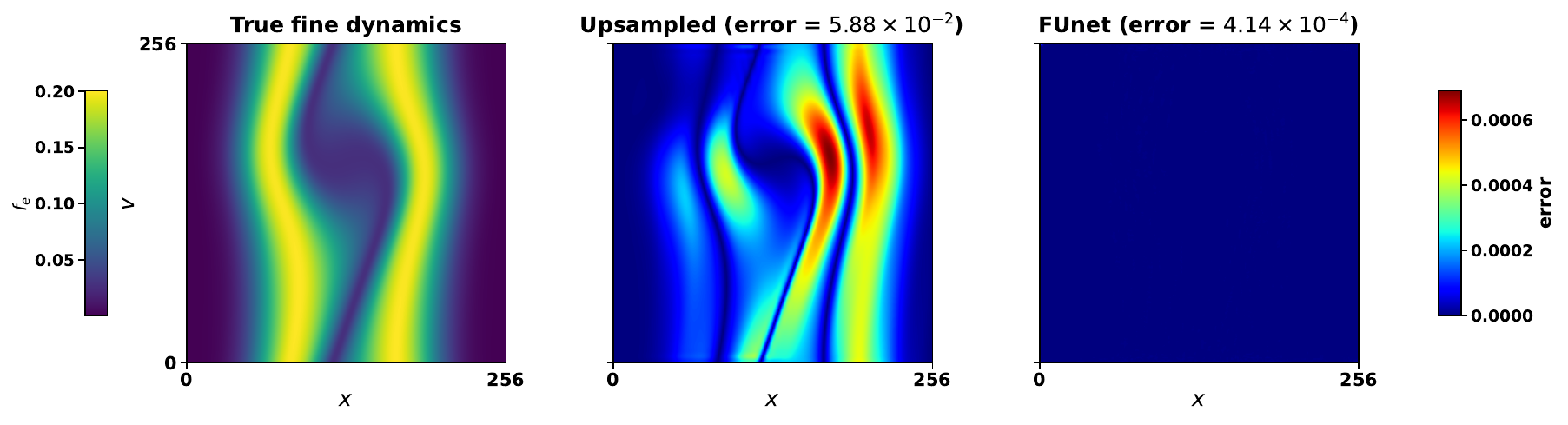}
    \end{minipage}
    \caption{\justifying  \textbf{Total Phase~1 super-resolution for two-stream instability.}}
    \label{fig:vlasov_phase1_two_stream_full}
\end{figure*}







\begin{figure*}[htbp]
    \centering

    \begin{minipage}[b]{\linewidth}
        \centering
        \includegraphics[width=\linewidth]{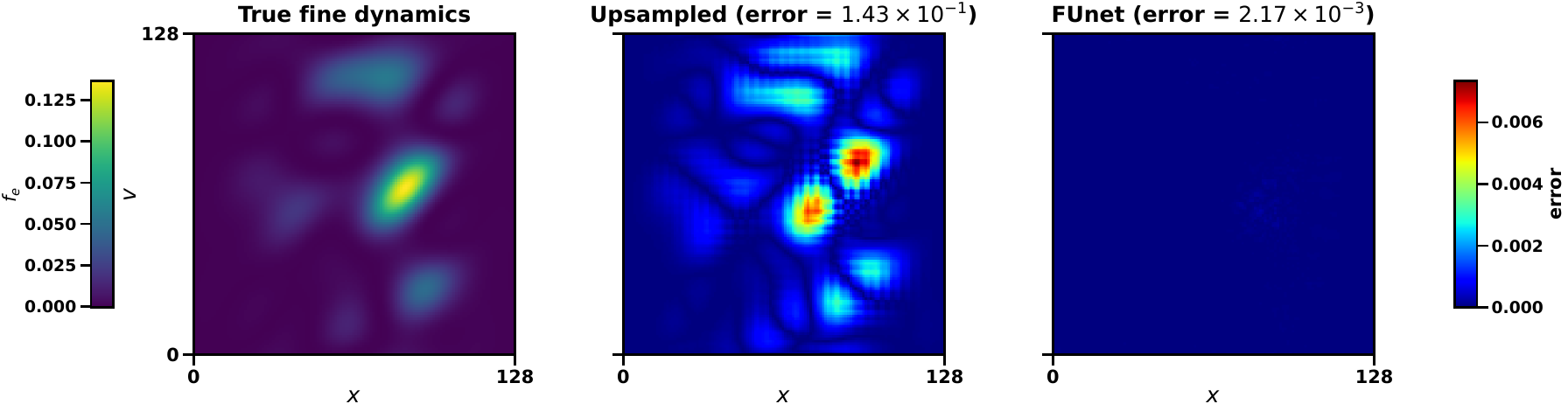}
    \end{minipage}

    \vspace{1em}
    \begin{minipage}[b]{\linewidth}
        \centering
        \includegraphics[width=\linewidth]{figures/2d_vlasov_test_results/randomized_step_050.pdf}
    \end{minipage}

    \vspace{1em}
        \begin{minipage}[b]{\linewidth}
        \centering
        \includegraphics[width=\linewidth]{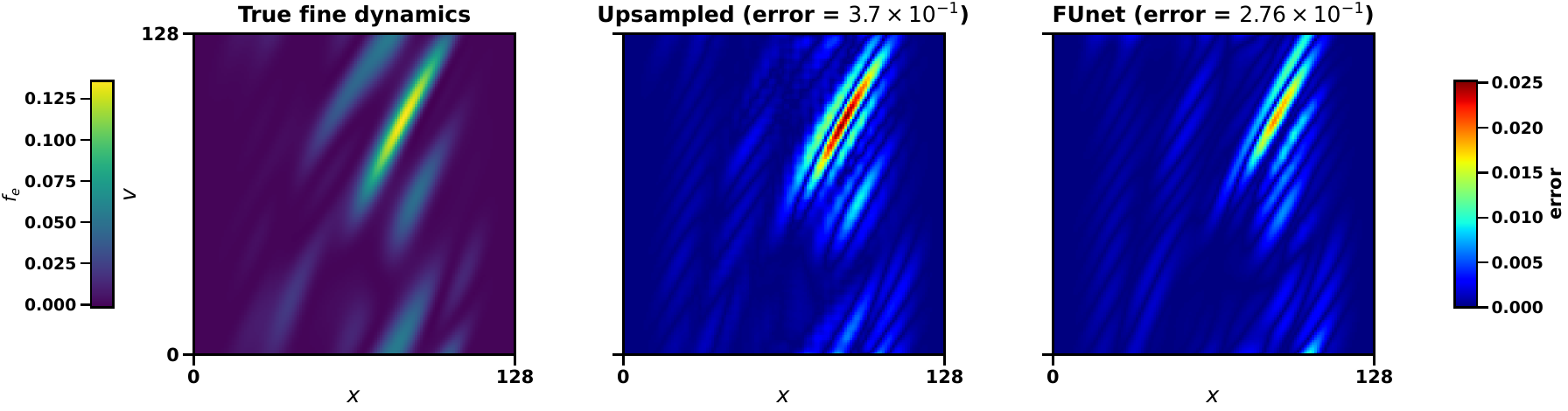}
    \end{minipage}

    \vspace{1em}
        \begin{minipage}[b]{\linewidth}
        \centering
        \includegraphics[width=\linewidth]{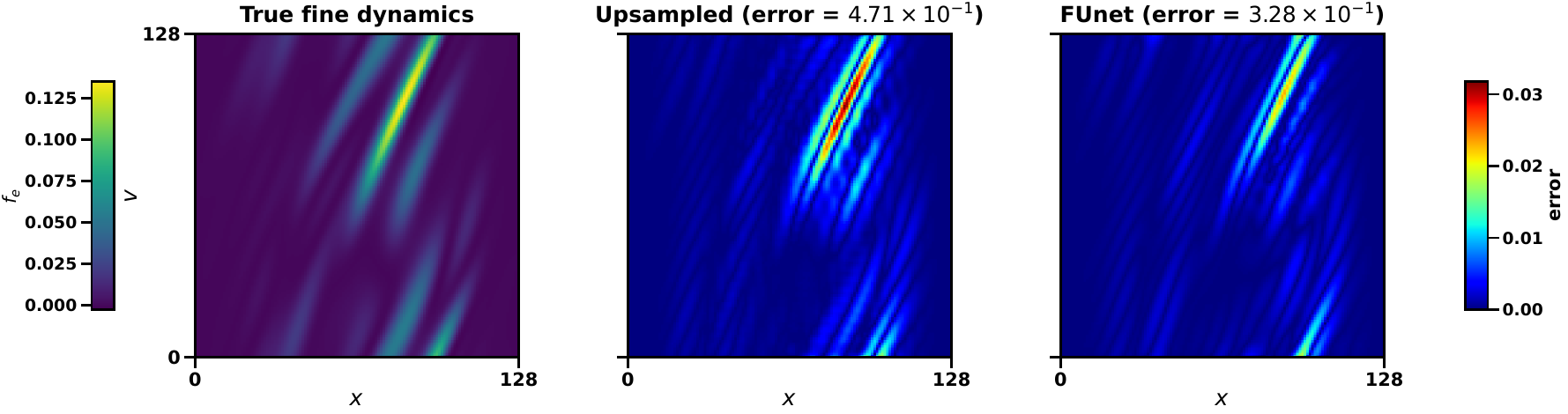}
    \end{minipage}

    \vspace{1em}
    \caption{\justifying \textbf{Total Phase~2 extrapolation for randomized dataset.}}
    \label{fig:vlasov_phase2_randomized_full}
\end{figure*}

\begin{figure*}[htbp]
    \centering

    \begin{minipage}[b]{\linewidth}
        \centering
        \includegraphics[width=\linewidth]{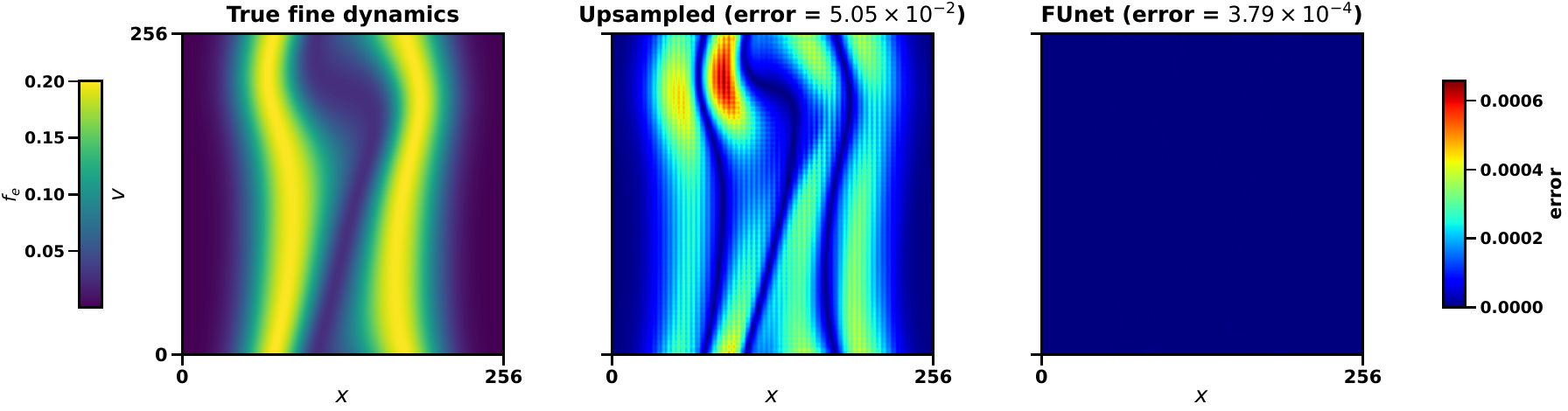}
    \end{minipage}

    \vspace{1em}
    \begin{minipage}[b]{\linewidth}
        \centering
        \includegraphics[width=\linewidth]{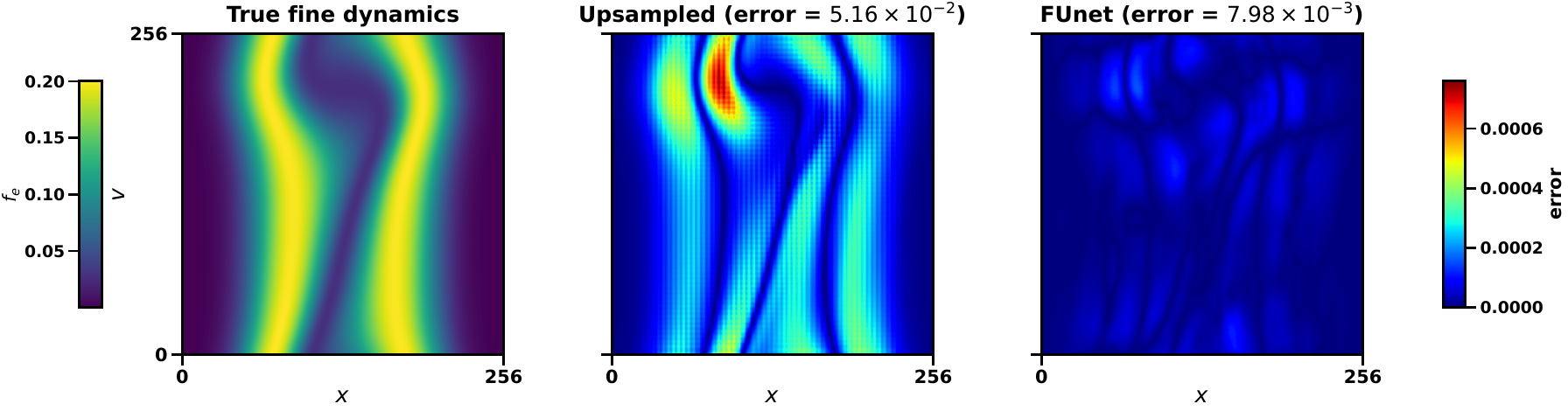}
    \end{minipage}

    \vspace{1em}
        \begin{minipage}[b]{\linewidth}
        \centering
        \includegraphics[width=\linewidth]{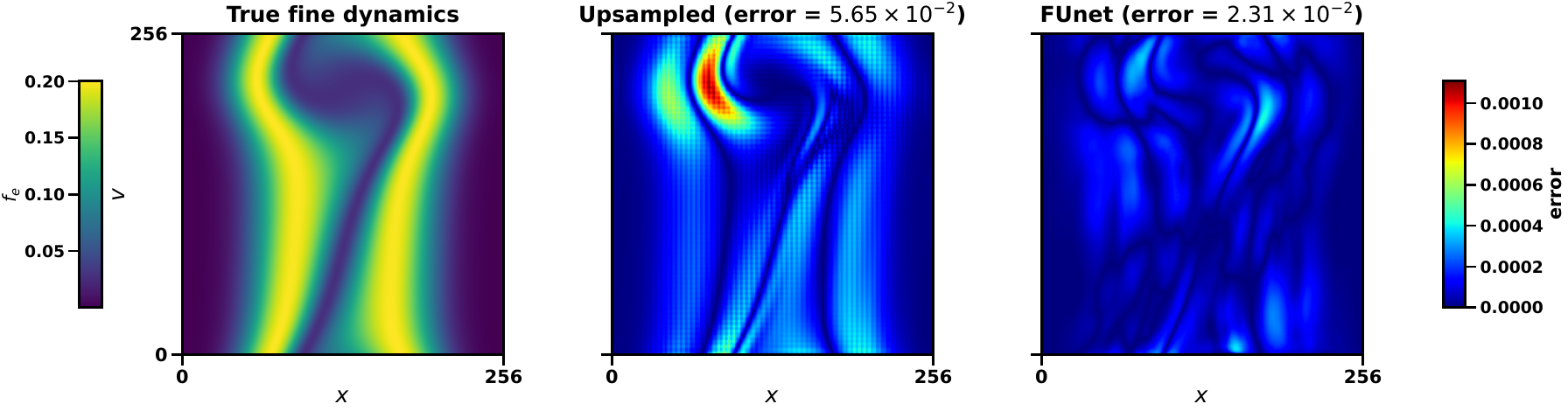}
    \end{minipage}

    \vspace{1em}
        \begin{minipage}[b]{\linewidth}
        \centering
        \includegraphics[width=\linewidth]{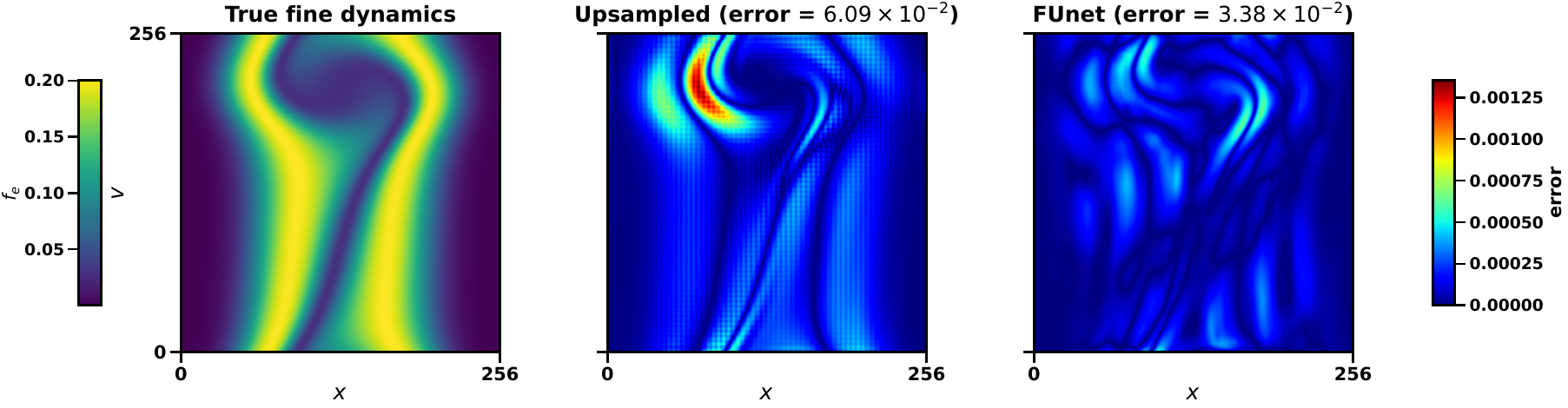}
    \end{minipage}

    \vspace{1em}
    \caption{\justifying  \textbf{Total Phase~2 extrapolation for two-stream instability.}}
    \label{fig:vlasov_phase2_two_stream_full}
\end{figure*}

\begin{figure*}[htbp]
    \centering

    \begin{minipage}[b]{\linewidth}
        \centering
        \includegraphics[width=\linewidth]{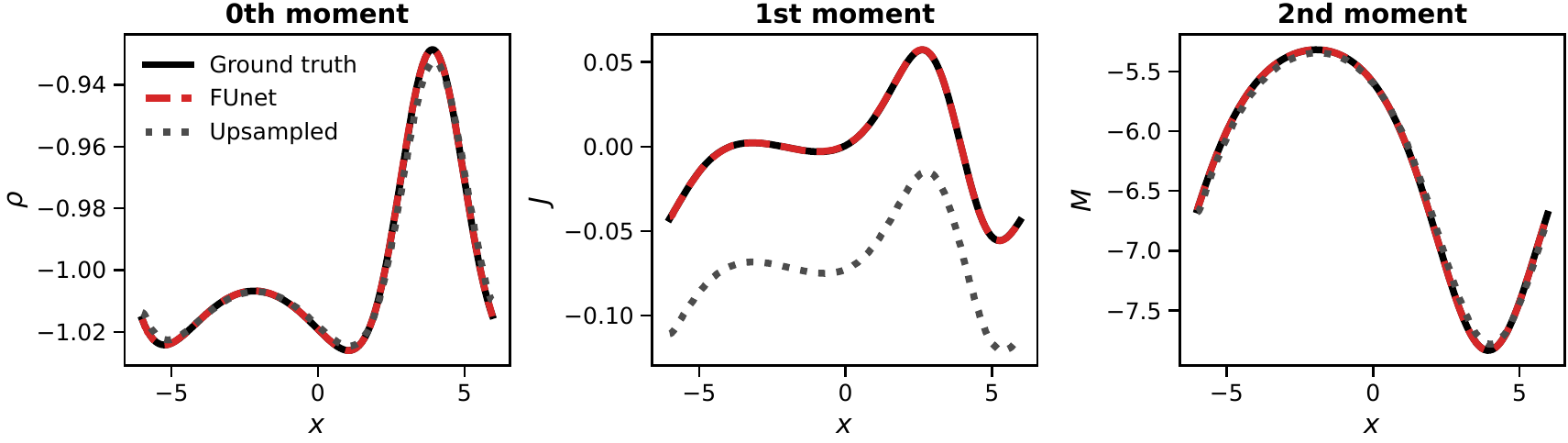}
    \end{minipage}

    \vspace{1em}
    \begin{minipage}[b]{\linewidth}
        \centering
        \includegraphics[width=\linewidth]{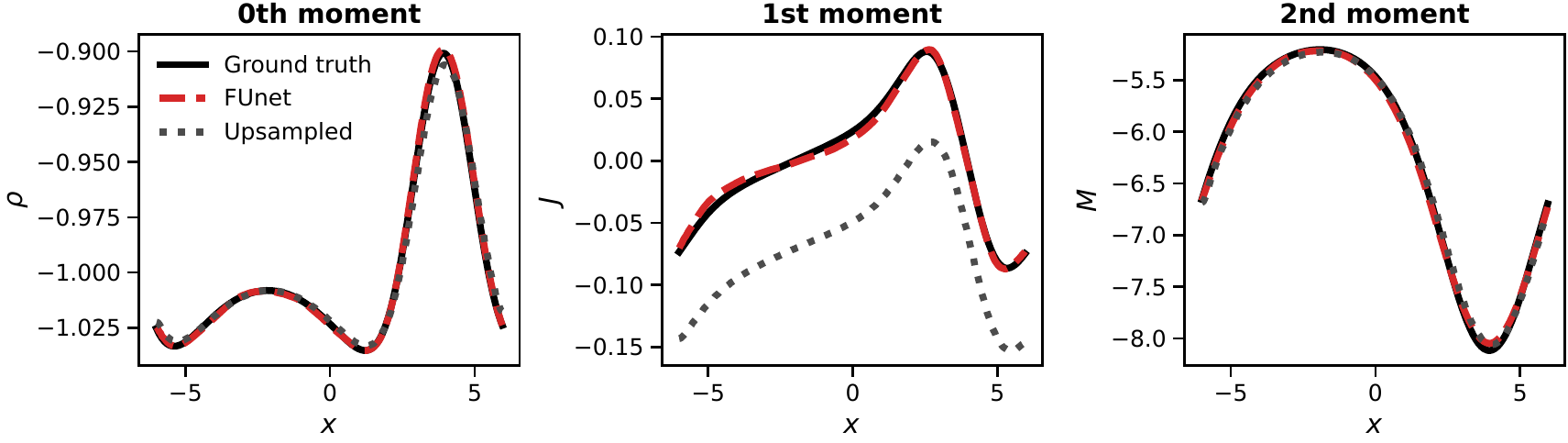}
    \end{minipage}

    \vspace{1em}
    \begin{minipage}[b]{\linewidth}
        \centering
        \includegraphics[width=\linewidth]{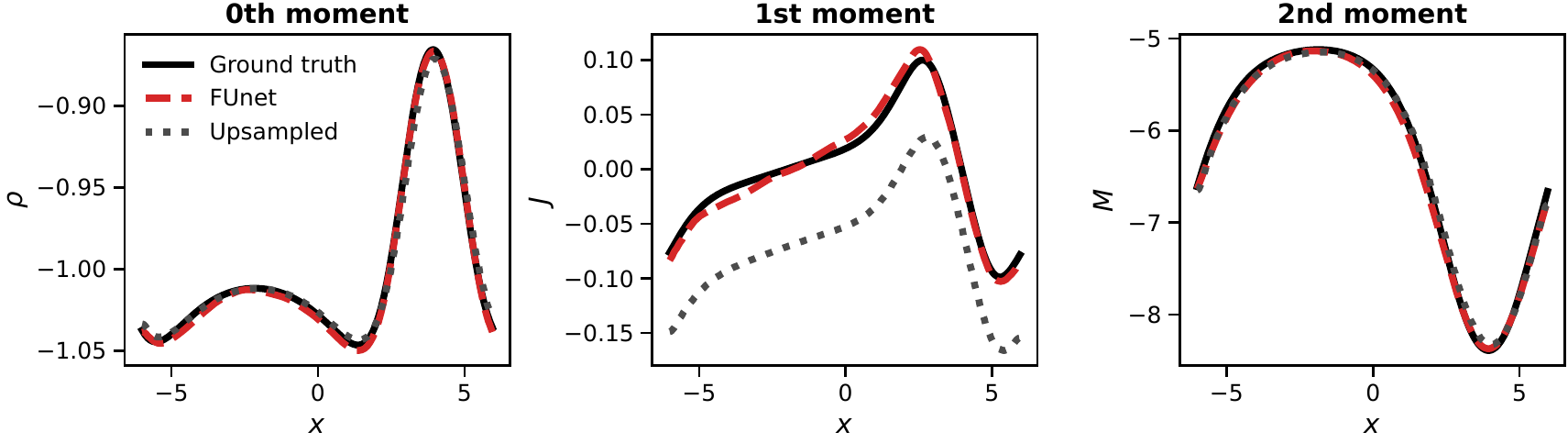}
    \end{minipage}

    \vspace{1em}
        \begin{minipage}[b]{\linewidth}
        \centering
        \includegraphics[width=\linewidth]{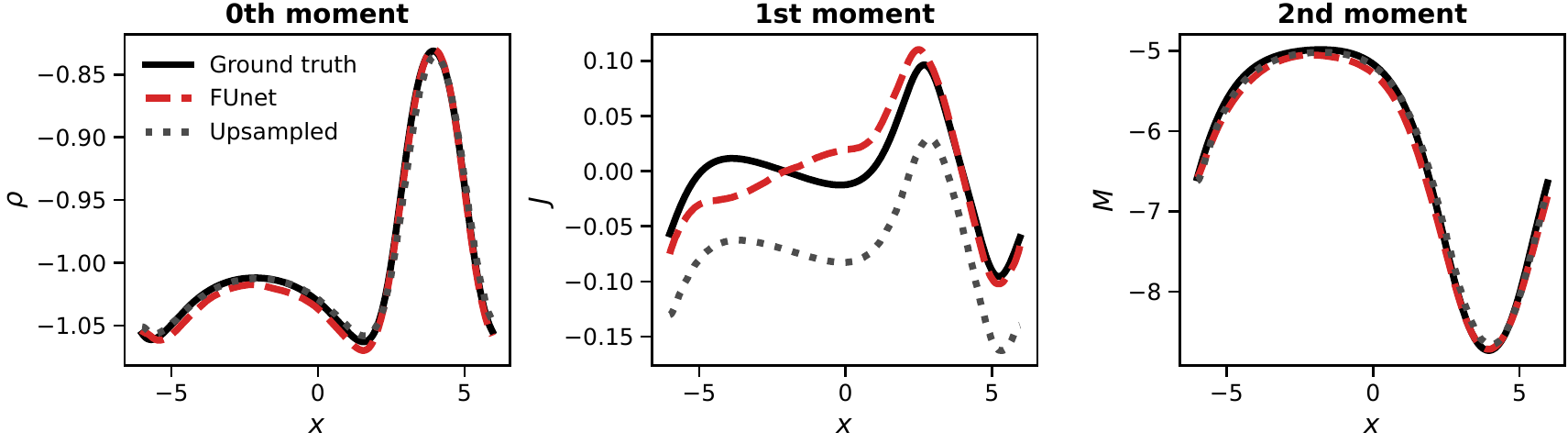}
    \end{minipage}

    \vspace{1em}
    \caption{\justifying  \textbf{Total Phase~2 moments extrapolation for two-stream instability.}}
    \label{fig:vlasov_phase2_two_stream_full_mom}
\end{figure*}

\begin{figure*}[htbp]
    \centering

    \vspace{1em}
    \begin{minipage}[b]{\linewidth}
        \centering
        \includegraphics[width=\linewidth]{figures/2d_vlasov_test_results/randomized_step_050_mom.pdf}
    \end{minipage}

    \vspace{1em}
    \begin{minipage}[b]{\linewidth}
        \centering
        \includegraphics[width=\linewidth]{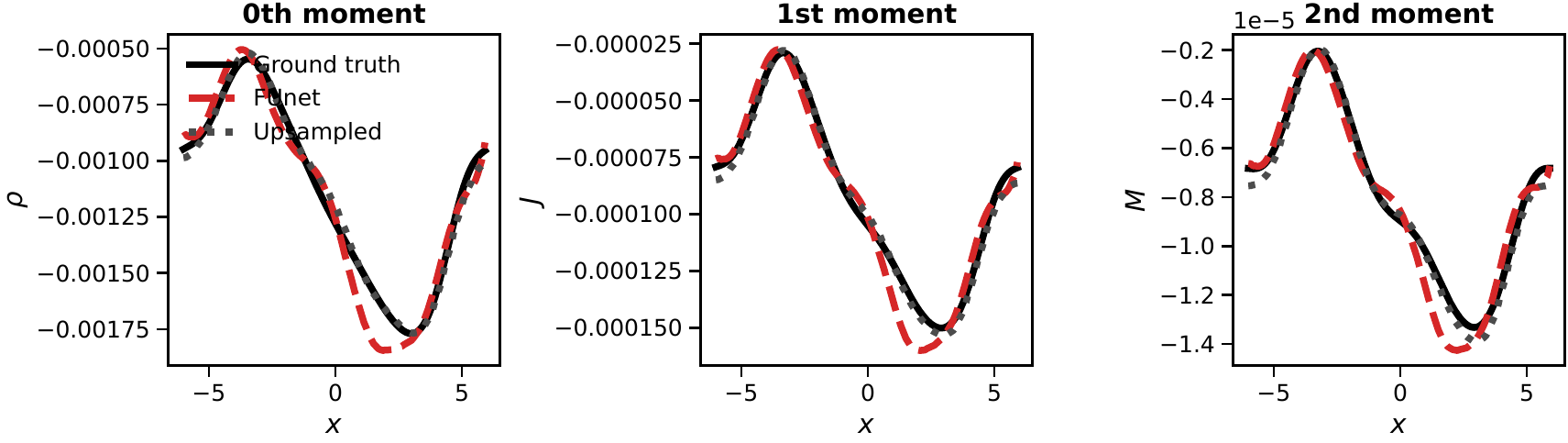}
    \end{minipage}

    \vspace{1em}
        \begin{minipage}[b]{\linewidth}
        \centering
        \includegraphics[width=\linewidth]{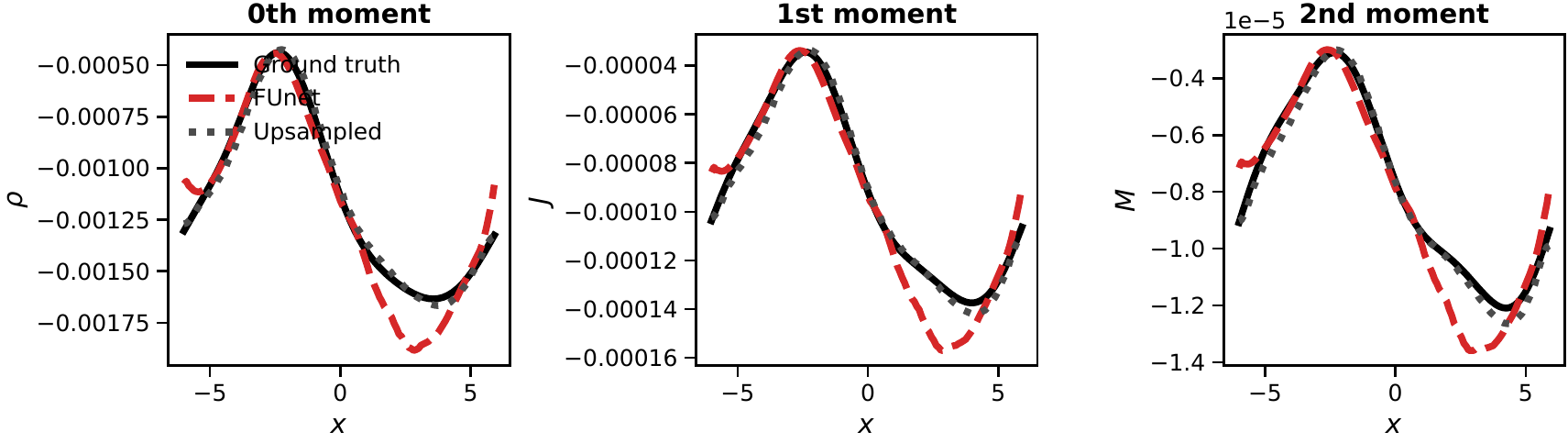}
    \end{minipage}

    \vspace{1em}
        \begin{minipage}[b]{\linewidth}
        \centering
        \includegraphics[width=\linewidth]{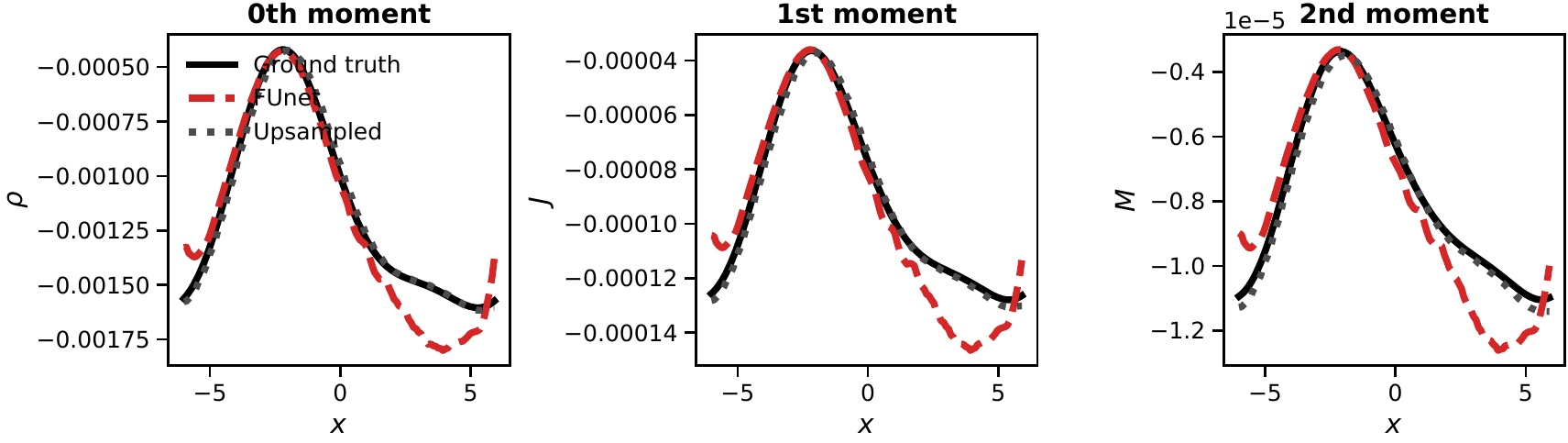}
    \end{minipage}

    \vspace{1em}
    \caption{\justifying  \textbf{Total Phase~2 moments extrapolation for randomized dataset.}}
    \label{fig:vlasov_phase2_randomized_full_mom}
\end{figure*}

\end{document}